\pdfoutput=1
\documentclass[11pt]{article}
\usepackage[left=1in,right=1in,top=1in,bottom=1in]{geometry}
\usepackage{times}
\usepackage{expl3}
\usepackage{cite}
\usepackage[table]{xcolor}
\usepackage{multirow}
\usepackage{stackengine} 
\usepackage{hhline}
\usepackage{lipsum}
\usepackage{titlesec}
\usepackage{wrapfig}
\usepackage{enumerate}
\usepackage{epsfig}
\usepackage{amsmath}
\usepackage{tabularx}
\usepackage{array}
\usepackage{booktabs}
\usepackage{enumitem}
\usepackage{bbm}
\usepackage{calc}
\usepackage{graphicx}
\usepackage{amsmath}
\usepackage[title]{appendix}
\usepackage{amssymb}
\usepackage{epstopdf}
\usepackage{boldline}
\usepackage{arydshln}
\usepackage{calligra}
\usepackage{bm}
\usepackage{url}
\usepackage{blindtext}
\usepackage{accents}

\newcommand{\define}{\stackrel{\mbox{\tiny def}}{=}}

\newtheorem{definition}{Definition}
\newtheorem{theorem}{Theorem}
\newtheorem{proposition}{Proposition}
\newtheorem{corollary}{Corollary}
\newtheorem{lemma}{Lemma}

\newtheorem{example}{Example}
\newtheorem{remark}{Remark}
\usepackage{mathtools}
\usepackage{epstopdf}
\usepackage{balance}
\usepackage{thmtools}
\usepackage{thm-restate}
\usepackage{hyperref}
\usepackage{cleveref}
\usepackage[mathscr]{euscript}

\usepackage[ruled,vlined]{algorithm2e}
\include{pythonlisting}
\newcommand{\ostar}{\mathbin{\mathpalette\make@circled\star}}

\makeatletter
\newcommand{\removelatexerror}{\let\@latex@error\@gobble}
\makeatother
\makeatletter
\newcommand*{\rom}[1]{\expandafter\@slowromancap\romannumeral #1@}
\makeatother

\ExplSyntaxOn
\newcommand\latinabbrev[1]{
  \peek_meaning:NTF . {
    #1\@}%
  { \peek_catcode:NTF a {
      #1.\@ }%
    {#1.\@}}}
\ExplSyntaxOff

\titleclass{\subsubsubsection}{straight}[\subsubsection]

\begin{document}
\vspace{1cm}
\title{Generalized Multivariate Hypercomplex Function Inequalities and Their Applications}
\vspace{1.8cm}
\author{Shih-Yu~Chang
\thanks{Shih Yu Chang is with the Department of Applied Data Science,
San Jose State University, San Jose, CA, U. S. A. (e-mail: {\tt
shihyu.chang@sjsu.edu}). 
           }}

\maketitle

\begin{abstract}
This work extends the Mond-Pecaric method to functions with multiple operators as arguments by providing arbitrarily close approximations of the original functions. Instead of using linear functions to establish lower and upper bounds for multivariate functions as in prior work, we apply sigmoid functions to achieve these bounds with any specified error threshold based on the multivariate function approximation method proposed by Cybenko. This approach allows us to derive fundamental inequalities for multivariate hypercomplex functions, leading to new inequalities based on ratio and difference kinds. For applications about these new derived inequalities for multivariate hypercomplex functions, we first introduce a new concept called W-boundedness for hypercomplex functions by applying ratio kind multivariate hypercomplex inequalities. W-boundedness generalizes R-boundedness for norm mappings with input from Banach space. Additionally, we develop an approximation theory for multivariate hypercomplex functions and establish bounds algebra, including operator bounds and tail bounds algebra for multivariate random tensors. 
\end{abstract}

\begin{keywords}
Operator inequality, hypercomplex function, sigmoid function, Kantorovich function, Banach space, random tensors. 
\end{keywords}

\section{Introduction}\label{sec: Introduction}

Operator inequalities are mathematical statements that establish relationships between operators in terms of ordering relations among operators (inequalities). Common ordering relations among operators include: Lowner ordering, star ordering, majorization ordering, etc~\cite{chang2024generalizedCDJ}. These operators can be linear, bounded, self-adjoint, or positive, and they often act on elements within a Hilbert or Banach space. Operator inequalities generalize many classical inequalities to the setting of operators, offering a powerful framework for analysis in the following fields: optimization, control theory, functional analysis and quantum mechanics~\cite{hutson2005applications,naylor1982linear,aron2022operator}.

Mond and Pecaric introduced the Mond-Pecaric (MP) method to address operator inequalities, requiring a real-valued function to be convex or concave and a normalized positive linear map between Hilbert spaces~\cite{furuta2005mond,fujii2012recentMP}. The author in~\cite{chang2024generalizedJensen} extends the MP method to non-convex/non-concave functions with single input variable and nonlinear mappings. Using the Stone–Weierstrass theorem and Kantorovich function, the MP method is generalized. Examples demonstrate traditional MP method inequalities with convex functions and normalized positive linear maps. New inequalities for hypercomplex functions are derived and applied to approximate these functions using ratio and difference criteria. Additionally, these inequalities establish bounds for hypercomplex functions algebra and derive tail bounds for random tensor ensembles. Note that all hypercomplex functions discussed in~\cite{chang2024generalizedJensen} are assumed with a single-input argument. 

In~\cite{chang2024multivariateMP}, we extend the Mond-Pecaric method from single-input to multiple-input operators. We begin by defining normalized positive linear maps with examples. Using the Mond-Pecaric method, we derive fundamental inequalities for multivariate hypercomplex functions bounded by linear functions. These serve as the basis for several inequalities focusing on ratio relationships and similar results based on difference relationships. Finally, we apply these inequalities to establish Sobolev embedding via the Sobolev inequality for hypercomplex functions with operator inputs.

This work is the extension of the Mond-Pecaric method to functions with multiple input operators by providing arbitrary close approximations of the original functions. Instead of utilizing linear functions to provide lower and upper bounds for multivariate functions in~\cite{chang2024multivariateMP}, in this work, we apply sigmoid functions to give lower and upper bounds for multivariate functions with any specific error difference parameter $\epsilon$ according to the multivariate function approximation method proposed by~\cite{cybenko1989approximation}. Such multivariate function approximation enables us to derive fundamental inequalities for multivariate hypercomplex functions. These fundamental inequalities establish new inequalities of hypercomplex functions based on ratio and difference kinds. A new concept called $W$-boundedness for hypercomplex functions is introduced by applying different kinds of multivariate hypercomplex inequalities. $W$-boundedness is a generalization of $R$-boundedness for real-valued mapping for inputs from Banach space. The second application of derived multivaraite hypercomplex function inequalties is the approximation theory for multivariate hypercomplex functions from ratio and difference comparisons. The third application of multivariate hypercomplex function inequalties is the development of bounds algebra, including operator bounds and tail bounds algebra for multivariate random tensors~\cite{gozlan2019high,chang2022convenient,chang2022randomPDT,chang2022randomMOI,chang2022randomDTI,chang2022generalMaj,chang2023tailTRP,chang2023BiRanTenPartI,chang2023tailMulRanTenMeans,chang2023algebraicConn}.

The rest of this paper is structured by the following sections. We show that a linear combination of sigmoid functions can be used to approximate any continuous multivariate functions in Section~\ref{sec: Multivariate Hypercomplex Functions Dense Representation By Self-Adjoint Hypercomplex Functions}. In Section~\ref{sec: Fundamental Inequalities for Multivariate Hypercomplex Functions with Nonlinear Map}, fundamental inequalities for multivariate hypercomplex functions using sigmoid functions are derived. In Section~\ref{sec: Multivariate Hypercomplex Functions Inequalities: Ratio Kind}, multivariate hypercomplex functions inequalities based on ratio kind are derived. In Section~\ref{sec: Multivariate Hypercomplex Functions Inequalities: Difference Kind}, multivariate hypercomplex functions inequalities based on difference kind are established. Following sections are applications of derived multivariate hypercomplex functions inequalities from Section~\ref{sec: Fundamental Inequalities for Multivariate Hypercomplex Functions with Nonlinear Map} to Section~\ref{sec: Multivariate Hypercomplex Functions Inequalities: Difference Kind}. New concept, named as $W$-boundedness for hypercomplex functions, is introduced in Section~\ref{sec: W-boundedness with Loewner Order}. In Section~\ref{sec: Multivariate Hypercomplex Function Approximations}, the approximation theory about multivariate hypercomplex functions is presented. Finally, bounds algebra for multivariate function of operator bounds and tail bounds for hypercomplex functions with multivariate random tensors are discussed in Section~\ref{sec: Bounds Algebra}.

\textbf{Nomenclature:} 
Inequalities \( \geq, >, \leq, \) and \( < \) for operators follow the Loewner ordering. The symbol \( \bm{0} \) represents the zero operator. \( \Lambda(\bm{A}) \) denotes the spectrum of the operator \( \bm{A} \), i.e., the set of eigenvalues of \( \bm{A} \). If \( \Lambda(\bm{A}) \) consists of real numbers, \( \min(\Lambda(\bm{A})) \) and \( \max(\Lambda(\bm{A})) \) represent the minimum and maximum values within \( \Lambda(\bm{A}) \), respectively. For given \( M > m > 0 \) and any \( r \in \mathbb{R} \) where \( r \neq 1 \), the Kantorovich function with respect to $m$, $M$, and $r$ is defined as follows:
\begin{eqnarray}\label{eq: Kantorovich function}
\mathscr{K}(m,M,r)&=&\frac{(mM^r - Mm^r)}{(r-1)(M-m)}\left[\frac{(r-1)(M^r-m^r)}{r(mM^r - Mm^r)}\right]^r.
\end{eqnarray}
A vector of variables $x_1,x_2,\ldots,x_n$ will be denoted by $\bm{x}$. A vector of indices $k_1,\ldots,k_n$ will be denoted by $\bm{k}$. A vector of indices $j_1,\ldots,j_n$ will be denoted by $\bm{j}$. We use $\bm{j}=\bm{1}$ to represent $j_1=1, j_2=1,\ldots,j_n=1$. A vector of operators $\bm{X}_1,\bm{X}_2,\ldots,\bm{X}_n$ will be denoted by $\underline{\bm{X}}$. A vector of operators $\bm{A}_1,\bm{X}_2,\ldots,\bm{X}_n$ will be denoted by $\underline{\bm{A}}$. A vector of operators $\bm{A}_{j_1},\bm{A}_{j_2},\ldots,\bm{A}_{j_n}$ will be represented by $\underline{\bm{A}}_{\bm{j}}$. The product of weights (nonnegative real number between 0 and 1)$w_{j_1}\cdot w_{j_2} \cdots  w_{j_n}$ will be represented by $\bm{w}^{\times}_{\bm{j}}$. Given $n$ real intervals $I_1, I_2, \ldots, I_n$, the Cartesian product of these $n$ intervals $\bigtimes\limits_{i=1}^n I_i$ will be represented by $\bm{I}^{\bigtimes}$~\cite{floudas2008encyclopedia,anastassiou2019frontiers}.

\section{Multivariate Hypercomplex Functions Dense Representation By Self-Adjoint Hypercomplex Functions}\label{sec: Multivariate Hypercomplex Functions Dense Representation By Self-Adjoint Hypercomplex Functions}

Given a multivariate operator-valued function $f(\underline{\bm{X}})$ with $n$ self-adjoint operators $\bm{X}_1, \bm{X}_2, \ldots,\bm{X}_n$ as input operators, the multivaraite function $f$ will be called as a \emph{self-adjoint hypercomplex function} if $f(\underline{\bm{X}})$ is a self-adjoint operator. As the product of two self-adjoint operators may not be a self-adjoint operator, we will provide the follownig Theorem~\ref{thm: Represent Multivariate Hypercomplex Function} to approximate any multivariate operator-valued function by a self-adjoint hypercomplex function.

\begin{theorem}\label{thm: Represent Multivariate Hypercomplex Function}
Given a multivariate operator-valued function $f(\underline{\bm{X}})$ with $n$ self-adjoint operators $\bm{X}_1, \bm{X}_2, \ldots,\bm{X}_n$ as input operators with their spectrum $\Lambda(\bm{X}_1) \times \Lambda(\bm{X}_2) \times \ldots \Lambda(\bm{X}_n) \subset \bm{I}^{\bigtimes}=\bigtimes\limits_{i=1}^n I_i$, where $I_i$ are real intervals, we can approxiate this operator-valued function $f(\underline{\bm{X}})$ by the following operator-valued function $\Psi(\underline{\bm{X}})$ with the format as:
\begin{eqnarray}\label{eq1: thm: Represent Multivariate Hypercomplex Function}
\Psi(\underline{\bm{X}})&=&\sum\limits_{i=1}^N A_i \sigma\left(\sum\limits_{j=1}^n C_j \bm{X}_j  + B_i \bm{I}\right), 
\end{eqnarray}
where $A_i, B_i, C_j$ are real scalars, and $\sigma(\bm{Y})$ is a sigmoid operator-valued function expressed by
\begin{eqnarray}\label{eq2: thm: Represent Multivariate Hypercomplex Function}
\sigma(\bm{Y})=\left(\bm{I}+\exp(-\bm{Y})\right)^{-1};
\end{eqnarray}
such that 
\begin{eqnarray}\label{eq3: thm: Represent Multivariate Hypercomplex Function}
|\Psi(\bm{x}) - f(\bm{x})| < \epsilon,
\end{eqnarray}
for any $\epsilon > 0$ and $\bm{x} \in \bm{I}^{\bigtimes}$, where 
\begin{eqnarray}\label{eq3.1: thm: Represent Multivariate Hypercomplex Function}
\Psi(\bm{x})=\sum\limits_{i=1}^N A_i \sigma\left(\sum\limits_{j=1}^n C_j x_j  + B_i\right).
\end{eqnarray}
Moreover, $\Psi(\underline{\bm{X}})$ is a \emph{self-adjoint hypercomplex function}.
\end{theorem}
\textbf{Proof:}
Our proof is based on a contradiction argument. 

Let $\mathfrak{C}\left(\bm{I}^{\bigtimes}\right)$ be a set of continuous functions defined on the domain $\bm{I}^{\bigtimes}$ and let $S_\Psi\left(\bm{I}^{\bigtimes}\right)$ be a set of continuous functions with the format expressed by Eq.~\eqref{eq3.1: thm: Represent Multivariate Hypercomplex Function}. The closure of the space $S_\Psi\left(\bm{I}^{\bigtimes}\right)$ is denoted by $\overline{S_\Psi\left(\bm{I}^{\bigtimes}\right)}$. The goal is to prove
\begin{eqnarray}\label{eq4: thm: Represent Multivariate Hypercomplex Function}
\overline{S_\Psi\left(\bm{I}^{\bigtimes}\right)}&=&\mathfrak{C}\left(\bm{I}^{\bigtimes}\right).
\end{eqnarray}
If we assume that Eq.~\eqref{eq4: thm: Represent Multivariate Hypercomplex Function} is not true, then we have $\overline{S_\Psi\left(\bm{I}^{\bigtimes}\right)} \subsetneqq \mathfrak{C}\left(\bm{I}^{\bigtimes}\right)$, i.e., it is a closed but proper subspace of $\mathfrak{C}\left(\bm{I}^{\bigtimes}\right)$. From the Hahn-Banach theorem, there is a bounded
linear functional on $\mathfrak{C}\left(\bm{I}^{\bigtimes}\right)$, named as $L$, such that $L\left(S_\Psi\left(\bm{I}^{\bigtimes}\right)\right)=0$ and $L\left(\overline{S_\Psi\left(\bm{I}^{\bigtimes}\right)}\right)=0$ with $L \neq 0$. 

Let $\mathfrak{B}\left(\bm{I}^{\bigtimes}\right)$ be a finite, signed regular Borel measures on the domain $\bm{I}^{\bigtimes}$, from Riesz representation theorem, the bounded linear functional $L$ can be expressed as
\begin{eqnarray}\label{eq5: thm: Represent Multivariate Hypercomplex Function}
L(f)&=&\int_{\bm{I}^{\bigtimes}} f(\bm{x}) d \mu(\bm{x}),
\end{eqnarray}
for some $\mu \in \mathfrak{B}\left(\bm{I}^{\bigtimes}\right)$ and all $f(\bm{x}) \in \mathfrak{C}\left(\bm{I}^{\bigtimes}\right)$. Because the function $\sigma\left(\sum\limits_{j=1}^n C_j x_j  + B_i\right)$ is in $\overline{S_\Psi\left(\bm{I}^{\bigtimes}\right)}$ with respect to all $C_j$ and $B_i$, we have
\begin{eqnarray}\label{eq6: thm: Represent Multivariate Hypercomplex Function}
\int_{\bm{I}^{\bigtimes}}\sigma\left(\sum\limits_{j=1}^n C_j x_j  + B_i\right)d \mu(\bm{x})&=&0,
\end{eqnarray}
for all $C_j$ and $B_i$. From Lemma 1 in~\cite{cybenko1989approximation}, we know that the function $\sigma$ is a discriminatory function. Then, we must have $\mu(\bm{x})=0$ from Eq.~\eqref{eq6: thm: Represent Multivariate Hypercomplex Function}. However, we know that $L \neq 0$ and a contradiction is introduced. Therefore, Eq.~\eqref{eq4: thm: Represent Multivariate Hypercomplex Function} is true.

For the part about self-adjoint hypercomplex function of $\Psi(\underline{\bm{X}})$, this statement is true due to that the linear combination, $\exp$ and the inversion operations of a self-adjoint operator is still self-adjoint.
$\hfill \Box$

Our proof of Theorem~\ref{thm: Represent Multivariate Hypercomplex Function} is modified from proofs in~\cite{cybenko1989approximation} by considering hypercomplex  functions instead of conventional functions with real numbers as input arguments. 

\section{Fundamental Inequalities for Multivariate Hypercomplex Functions with Nonlinear Map}\label{sec: Fundamental Inequalities for Multivariate Hypercomplex Functions with Nonlinear Map}

\subsection{Upper Bound and Lower Bound of Multivaraite Function $f(x_1,x_2,\ldots,x_n)$}\label{sec: Upper Bound and Lower Bound of Function f}

From Theorem~\ref{thm: Represent Multivariate Hypercomplex Function}, given any $\epsilon>0$, we have a self-adjoint hypercomplex function $\Phi_{\mathscr{U}}(\bm{x})$ (as upper bound) and a self-adjoint hypercomplex function $\Phi_{\mathscr{L}}(\bm{x})$ (as lower bound) satisfying:
\begin{eqnarray}\label{eq:lower and upper Phi}
0 &\leq& \Psi_{\mathscr{U}}(\bm{x}) - f(\bm{x})~~\leq~~\epsilon, \nonumber \\
0 &\leq& f(\bm{x})-\Psi_{\mathscr{L}}(\bm{x})~~\leq~~\epsilon,
\end{eqnarray}
for $\bm{x} \in \bm{I}^{\bigtimes}$.

\subsection{Polynomial Map $\Phi$}\label{sec: Polynomial Maps Phi}

Consider two Hilbert spaces, $\mathfrak{H}$ and $\mathfrak{K}$. The sets $\mathbb{B}(\mathfrak{H})$ and $\mathbb{B}(\mathfrak{K})$ denote the algebras of all bounded linear operators on these respective spaces. In the context of the Choi-Davis-Jensen inequality, the mapping $\Phi: \mathbb{B}(\mathfrak{H}) \rightarrow \mathbb{B}(\mathfrak{K})$ is defined as a normalized positive linear map. This map is precisely characterized by Definition~\ref{def: normalized positive linear map} below~\cite{furuta2005mond,fujii2012recentMP}.
\begin{definition}\label{def: normalized positive linear map}
A map $\Phi: \mathbb{B}(\mathfrak{H}) \rightarrow \mathbb{B}(\mathfrak{K})$ is considered a normalized positive linear map if it satisfies the following conditions:
\begin{itemize}
\item Linearity: $\Phi(a\bm{X}+b\bm{Y})=a\Phi(\bm{X})+b\Phi(\bm{Y})$ for any $a,b \in \mathbb{C}$ and any $\bm{X}, \bm{Y} \in \mathbb{B}(\mathfrak{H})$.
\item Positivity: If $\bm{X}\geq\bm{Y}$, then $\Phi(\bm{X})\geq\Phi(\bm{Y})$.
\item Normalization: $\Phi(\bm{I}_{\mathfrak{H}})=\bm{I}_{\mathfrak{K}}$, where $\bm{I}_{\mathfrak{H}}$ and $\bm{I}_{\mathfrak{K}}$ are the identity operators of the Hilbert spaces $\mathfrak{H}$ and $\mathfrak{K}$, respectively. 
\end{itemize}
\end{definition}

In this paper, we will consider a broader class of $\Phi$ by defining $\Phi$ as follows: 
\begin{eqnarray}\label{eq: new phi def}
\Phi(\bm{X}) &=& \bm{V}^{\ast}\left(\sum\limits_{i=0}^{I}a_{i}\bm{X}^{i}\right)\bm{V}\nonumber \\
&=& \bm{V}^{\ast}\left(\sum\limits_{i_{+}\in S_{I_{+}}}a_{i_{+}}\bm{X}^{i_{+}}+\sum\limits_{i_{-}\in S_{I_{-}}}a_{i_{-}}\bm{X}^{i_{-}}\right)\bm{V},
\end{eqnarray}
where, $\bm{V}$ stands as an isometry within $\mathfrak{H}$, adhering to $\bm{V}^{\ast}\bm{V}=\bm{I}_{\mathfrak{H}}$. The definition about $\Phi(\bm{X})$ is also used in~\cite{chang2024generalizedJensen} but we present here again for a self-contained presentation purpose. In $a_i$, the coefficients $a_{i_{+}}$ are nonnegative, while $a_{i_{-}}$ are negative. The indices of the positive coefficients form the set $S_{I_{+}}$, and those of the negative coefficients form $S_{I_{-}}$. Notably, no constraints on linearity, positivity, or normalization are imposed on $\Phi$ as specified in Eq.~\eqref{eq: new phi def}. Under this premise, the standard notion of a normalized positive linear map, described in Definition~\ref{def: normalized positive linear map}, becomes a special case. This is achieved by setting the polynomial $\sum_{i=0}^{I}a_i\bm{X}^i$ to be the identity map, with all coefficients $a_i$ being zero except for $a_1$.

We require the following Lemma~\ref{lma: lower and upper bound for f(A)} to provide the lower and upper bounds for $f(\underline{\bm{A}})$ if $\underline{\bm{A}}=\bm{A}_1,\bm{A}_2,\ldots,\bm{A}_n$ are $n$ self-adjoint operators. 
\begin{lemma}\label{lma: lower and upper bound for f(A)}
Given $n$ self-adjoint operators $\bm{A}_1,\bm{A}_2,\ldots,\bm{A}_n$ with their spectrums in $\bm{I}^{\bigtimes}$, such that 
\begin{eqnarray}\label{eq1: lma: lower and upper bound for f(A)}
0 &\leq& \Psi_{\mathscr{U}}(\bm{x}) - f(\bm{x})~~\leq~~\epsilon, \nonumber \\
0 &\leq& f(\bm{x})-\Psi_{\mathscr{L}}(\bm{x})~~\leq~~\epsilon,
\end{eqnarray}
for $\bm{x} \in \bm{I}^{\bigtimes}$ with the linear combination of sigmoidal functions $\Phi_{\mathscr{U}}(\bm{x})$ and $\Phi_{\mathscr{L}}(\bm{x})$ expressed by 
\begin{eqnarray}\label{eq: lower poly formats}
\Psi_{\mathscr{L}}(\bm{x})&=&\sum\limits_{i=1}^{N_{\mathscr{L}}} A_{\mathscr{L},i} \sigma\left(\sum\limits_{j=1}^n C_{\mathscr{L},j} x_j  + B_{\mathscr{L},i}\right),\nonumber \\
\Psi_{\mathscr{U}}(\bm{x})&=&\sum\limits_{i=1}^{N_{\mathscr{U}}} A_{\mathscr{U},i} \sigma\left(\sum\limits_{j=1}^n C_{\mathscr{U},j} x_j  + B_{\mathscr{U},i}\right). 
\end{eqnarray}

We use $\min(\Lambda(\Psi_{\mathscr{L}}(\underline{\bm{A}})))$ and $\max(\Lambda(\Psi_{\mathscr{L}}(\underline{\bm{A}})))$ to represent the minimum and the maximum values for eigenvalues of the operator $\Psi_{\mathscr{L}}(\underline{\bm{A}})$.  Similarly, we also use $\min(\Lambda(\Psi_{\mathscr{U}}(\underline{\bm{A}})))$ and $\max(\Lambda(\Psi_{\mathscr{U}}(\underline{\bm{A}})))$ to represent the minimum and the maximum values for eigenvalues of the operator $\Psi_{\mathscr{U}}(\underline{\bm{A}})$. 

By further assuming that $\Psi_{\mathscr{L}}(\underline{\bm{A}})\geq \bm{0}$, we have
\begin{eqnarray}\label{eq2-1: lma: lower and upper bound for f(A)}
 f^{i_+}(\underline{\bm{A}})&\geq&\mathscr{K}^{-1}(\min(\Lambda(\Psi_{\mathscr{L}}(\underline{\bm{A}}))),\max(\Lambda(\Psi_{\mathscr{L}}(\underline{\bm{A}}))),i_+)\Psi^{i_+}_{\mathscr{L}}(\underline{\bm{A}}),
\end{eqnarray}
and
\begin{eqnarray}\label{eq2-2: lma: lower and upper bound for f(A)}
f^{i_+}(\underline{\bm{A}})&\leq&\mathscr{K}(\min(\Lambda(\Psi_{\mathscr{U}}(\underline{\bm{A}}))),\max(\Lambda(\Psi_{\mathscr{U}}(\underline{\bm{A}}))),i_+)\Psi^{i_+}_{\mathscr{U}}(\underline{\bm{A}}),
\end{eqnarray}
where Kantorovich functions $\mathscr{K}^{-1}(\min(\Lambda(\Psi_{\mathscr{L}}(\underline{\bm{A}}))),\max(\Lambda(\Psi_{\mathscr{L}}(\underline{\bm{A}}))),i_+)$ and \\$\mathscr{K}(\min(\Lambda(\Psi_{\mathscr{U}}(\underline{\bm{A}}))),\max(\Lambda(\Psi_{\mathscr{U}}(\underline{\bm{A}}))),i_+)$ are defined by Eq.~\eqref{eq: Kantorovich function}. 

Moreover, we also have 
\begin{eqnarray}\label{eq3-1: lma: lower and upper bound for f(A)}
f^{i_-}(\underline{\bm{A}})&\leq&\mathscr{K}(\min(\Lambda(\Psi_{\mathscr{U}}(\underline{\bm{A}}))),\max(\Lambda(\Psi_{\mathscr{U}}(\underline{\bm{A}}))),i_-)\Psi^{i_-}_{\mathscr{U}}(\underline{\bm{A}}),
\end{eqnarray}
and
\begin{eqnarray}\label{eq3-2: lma: lower and upper bound for f(A)}
f^{i_-}(\underline{\bm{A}})&\geq&\mathscr{K}^{-1}(\min(\Lambda(\Psi_{\mathscr{L}}(\underline{\bm{A}}))),\max(\Lambda(\Psi_{\mathscr{L}}(\underline{\bm{A}}))),i_-)\Psi^{i_-}_{\mathscr{L}}(\underline{\bm{A}}),
\end{eqnarray}
where Kantorovich function $\mathscr{K}(\min(\Lambda(\Psi_{\mathscr{U}}(\underline{\bm{A}}))),\max(\Lambda(\Psi_{\mathscr{U}}(\underline{\bm{A}}))),i_-)$ and \\$\mathscr{K}^{-1}(\min(\Lambda(\Psi_{\mathscr{L}}(\underline{\bm{A}}))),\max(\Lambda(\Psi_{\mathscr{L}}(\underline{\bm{A}}))),i_-)$ are defined by Eq.~\eqref{eq: Kantorovich function}.
\end{lemma}
\textbf{Proof:}
From Eq.~\eqref{eq1: lma: lower and upper bound for f(A)} and spectrum mapping theorem of the self-adjoint operator $\bm{A}$, we have 
\begin{eqnarray}\label{eq4: lma: lower and upper bound for f(A)}
f(\underline{\bm{A}})&\leq&\Psi_{\mathscr{U}}(\underline{\bm{A}}),
\end{eqnarray}
and 
\begin{eqnarray}\label{eq5: lma: lower and upper bound for f(A)}
\Psi_{\mathscr{L}}(\underline{\bm{A}})&\leq&f(\underline{\bm{A}}).
\end{eqnarray}
From Theorem 8.3 in~\cite{furuta2005mond}, we have Eq.~\eqref{eq2-1: lma: lower and upper bound for f(A)} from Eq.~\eqref{eq5: lma: lower and upper bound for f(A)}. Similarly, we have Eq.~\eqref{eq2-2: lma: lower and upper bound for f(A)} from Eq.~\eqref{eq4: lma: lower and upper bound for f(A)}. 

Again, from Theorem 8.3 in~\cite{furuta2005mond}, we have Eq.~\eqref{eq3-1: lma: lower and upper bound for f(A)} from Eq.~\eqref{eq4: lma: lower and upper bound for f(A)}. Similarly, we have Eq.~\eqref{eq3-2: lma: lower and upper bound for f(A)} from Eq.~\eqref{eq5: lma: lower and upper bound for f(A)}. 
$\hfill \Box$

Our next Lemma~\ref{lma: phi(f(A)) bounds} is about the upper and the lower bounds for $\Phi(f(\underline{\bm{A}}))$ based on Lemma~\ref{lma: lower and upper bound for f(A)}.
\begin{lemma}\label{lma: phi(f(A)) bounds}
Under the definition of $\Phi$ provided by Eq.~\eqref{eq: new phi def} and same conditions provided by  Lemma~\ref{lma: lower and upper bound for f(A)}, we have
\begin{eqnarray}\label{eq1:lma: phi(f(A)) bounds}
\lefteqn{\Phi(f(\underline{\bm{A}}))}\nonumber \\
&\leq& \bm{V}^{\ast}\left\{\sum\limits_{i_{+}\in S_{I_{+}}}a_{i_{+}}\mathscr{K}(\min(\Lambda(\Psi_{\mathscr{U}}(\underline{\bm{A}}))),\max(\Lambda(\Psi_{\mathscr{U}}(\underline{\bm{A}}))),i_+)\Psi^{i_+}_{\mathscr{U}}(\underline{\bm{A}})\right. \nonumber \\
&  &\left.+\sum\limits_{i_{-}\in S_{I_{-}}}a_{i_{-}}\mathscr{K}^{-1}(\min(\Lambda(\Psi_{\mathscr{L}}\underline{\bm{A}}))),\max(\Lambda(\Psi_{\mathscr{L}}(\underline{\bm{A}}))),i_-)\Psi^{i_-}_{\mathscr{L}}(\underline{\bm{A}})\right\}\bm{V} \nonumber \\
&\define&  \bm{V}^{\ast}\mbox{Poly}_{f,\Psi,\mathscr{U}}(\underline{\bm{A}})\bm{V}.
\end{eqnarray}
On the other hand, we also have
\begin{eqnarray}\label{eq2:lma: phi(f(A)) bounds}
\lefteqn{\Phi(f(\underline{\bm{A}}))}\nonumber \\
&\geq& \bm{V}^{\ast}\left\{\sum\limits_{i_{+}\in S_{I_{+}}}a_{i_{+}}\mathscr{K}^{-1}(\min(\Lambda(\Psi_{\mathscr{L}}(\underline{\bm{A}}))),\max(\Lambda(\Psi_{\mathscr{L}}(\underline{\bm{A}}))),i_+)\Psi^{i_+}_{\mathscr{L}}(\underline{\bm{A}}) \right. \nonumber \\
&  &\left.+\sum\limits_{i_{-}\in S_{I_{-}}}a_{i_{-}}\mathscr{K}(\min(\Lambda(\Psi_{\mathscr{U}}(\underline{\bm{A}}))),\max(\Lambda(\Psi_{\mathscr{U}}(\underline{\bm{A}}))),i_-)\Psi^{i_-}_{\mathscr{U}}(\underline{\bm{A}})\right\}\bm{V}\nonumber \\
&\define&  \bm{V}^{\ast}\mbox{Poly}_{f,\Psi,\mathscr{L}}(\underline{\bm{A}})\bm{V}.
\end{eqnarray}
\end{lemma}
\textbf{Proof:}
This Lemma is proved by applying Lemma~\ref{lma: lower and upper bound for f(A)} to the definition of $\Phi$ provided by Eq.~\eqref{eq: new phi def}.
$\hfill \Box$

Given $\Lambda(\bm{A}_1)\times\ldots\times\Lambda(\bm{A}_n)\in \bm{I}^{\bigtimes}$, the eigenvalues range in $\mathbb{R}$ for the operator $\mbox{Poly}_{f,\Psi,\mathscr{L}}(\underline{\bm{A}})$ is represented by $\widetilde{\mbox{Poly}}_{f,\Psi,\mathscr{L}}\left(\bm{I}^{\bigtimes}\right)$. Similarly, the eigenvalues range in $\mathbb{R}$ for the operator $\mbox{Poly}_{f,\Psi,\mathscr{U}}(\underline{\bm{A}})$ is represented by $\widetilde{\mbox{Poly}}_{f,\Psi,\mathscr{U}}\left(\bm{I}^{\bigtimes}\right)$. 

\subsection{Fundamental Inequalities}\label{sec: Fundamental Inequalities}

The main theorem of this paper is presented below. Theorem~\ref{thm: main 2.3} will give operator inequalties, lower and uppber bounds, of functional with respect to $\Phi(f(\underline{\bm{A}}))$.     
\begin{theorem}\label{thm: main 2.3}
Let $\bm{A}_{j_i}$ be self-adjoint operators with $\Lambda(\bm{A}_{j_i}) \in [m_i, M_i]$ for real scalars $m_i <  M_i$. The mappings $\Phi_{j_1,\ldots,j_n}: \mathscr{B}(\mathfrak{H}) \rightarrow \mathscr{B}(\mathfrak{K})$ are defined by Eq.~\eqref{eq: new phi def}, where $j_i=1,2,\ldots,k_i$ for $i=1,2,\ldots,n$. We have $n$ probability vectors $\bm{w}_i =[w_{i,1},w_{i,2},\cdots, w_{i,k_i}]$ with the dimension $k_i$ for $i=1,2,\ldots,n$, i.e., $\sum\limits_{\ell=1}^{k_i}w_{i,\ell} = 1$. Let $f(\bm{x})$ be any real valued continuous functions with $n$ variables defined on the range $\bigtimes\limits_{i=1}^n [m_i, M_i] \in \mathbb{R}^n$, where $\bigtimes$ is the Cartesian product. Besides, given any $\epsilon>0$, we assume that the function $f(x_1,x_2,\ldots,x_n)$ satisfies the following:
\begin{eqnarray}\label{eq: lower and upper Psi thm: main 2.3}
0 &\leq& \Psi_{\mathscr{U}}(\bm{x}) - f(\bm{x})~~\leq~~\epsilon, \nonumber \\
0 &\leq& f(\bm{x})-\Psi_{\mathscr{L}}(\bm{x})~~\leq~~\epsilon,
\end{eqnarray}
for $\bm{x}\in \bm{I}^{\bigtimes}$.

The function $g(\bm{x})$ is also a real valued continuous function with $n$ variables defined on the range $\bigtimes\limits_{i=1}^n [m_i, M_i]$~\footnote{In this work, we also assume that $g(\underline{\bm{X}})$ is a self-adjoint hypercomplex function given $\underline{\bm{X}}=\bm{X}_1, \bm{X}_2, \ldots, \bm{X}_n$ are $n$ self-adjoint operators.}. We also have a real valued function $F(u,v)$ with operator monotone on the first variable $u$ defined on $U \times V$  such that $f(\bm{I}^{\bigtimes}) \subset U$, and \\
$g\left(\bigtimes\limits_{i=1}^n\left(\bigcup\limits_{\bm{j}=\bm{1}}^{\bm{k}}\bm{w}^{\times}_{\bm{j}}\widetilde{\mbox{Poly}}_{f,\Psi_{\bm{j}},\mathscr{U}}\left(\bm{I}^{\bigtimes}_i\right)\right)\bigcup\bigtimes\limits_{i=1}^n\left(\bigcup\limits_{\bm{j}=\bm{1}}^{\bm{k}}\bm{w}^{\times}_{\bm{j}}\widetilde{\mbox{Poly}}_{f,\Psi_{\bm{j}},\mathscr{L}}\left(\bm{I}^{\bigtimes}_i\right)\right)\right) \subset V$. 

Then, we have the following upper bound:
\begin{eqnarray}\label{eq UB: thm:main 2.3}
F\Bigg(\sum\limits_{\bm{j}=\bm{1}}^{\bm{k}}\bm{w}^{\times}_{\bm{j}}\Phi_{\bm{j}}(f(\underline{\bm{A}}_{\bm{j}})),
g\Bigg(\sum\limits_{\bm{j}=\bm{1}}^{\bm{k}}\bm{w}^{\times}_{\bm{j}}\bm{V}^{\ast}\mbox{Poly}_{f,\Psi_{\bm{j}},\mathscr{U}}(\underline{\bm{A}}_{j_1})\bm{V},\ldots,\sum\limits_{\bm{j}=\bm{1}}^{\bm{k}}\bm{w}^{\times}_{\bm{j}}\bm{V}^{\ast}\mbox{Poly}_{f,\Psi_{\bm{j}},\mathscr{U}}(\underline{\bm{A}}_{j_n})\bm{V}\Bigg)
\Bigg)\nonumber \\
\leq
\max\limits_{x \in \left(\bigcup\limits_{\bm{j}=\bm{1}}^{\bm{k}}\bm{w}^{\times}_{\bm{j}}\widetilde{\mbox{Poly}}_{f,\Psi_{\bm{j}},\mathscr{U}}\left(\bm{I}^{\bigtimes}\right)\right), \bm{x} \in \bigtimes\limits_{i=1}^n\left(\bigcup\limits_{\bm{j}=\bm{1}}^{\bm{k}}\bm{w}^{\times}_{\bm{j}}\widetilde{\mbox{Poly}}_{f,\Psi_{\bm{j}},\mathscr{U}}\left(\bm{I}^{\bigtimes}_i\right)\right)}F\Big(x, g(\bm{x})\Big)\bm{I}_{\mathfrak{K}}.~~~~~~~~~~~~~~
\end{eqnarray}
where $\underline{\bm{A}}_{j_\ell}$ represents a vector of $n$ operaors $\bm{A}_{j_\ell},\bm{A}_{j_\ell},\ldots,\bm{A}_{j_\ell}$ for $\ell=1,2,\ldots,n$, and $\bm{I}^{\bigtimes}_i = \bigtimes_{j=1}^n I_i$. Similarly, we also have the following lower bound:
\begin{eqnarray}\label{eq LB: thm:main 2.3}
F\Bigg(\sum\limits_{\bm{j}=\bm{1}}^{\bm{k}}\bm{w}^{\times}_{\bm{j}}\Phi_{\bm{j}}(f(\underline{\bm{A}}_{\bm{j}})),
g\Bigg(\sum\limits_{\bm{j}=\bm{1}}^{\bm{k}}\bm{w}^{\times}_{\bm{j}}\bm{V}^{\ast}\mbox{Poly}_{f,\Psi_{\bm{j}},\mathscr{L}}(\underline{\bm{A}}_{j_1})\bm{V},\ldots,\sum\limits_{\bm{j}=\bm{1}}^{\bm{k}}\bm{w}^{\times}_{\bm{j}}\bm{V}^{\ast}\mbox{Poly}_{f,\Psi_{\bm{j}},\mathscr{L}}(\underline{\bm{A}}_{j_n})\bm{V}\Bigg)
\Bigg)\nonumber \\
\geq
\min\limits_{x \in \left(\bigcup\limits_{\bm{j}=\bm{1}}^{\bm{k}}\bm{w}^{\times}_{\bm{j}}\widetilde{\mbox{Poly}}_{f,\Psi_{\bm{j}},\mathscr{L}}\left(\bm{I}^{\bigtimes}\right)\right), \bm{x} \in \bigtimes\limits_{i=1}^n\left(\bigcup\limits_{\bm{j}=\bm{1}}^{\bm{k}}\bm{w}^{\times}_{\bm{j}}\widetilde{\mbox{Poly}}_{f,\Psi_{\bm{j}},\mathscr{L}}\left(\bm{I}^{\bigtimes}_i\right)\right)}F\Big(x, g(\bm{x})\Big)\bm{I}_{\mathfrak{K}}.~~~~~~~~~~~~~~
\end{eqnarray}
\end{theorem}
\textbf{Proof:}
We begin with the proof for the upper bound provided by Eq.~\eqref{eq UB: thm:main 2.3}. From Lemma~\ref{lma: phi(f(A)) bounds}, we have 
\begin{eqnarray}\label{eq1: thm:main 2.3}
\Phi(f(\underline{\bm{A}})) \leq \bm{V}^{\ast}\mbox{Poly}_{f,\Psi,\mathscr{U}}(\underline{\bm{A}})\bm{V}.
\end{eqnarray}
By replacing $\underline{\bm{A}}_{\bm{j}}$ with $\underline{\bm{A}}$ in Eq.~\eqref{eq1: thm:main 2.3} and applying $\bm{w}^{\times}_{\bm{j}}$ with respect to each $\underline{\bm{A}}_{\bm{j}}$, we have 
\begin{eqnarray}\label{eq2: thm:main 2.3}
\sum\limits_{\bm{j}=\bm{1}}^{\bm{k}}\bm{w}^{\times}_{\bm{j}}\Phi_{\bm{j}}(f(\underline{\bm{A}}_{\bm{j}}))&\leq&\sum\limits_{\bm{j}=\bm{1}}^{\bm{k}}\bm{w}^{\times}_{\bm{j}}\bm{V}^{\ast}\mbox{Poly}_{f,\Psi_{\bm{j}},\mathscr{U}}(\underline{\bm{A}}_{\bm{j}})\bm{V}.
\end{eqnarray}
According to the function $F(u,v)$ condition, we have
\begin{eqnarray}\label{eq3: thm:main 2.3}
\lefteqn{F\left(\sum\limits_{\bm{j}=\bm{1}}^{\bm{k}}\bm{w}^{\times}_{\bm{j}}\Phi_{\bm{j}}(f(\underline{\bm{A}}_{\bm{j}})),g\Bigg(\sum\limits_{\bm{j}=\bm{1}}^{\bm{k}}\bm{w}^{\times}_{\bm{j}}\bm{V}^{\ast}\mbox{Poly}_{f,\Psi_{\bm{j}},\mathscr{U}}(\underline{\bm{A}}_{j_1})\bm{V},\ldots,\sum\limits_{\bm{j}=\bm{1}}^{\bm{k}}\bm{w}^{\times}_{\bm{j}}\bm{V}^{\ast}\mbox{Poly}_{f,\Psi_{\bm{j}},\mathscr{U}}(\underline{\bm{A}}_{j_n})\bm{V}\Bigg)\right)}\nonumber \\
&\leq&F\left(\sum\limits_{\bm{j}=\bm{1}}^{\bm{k}}\bm{w}^{\times}_{\bm{j}}\bm{V}^{\ast}\mbox{Poly}_{f,\Psi_{\bm{j}},\mathscr{U}}(\underline{\bm{A}}_{\bm{j}})\bm{V},\right. \nonumber \\
& &\left. g\Bigg(\sum\limits_{\bm{j}=\bm{1}}^{\bm{k}}\bm{w}^{\times}_{\bm{j}}\bm{V}^{\ast}\mbox{Poly}_{f,\Psi_{\bm{j}},\mathscr{U}}(\underline{\bm{A}}_{j_1})\bm{V},\ldots,\sum\limits_{\bm{j}=\bm{1}}^{\bm{k}}\bm{w}^{\times}_{\bm{j}}\bm{V}^{\ast}\mbox{Poly}_{f,\Psi_{\bm{j}},\mathscr{U}}(\underline{\bm{A}}_{j_n})\bm{V}\Bigg)\right)\nonumber \\
&\leq&\max\limits_{x \in \left(\bigcup\limits_{\bm{j}=\bm{1}}^{\bm{k}}\bm{w}^{\times}_{\bm{j}}\widetilde{\mbox{Poly}}_{f,\Psi_{\bm{j}},\mathscr{U}}\left(\bm{I}^{\bigtimes}\right)\right), \bm{x} \in \bigtimes\limits_{i=1}^n\left(\bigcup\limits_{\bm{j}=\bm{1}}^{\bm{k}}\bm{w}^{\times}_{\bm{j}}\widetilde{\mbox{Poly}}_{f,\Psi_{\bm{j}},\mathscr{U}}\left(\bm{I}^{\bigtimes}_i\right)\right)}F\Big(x, g(\bm{x})\Big)\bm{I}_{\mathfrak{K}},
\end{eqnarray}
where the last inequality comes from that the spectrum $\Lambda\left(\sum\limits_{\bm{j}=\bm{1}}^{\bm{k}}\bm{w}^{\times}_{\bm{j}}\bm{V}^{\ast}\mbox{Poly}_{f,\Psi_{\bm{j}},\mathscr{U}}(\underline{\bm{A}}_{\bm{j}})\bm{V}\right)$ is in the range of \\ $\left(\bigcup\limits_{\bm{j}=\bm{1}}^{\bm{k}}\bm{w}^{\times}_{\bm{j}}\widetilde{\mbox{Poly}}_{f,\Psi_{\bm{j}},\mathscr{U}}\left(\bm{I}^{\bigtimes}\right)\right)$, and the spectrum $\Lambda\left(\sum\limits_{\bm{j}=\bm{1}}^{\bm{k}}\bm{w}^{\times}_{\bm{j}}\bm{V}^{\ast}\mbox{Poly}_{f,\Psi_{\bm{j}},\mathscr{U}}(\underline{\bm{A}}_{j_i})\bm{V}\right)$ is in the range of \\ $\left(\bigcup\limits_{\bm{j}=\bm{1}}^{\bm{k}}\bm{w}^{\times}_{\bm{j}}\widetilde{\mbox{Poly}}_{f,\Psi_{\bm{j}},\mathscr{U}}\left(\bm{I}_i^{\bigtimes}\right)\right)$ for $i=1,2,\ldots,n$. The desired inequality provided by Eq.~\eqref{eq UB: thm:main 2.3} is established. 

Now, we will prove the lower bound provided by Eq.~\eqref{eq LB: thm:main 2.3}. From Lemma~\ref{lma: phi(f(A)) bounds}, we have 
\begin{eqnarray}\label{eq4: thm:main 2.3}
\Phi(f(\underline{\bm{A}})) \geq \bm{V}^{\ast}\mbox{Poly}_{f,\Psi,\mathscr{L}}(\underline{\bm{A}})\bm{V}.
\end{eqnarray}
By replacing $\underline{\bm{A}}$ with $\underline{\bm{A}}_{\bm{j}}$ in Eq.~\eqref{eq4: thm:main 2.3} and applying $\bm{w}^{\times}_{\bm{j}}$ with respect to each $\underline{\bm{A}}_{\bm{j}}$, we have 
\begin{eqnarray}\label{eq5: thm:main 2.3}
\sum\limits_{\bm{j}=\bm{1}}^{\bm{k}}\bm{w}^{\times}_{\bm{j}}\Phi_{\bm{j}}(f(\underline{\bm{A}}_{\bm{j}}))&\geq&\sum\limits_{\bm{j}=\bm{1}}^{\bm{k}}\bm{w}^{\times}_{\bm{j}}\bm{V}^{\ast}\mbox{Poly}_{f,\Psi_{\bm{j}},\mathscr{L}}(\underline{\bm{A}}_{\bm{j}})\bm{V}.
\end{eqnarray}
From the function $F(u,v)$ condition, we also have
\begin{eqnarray}\label{eq3: thm:main 2.3}
\lefteqn{F\left(\sum\limits_{\bm{j}=\bm{1}}^{\bm{k}}\bm{w}^{\times}_{\bm{j}}\Phi_{\bm{j}}(f(\underline{\bm{A}}_{\bm{j}})),g\Bigg(\sum\limits_{\bm{j}=\bm{1}}^{\bm{k}}\bm{w}^{\times}_{\bm{j}}\bm{V}^{\ast}\mbox{Poly}_{f,\Psi_{\bm{j}},\mathscr{L}}(\underline{\bm{A}}_{j_1})\bm{V},\ldots,\sum\limits_{\bm{j}=\bm{1}}^{\bm{k}}\bm{w}^{\times}_{\bm{j}}\bm{V}^{\ast}\mbox{Poly}_{f,\Psi_{\bm{j}},\mathscr{L}}(\underline{\bm{A}}_{j_n})\bm{V}\Bigg)\right)}\nonumber \\
&\geq&F\left(\sum\limits_{\bm{j}=\bm{1}}^{\bm{k}}\bm{w}^{\times}_{\bm{j}}\bm{V}^{\ast}\mbox{Poly}_{f,\Psi_{\bm{j}},\mathscr{L}}(\underline{\bm{A}}_{\bm{j}})\bm{V},\right. \nonumber \\
& &\left. g\Bigg(\sum\limits_{\bm{j}=\bm{1}}^{\bm{k}}\bm{w}^{\times}_{\bm{j}}\bm{V}^{\ast}\mbox{Poly}_{f,\Psi_{\bm{j}},\mathscr{L}}(\underline{\bm{A}}_{j_1})\bm{V},\ldots,\sum\limits_{\bm{j}=\bm{1}}^{\bm{k}}\bm{w}^{\times}_{\bm{j}}\bm{V}^{\ast}\mbox{Poly}_{f,\Psi_{\bm{j}},\mathscr{L}}(\underline{\bm{A}}_{j_n})\bm{V}\Bigg)\right)\nonumber \\
&\geq&\min\limits_{x \in \left(\bigcup\limits_{\bm{j}=\bm{1}}^{\bm{k}}\bm{w}^{\times}_{\bm{j}}\widetilde{\mbox{Poly}}_{f,\Psi_{\bm{j}},\mathscr{L}}\left(\bm{I}^{\bigtimes}\right)\right), \bm{x} \in \bigtimes\limits_{i=1}^n\left(\bigcup\limits_{\bm{j}=\bm{1}}^{\bm{k}}\bm{w}^{\times}_{\bm{j}}\widetilde{\mbox{Poly}}_{f,\Psi_{\bm{j}},\mathscr{L}}\left(\bm{I}^{\bigtimes}_i\right)\right)}F\Big(x, g(\bm{x})\Big)\bm{I}_{\mathfrak{K}},
\end{eqnarray}
where the last inequality comes from that the spectrum $\Lambda\left(\sum\limits_{\bm{j}=\bm{1}}^{\bm{k}}\bm{w}^{\times}_{\bm{j}}\bm{V}^{\ast}\mbox{Poly}_{f,\Psi_{\bm{j}},\mathscr{L}}(\underline{\bm{A}}_{\bm{j}})\bm{V}\right)$ is in the range of \\ $\left(\bigcup\limits_{\bm{j}=\bm{1}}^{\bm{k}}\bm{w}^{\times}_{\bm{j}}\widetilde{\mbox{Poly}}_{f,\Psi_{\bm{j}},\mathscr{L}}\left(\bm{I}^{\bigtimes}\right)\right)$, and the spectrum $\Lambda\left(\sum\limits_{\bm{j}=\bm{1}}^{\bm{k}}\bm{w}^{\times}_{\bm{j}}\bm{V}^{\ast}\mbox{Poly}_{f,\Psi_{\bm{j}},\mathscr{L}}(\underline{\bm{A}}_{j_i})\bm{V}\right)$ is in the range of \\ $\left(\bigcup\limits_{\bm{j}=\bm{1}}^{\bm{k}}\bm{w}^{\times}_{\bm{j}}\widetilde{\mbox{Poly}}_{f,\Psi_{\bm{j}},\mathscr{L}}\left(\bm{I}_i^{\bigtimes}\right)\right)$ for $i=1,2,\ldots,n$. The desired inequality provided by Eq.~\eqref{eq LB: thm:main 2.3} is established. 
$\hfill \Box$

To establish Theorem~\ref{thm:main 2.4}, we define the function $F(u,v)$ as follows:
\begin{eqnarray}\label{eq: F u v }
F(u,v) &=& u - \alpha v,
\end{eqnarray} 
where $\alpha \in \mathbb{R}$. This theorem sets the groundwork for inequalities involving multivariate hypercomplex functions, addressing both ratio and difference variations, which will be detailed in the following sections.

\begin{theorem}\label{thm:main 2.4}
Let $\bm{A}_{j_i}$ be self-adjoint operators with $\Lambda(\bm{A}_{j_i}) \in [m_i, M_i]$ for real scalars $m_i <  M_i$. The mappings $\Phi_{j_1,\ldots,j_n}: \mathscr{B}(\mathfrak{H}) \rightarrow \mathscr{B}(\mathfrak{K})$ are defined by Eq.~\eqref{eq: new phi def}, where $j_i=1,2,\ldots,k_i$ for $i=1,2,\ldots,n$. We have $n$ probability vectors $\bm{w}_i =[w_{i,1},w_{i,2},\cdots, w_{i,k_i}]$ with the dimension $k_i$ for $i=1,2,\ldots,n$, i.e., $\sum\limits_{\ell=1}^{k_i}w_{i,\ell} = 1$. Let $f(\bm{x})$ be any real valued continuous functions with $n$ variables defined on the range $\bigtimes\limits_{i=1}^n [m_i, M_i] \in \mathbb{R}^n$, where $\bigtimes$ is the Cartesian product. Besides, given any $\epsilon>0$, we assume that the function $f(\bm{x})$ satisfies the following:
\begin{eqnarray}\label{eq: lower and upper Psi thm: main 2.4}
0 &\leq& \Psi_{\mathscr{U}}(\bm{x}) - f(\bm{x})~~\leq~~\epsilon, \nonumber \\
0 &\leq& f(\bm{x})-\Psi_{\mathscr{L}}(\bm{x})~~\leq~~\epsilon,
\end{eqnarray}
for $\bm{x}\in \bm{I}^{\bigtimes}$. The function $g(\bm{x})$ is also a real valued continuous function with $n$ variables defined on the range $\bigtimes\limits_{i=1}^n [m_i, M_i]$. We also have a real valued function $F(u,v)$ defined by Eq.~\eqref{eq: F u v } with support domain on $U \times V$ such that $f(\bm{I}^{\bigtimes}) \subset U$, and \\
$g\left(\bigtimes\limits_{i=1}^n\left(\bigcup\limits_{\bm{j}=\bm{1}}^{\bm{k}}\bm{w}^{\times}_{\bm{j}}\widetilde{\mbox{Poly}}_{f,\Psi_{\bm{j}},\mathscr{U}}\left(\bm{I}^{\bigtimes}_i\right)\right)\bigcup\bigtimes\limits_{i=1}^n\left(\bigcup\limits_{\bm{j}=\bm{1}}^{\bm{k}}\bm{w}^{\times}_{\bm{j}}\widetilde{\mbox{Poly}}_{f,\Psi_{\bm{j}},\mathscr{L}}\left(\bm{I}^{\bigtimes}_i\right)\right)\right) \subset V$. 

Then, we have the following upper bound:
\begin{eqnarray}\label{eq UB: thm:main 2.4}
\lefteqn{\sum\limits_{\bm{j}=\bm{1}}^{\bm{k}}\bm{w}^{\times}_{\bm{j}}\Phi_{\bm{j}}(f(\underline{\bm{A}}_{\bm{j}}))}\nonumber \\
&\leq&\alpha g\Bigg(\sum\limits_{\bm{j}=\bm{1}}^{\bm{k}}\bm{w}^{\times}_{\bm{j}}\bm{V}^{\ast}\mbox{Poly}_{f,\Psi_{\bm{j}},\mathscr{U}}(\underline{\bm{A}}_{j_1})\bm{V},\ldots,\sum\limits_{\bm{j}=\bm{1}}^{\bm{k}}\bm{w}^{\times}_{\bm{j}}\bm{V}^{\ast}\mbox{Poly}_{f,\Psi_{\bm{j}},\mathscr{U}}(\underline{\bm{A}}_{j_n})\bm{V}\Bigg) \nonumber \\
& &+ 
\max\limits_{x \in \left(\bigcup\limits_{\bm{j}=\bm{1}}^{\bm{k}}\bm{w}^{\times}_{\bm{j}}\widetilde{\mbox{Poly}}_{f,\Psi_{\bm{j}},\mathscr{U}}\left(\bm{I}^{\bigtimes}\right)\right), \bm{x} \in \bigtimes\limits_{i=1}^n\left(\bigcup\limits_{\bm{j}=\bm{1}}^{\bm{k}}\bm{w}^{\times}_{\bm{j}}\widetilde{\mbox{Poly}}_{f,\Psi_{\bm{j}},\mathscr{U}}\left(\bm{I}^{\bigtimes}_i\right)\right)}(x-\alpha g(\bm{x}))\bm{I}_{\mathfrak{K}}.
\end{eqnarray}
Similarly, we also have the following lower bound:
\begin{eqnarray}\label{eq LB: thm:main 2.4}
\lefteqn{\sum\limits_{\bm{j}=\bm{1}}^{\bm{k}}\bm{w}^{\times}_{\bm{j}}\Phi_{\bm{j}}(f(\underline{\bm{A}}_{\bm{j}}))}\nonumber \\
&\geq&\alpha g\Bigg(\sum\limits_{\bm{j}=\bm{1}}^{\bm{k}}\bm{w}^{\times}_{\bm{j}}\bm{V}^{\ast}\mbox{Poly}_{f,\Psi_{\bm{j}},\mathscr{L}}(\underline{\bm{A}}_{j_1})\bm{V},\ldots,\sum\limits_{\bm{j}=\bm{1}}^{\bm{k}}\bm{w}^{\times}_{\bm{j}}\bm{V}^{\ast}\mbox{Poly}_{f,\Psi_{\bm{j}},\mathscr{L}}(\underline{\bm{A}}_{j_n})\bm{V}\Bigg) \nonumber \\
& &+ 
\min\limits_{x \in \left(\bigcup\limits_{\bm{j}=\bm{1}}^{\bm{k}}\bm{w}^{\times}_{\bm{j}}\widetilde{\mbox{Poly}}_{f,\Psi_{\bm{j}},\mathscr{L}}\left(\bm{I}^{\bigtimes}\right)\right), \bm{x} \in \bigtimes\limits_{i=1}^n\left(\bigcup\limits_{\bm{j}=\bm{1}}^{\bm{k}}\bm{w}^{\times}_{\bm{j}}\widetilde{\mbox{Poly}}_{f,\Psi_{\bm{j}},\mathscr{L}}\left(\bm{I}^{\bigtimes}_i\right)\right)}(x-\alpha g(\bm{x}))\bm{I}_{\mathfrak{K}}.
\end{eqnarray}
\end{theorem} 
\textbf{Proof:}
By setting $F(u,v)=u - \alpha v$ in Eq.~\eqref{eq UB: thm:main 2.3} in Theorem~\ref{thm: main 2.3}, we have the desired inequality provided by Eq.~\eqref{eq UB: thm:main 2.4}. Similarly, By setting $F(u,v)=u - \alpha v$ in Eq.~\eqref{eq LB: thm:main 2.3} in Theorem~\ref{thm: main 2.3}, we have the desired inequality provided by Eq.~\eqref{eq LB: thm:main 2.4}.
$\hfill \Box$

Corollary~\ref{cor: special cases of g func} below is provided to give upper and lower bounds for special types of the function $g$ by applying Theorem~\ref{thm:main 2.4}.

\begin{corollary}\label{cor: special cases of g func}
Let $\bm{A}_{j_i}$ be self-adjoint operators with $\Lambda(\bm{A}_{j_i}) \in [m_i, M_i]$ for real scalars $m_i <  M_i$. The mappings $\Phi_{j_1,\ldots,j_n}: \mathscr{B}(\mathfrak{H}) \rightarrow \mathscr{B}(\mathfrak{K})$ are normalized positive linear maps, where $j_i=1,2,\ldots,k_i$ for $i=1,2,\ldots,n$. We have $n$ probability vectors $\bm{w}_i =[w_{i,1},w_{i,2},\cdots, w_{i,k_i}]$ with the dimension $k_i$ for $i=1,2,\ldots,n$, i.e., $\sum\limits_{\ell=1}^{k_i}w_{i,\ell} = 1$. Let $f(\bm{x})$ be any real valued continuous functions with $n$ variables defined on the range $\bigtimes\limits_{i=1}^n [m_i, M_i] \in \mathbb{R}^n$, where $\times$ is the Cartesian product. Besides, given any $\epsilon>0$, we assume that the function $f(\bm{x})$ satisfies the following:
\begin{eqnarray}\label{eq: lower and upper f cor: special cases of g func} 
0 &\leq& \Psi_{\mathscr{U}}(\bm{x}) - f(\bm{x})~~\leq~~\epsilon, \nonumber \\
0 &\leq& f(\bm{x})-\Psi_{\mathscr{L}}(\bm{x})~~\leq~~\epsilon,
\end{eqnarray}
for $\bm{x}\in \bm{I}^{\bigtimes}$.

(I) If $g(\bm{x}) = \Big(\sum\limits_{i=1}^n \beta_i x_i\Big)^q$, where $\beta_i \geq 0$, $q \in \mathbb{R}$ with $\left(\sum\limits_{\bm{j}=\bm{1}}^{\bm{k}}\bm{w}^{\times}_{\bm{j}}\bm{V}^{\ast}\mbox{Poly}_{f,\Psi_{\bm{j}},\mathscr{L}}(\underline{\bm{A}}_{j_i})\bm{V}\right)\geq\bm{0}$ and $\left(\sum\limits_{\bm{j}=\bm{1}}^{\bm{k}}\bm{w}^{\times}_{\bm{j}}\bm{V}^{\ast}\mbox{Poly}_{f,\Psi_{\bm{j}},\mathscr{U}}(\underline{\bm{A}}_{j_i})\bm{V}\right)\geq\bm{0}$ for $i=1,2,\ldots,n$, we have the upper bound for $\sum\limits_{\bm{j}=\bm{1}}^{\bm{k}}\bm{w}^{\times}_{\bm{j}}\Phi_{\bm{j}}(f(\underline{\bm{A}}_{\bm{j}}))$:
\begin{eqnarray}\label{eq power U: cor: special cases of g func}
\sum\limits_{\bm{j}=\bm{1}}^{\bm{k}}\bm{w}^{\times}_{\bm{j}}\Phi_{\bm{j}}(f(\underline{\bm{A}}_{\bm{j}}))&\leq&\alpha \left(\sum\limits_{i=1}^n \beta_i\left(\sum\limits_{\bm{j}=\bm{1}}^{\bm{k}}\bm{w}^{\times}_{\bm{j}}\bm{V}^{\ast}\mbox{Poly}_{f,\Psi_{\bm{j}},\mathscr{U}}(\underline{\bm{A}}_{j_i})\bm{V}\right)\right)^q\nonumber \\
& &+
\max\limits_{x \in \left(\bigcup\limits_{\bm{j}=\bm{1}}^{\bm{k}}\bm{w}^{\times}_{\bm{j}}\widetilde{\mbox{Poly}}_{f,\Psi_{\bm{j}},\mathscr{U}}\left(\bm{I}^{\bigtimes}\right)\right),x_i \in \left(\bigcup\limits_{\bm{j}=\bm{1}}^{\bm{k}}\bm{w}^{\times}_{\bm{j}}\widetilde{\mbox{Poly}}_{f,\Psi_{\bm{j}},\mathscr{U}}\left(\bm{I}^{\bigtimes}_i\right)\right)}\Big(x-\alpha \Big(\sum\limits_{i=1}^n \beta_i x_i\Big)^q\Big)\bm{I}_{\mathfrak{K}},\nonumber \\
\end{eqnarray}
and we have the lower bound for $\sum\limits_{\bm{j}=\bm{1}}^{\bm{k}}\bm{w}^{\times}_{\bm{j}}\Phi_{\bm{j}}(f(\underline{\bm{A}}_{\bm{j}}))$:
\begin{eqnarray}\label{eq power L: cor: special cases of g func}
\sum\limits_{\bm{j}=\bm{1}}^{\bm{k}}\bm{w}^{\times}_{\bm{j}}\Phi_{\bm{j}}(f(\underline{\bm{A}}_{\bm{j}}))&\geq&\alpha \left(\sum\limits_{i=1}^n \beta_i \left(\sum\limits_{\bm{j}=\bm{1}}^{\bm{k}}\bm{w}^{\times}_{\bm{j}}\bm{V}^{\ast}\mbox{Poly}_{f,\Psi_{\bm{j}},\mathscr{L}}(\underline{\bm{A}}_{j_i})\bm{V}\right) \right)^q\nonumber \\
& &+
\min\limits_{x \in \left(\bigcup\limits_{\bm{j}=\bm{1}}^{\bm{k}}\bm{w}^{\times}_{\bm{j}}\widetilde{\mbox{Poly}}_{f,\Psi_{\bm{j}},\mathscr{L}}\left(\bm{I}^{\bigtimes}\right)\right),x_i \in \left(\bigcup\limits_{\bm{j}=\bm{1}}^{\bm{k}}\bm{w}^{\times}_{\bm{j}}\widetilde{\mbox{Poly}}_{f,\Psi_{\bm{j}},\mathscr{L}}\left(\bm{I}^{\bigtimes}_i\right)\right)}\Big(x-\alpha \Big(\sum\limits_{i=1}^n \beta_i x_i\Big)^q\Big)\bm{I}_{\mathfrak{K}}.\nonumber \\
\end{eqnarray}

(II) If $g(x_1,\ldots,x_n) = \log\Big(\sum\limits_{i=1}^n \beta_i x_i\Big)$, where $\beta_i \geq 0$ with $\left(\sum\limits_{\bm{j}=\bm{1}}^{\bm{k}}\bm{w}^{\times}_{\bm{j}}\bm{V}^{\ast}\mbox{Poly}_{f,\Psi_{\bm{j}},\mathscr{L}}(\underline{\bm{A}}_{j_i})\bm{V}\right)>\bm{0}$ and $\left(\sum\limits_{\bm{j}=\bm{1}}^{\bm{k}}\bm{w}^{\times}_{\bm{j}}\bm{V}^{\ast}\mbox{Poly}_{f,\Psi_{\bm{j}},\mathscr{U}}(\underline{\bm{A}}_{j_i})\bm{V}\right)>\bm{0}$ for $i=1,2,\ldots,n$, we have the upper bound for $\sum\limits_{\bm{j}=\bm{1}}^{\bm{k}}\bm{w}^{\times}_{\bm{j}}\Phi_{\bm{j}}(f(\underline{\bm{A}}_{\bm{j}}))$:
\begin{eqnarray}\label{eq log U: cor: special cases of g func}
\lefteqn{\sum\limits_{\bm{j}=\bm{1}}^{\bm{k}}\bm{w}^{\times}_{\bm{j}}\Phi_{\bm{j}}(f(\underline{\bm{A}}_{\bm{j}}))}\nonumber \\
&\leq&\alpha \log\left(\sum\limits_{i=1}^n \beta_i \left(\sum\limits_{\bm{j}=\bm{1}}^{\bm{k}}\bm{w}^{\times}_{\bm{j}}\bm{V}^{\ast}\mbox{Poly}_{f,\Psi_{\bm{j}},\mathscr{U}}(\underline{\bm{A}}_{j_i})\bm{V}\right)\right)\nonumber \\
& &+
\max\limits_{x \in \left(\bigcup\limits_{\bm{j}=\bm{1}}^{\bm{k}}\bm{w}^{\times}_{\bm{j}}\widetilde{\mbox{Poly}}_{f,\Psi_{\bm{j}},\mathscr{U}}\left(\bm{I}^{\bigtimes}\right)\right),x_i \in \left(\bigcup\limits_{\bm{j}=\bm{1}}^{\bm{k}}\bm{w}^{\times}_{\bm{j}}\widetilde{\mbox{Poly}}_{f,\Psi_{\bm{j}},\mathscr{U}}\left(\bm{I}^{\bigtimes}_i\right)\right)}\Big(x-\alpha \log\Big(\sum\limits_{i=1}^n \beta_i x_i\Big)\Big)\bm{I}_{\mathfrak{K}},\nonumber \\
\end{eqnarray}
and we have the lower bound for $\sum\limits_{\bm{j}=\bm{1}}^{\bm{k}}\bm{w}^{\times}_{\bm{j}}\Phi_{\bm{j}}(f(\underline{\bm{A}}_{\bm{j}}))$:
\begin{eqnarray}\label{eq log L: cor: special cases of g func}
\lefteqn{\sum\limits_{\bm{j}=\bm{1}}^{\bm{k}}\bm{w}^{\times}_{\bm{j}}\Phi_{\bm{j}}(f(\underline{\bm{A}}_{\bm{j}}))}\nonumber \\
&\geq&\alpha \log\left(\sum\limits_{i=1}^n \beta_i \left(\sum\limits_{\bm{j}=\bm{1}}^{\bm{k}}\bm{w}^{\times}_{\bm{j}}\bm{V}^{\ast}\mbox{Poly}_{f,\Psi_{\bm{j}},\mathscr{L}}(\underline{\bm{A}}_{j_i})\bm{V}\right)\right)\nonumber \\
& &+
\min\limits_{x \in \left(\bigcup\limits_{\bm{j}=\bm{1}}^{\bm{k}}\bm{w}^{\times}_{\bm{j}}\widetilde{\mbox{Poly}}_{f,\Psi_{\bm{j}},\mathscr{L}}\left(\bm{I}^{\bigtimes}\right)\right),x_i \in \left(\bigcup\limits_{\bm{j}=\bm{1}}^{\bm{k}}\bm{w}^{\times}_{\bm{j}}\widetilde{\mbox{Poly}}_{f,\Psi_{\bm{j}},\mathscr{L}}\left(\bm{I}^{\bigtimes}_i\right)\right)}\Big(x-\alpha \log\Big(\sum\limits_{i=1}^n \beta_i x_i\Big)\Big)\bm{I}_{\mathfrak{K}}.\nonumber \\
\end{eqnarray}

(III) If $g(x_1,\ldots,x_n) = \exp\Big(\sum\limits_{i=1}^n \beta_i x_i\Big)$, we have the upper bound for $\sum\limits_{\bm{j}=\bm{1}}^{\bm{k}}\bm{w}^{\times}_{\bm{j}}\Phi_{\bm{j}}(f(\underline{\bm{A}}_{\bm{j}}))$:
\begin{eqnarray}\label{eq exp U: cor: special cases of g func}
\lefteqn{\sum\limits_{\bm{j}=\bm{1}}^{\bm{k}}\bm{w}^{\times}_{\bm{j}}\Phi_{\bm{j}}(f(\underline{\bm{A}}_{\bm{j}}))}\nonumber \\
&\leq&\alpha \exp\left(\sum\limits_{i=1}^n \beta_i \left(\sum\limits_{\bm{j}=\bm{1}}^{\bm{k}}\bm{w}^{\times}_{\bm{j}}\bm{V}^{\ast}\mbox{Poly}_{f,\Psi_{\bm{j}},\mathscr{U}}(\underline{\bm{A}}_{j_i})\bm{V}\right)\right)\nonumber \\
& &+
\max\limits_{x \in \left(\bigcup\limits_{\bm{j}=\bm{1}}^{\bm{k}}\bm{w}^{\times}_{\bm{j}}\widetilde{\mbox{Poly}}_{f,\Psi_{\bm{j}},\mathscr{U}}\left(\bm{I}^{\bigtimes}\right)\right),x_i \in \left(\bigcup\limits_{\bm{j}=\bm{1}}^{\bm{k}}\bm{w}^{\times}_{\bm{j}}\widetilde{\mbox{Poly}}_{f,\Psi_{\bm{j}},\mathscr{U}}\left(\bm{I}^{\bigtimes}_i\right)\right)}\Big(x-\alpha \exp\Big(\sum\limits_{i=1}^n \beta_i x_i\Big)\Big)\bm{I}_{\mathfrak{K}},\nonumber \\
\end{eqnarray}
and we have the lower bound for $\sum\limits_{\bm{j}=\bm{1}}^{\bm{k}}\bm{w}^{\times}_{\bm{j}}\Phi_{\bm{j}}(f(\underline{\bm{A}}_{\bm{j}}))$:
\begin{eqnarray}\label{eq exp L: cor: special cases of g func}
\lefteqn{\sum\limits_{\bm{j}=\bm{1}}^{\bm{k}}\bm{w}^{\times}_{\bm{j}}\Phi_{\bm{j}}(f(\underline{\bm{A}}_{\bm{j}}))}\nonumber \\
&\geq&\alpha \exp\left(\sum\limits_{i=1}^n \beta_i \left(\sum\limits_{\bm{j}=\bm{1}}^{\bm{k}}\bm{w}^{\times}_{\bm{j}}\bm{V}^{\ast}\mbox{Poly}_{f,\Psi_{\bm{j}},\mathscr{L}}(\underline{\bm{A}}_{j_i})\bm{V}\right)\right)\nonumber \\
& &+
\min\limits_{x \in \left(\bigcup\limits_{\bm{j}=\bm{1}}^{\bm{k}}\bm{w}^{\times}_{\bm{j}}\widetilde{\mbox{Poly}}_{f,\Psi_{\bm{j}},\mathscr{L}}\left(\bm{I}^{\bigtimes}\right)\right),x_i \in \left(\bigcup\limits_{\bm{j}=\bm{1}}^{\bm{k}}\bm{w}^{\times}_{\bm{j}}\widetilde{\mbox{Poly}}_{f,\Psi_{\bm{j}},\mathscr{L}}\left(\bm{I}^{\bigtimes}_i\right)\right)}\Big(x-\alpha \exp\Big(\sum\limits_{i=1}^n \beta_i x_i\Big)\Big)\bm{I}_{\mathfrak{K}}.\nonumber \\
\end{eqnarray}
\end{corollary}
\textbf{Proof:}
This corollary is derived by applying Theorem~\ref{thm:main 2.4} to different functions:
\\
- Part (I) is obtained by using the function \( g(\bm{x}) = \left( \sum\limits_{i=1}^n \beta_i x_i \right)^q \).\\
- Part (II) is obtained by using the function \( g(\bm{x}) = \log\left( \sum\limits_{i=1}^n \beta_i x_i \right) \).\\
- Part (III) is obtained by using the function \( g(\bm{x}) = \exp\left( \sum\limits_{i=1}^n \beta_i x_i \right) \).
$\hfill \Box$

\section{Multivariate Hypercomplex Functions Inequalities: Ratio Kind}\label{sec: Multivariate Hypercomplex Functions Inequalities: Ratio Kind}

In this section, we will derive the lower and upper bounds for $\sum\limits_{\bm{j}=\bm{1}}^{\bm{k}}\bm{w}^{\times}_{\bm{j}}\Phi_{\bm{j}}(f(\underline{\bm{A}}_{\bm{j}}))$ in terms of ratio criteria related to the function $g$.

\begin{theorem}\label{thm: 2.9}
Let $\bm{A}_{j_i}$ be self-adjoint operators with $\Lambda(\bm{A}_{j_i}) \in [m_i, M_i]$ for real scalars $m_i <  M_i$. The mappings $\Phi_{j_1,\ldots,j_n}: \mathscr{B}(\mathfrak{H}) \rightarrow \mathscr{B}(\mathfrak{K})$ are defined by Eq.~\eqref{eq: new phi def}, where $j_i=1,2,\ldots,k_i$ for $i=1,2,\ldots,n$. We have $n$ probability vectors $\bm{w}_i =[w_{i,1},w_{i,2},\cdots, w_{i,k_i}]$ with the dimension $k_i$ for $i=1,2,\ldots,n$, i.e., $\sum\limits_{\ell=1}^{k_i}w_{i,\ell} = 1$. Let $f(\bm{x})$ be any real valued continuous functions with $n$ variables defined on the range $\bigtimes\limits_{i=1}^n [m_i, M_i] \in \mathbb{R}^n$, where $\bigtimes$ is the Cartesian product. Besides, given any $\epsilon>0$, we assume that the function $f(\bm{x})$ satisfies the following:
\begin{eqnarray}\label{eq: lower and upper Psi thm: main 2.4}
0 &\leq& \Psi_{\mathscr{U}}(\bm{x}) - f(\bm{x})~~\leq~~\epsilon, \nonumber \\
0 &\leq& f(\bm{x})-\Psi_{\mathscr{L}}(\bm{x})~~\leq~~\epsilon,
\end{eqnarray}
for $\bm{x}\in \bm{I}^{\bigtimes}$. The function $g(\bm{x})$ is also a real valued continuous function with $n$ variables defined on the range $\bigtimes\limits_{i=1}^n [m_i, M_i]$. The function $g$ is also a real valued continuous function defined on the range $\bigtimes\limits_{i=1}^n\left(\bigcup\limits_{\bm{j}=\bm{1}}^{\bm{k}}\bm{w}^{\times}_{\bm{j}}\widetilde{\mbox{Poly}}_{f,\Psi_{\bm{j}},\mathscr{L}}\left(\bm{I}^{\bigtimes}_i\right)\right) \bigcup \bigtimes\limits_{i=1}^n\left(\bigcup\limits_{\bm{j}=\bm{1}}^{\bm{k}}\bm{w}^{\times}_{\bm{j}}\widetilde{\mbox{Poly}}_{f,\Psi_{\bm{j}},\mathscr{U}}\left(\bm{I}^{\bigtimes}_i\right)\right)$, and $g(\bm{x})\neq 0$ for \\
$\bm{x} \in \bigtimes\limits_{i=1}^n\left(\bigcup\limits_{\bm{j}=\bm{1}}^{\bm{k}}\bm{w}^{\times}_{\bm{j}}\widetilde{\mbox{Poly}}_{f,\Psi_{\bm{j}},\mathscr{L}}\left(\bm{I}^{\bigtimes}_i\right)\right) \bigcup \bigtimes\limits_{i=1}^n\left(\bigcup\limits_{\bm{j}=\bm{1}}^{\bm{k}}\bm{w}^{\times}_{\bm{j}}\widetilde{\mbox{Poly}}_{f,\Psi_{\bm{j}},\mathscr{U}}\left(\bm{I}^{\bigtimes}_i\right)\right)$. 

(I) If we also assume that $g\Bigg(\sum\limits_{\bm{j}=\bm{1}}^{\bm{k}}\bm{w}^{\times}_{\bm{j}}\bm{V}^{\ast}\mbox{Poly}_{f,\Psi_{\bm{j}},\mathscr{U}}(\underline{\bm{A}}_{j_1})\bm{V},\ldots,\sum\limits_{\bm{j}=\bm{1}}^{\bm{k}}\bm{w}^{\times}_{\bm{j}}\bm{V}^{\ast}\mbox{Poly}_{f,\Psi_{\bm{j}},\mathscr{U}}(\underline{\bm{A}}_{j_n})\bm{V}\Bigg) > \bm{0}$ and $g\Bigg(\sum\limits_{\bm{j}=\bm{1}}^{\bm{k}}\bm{w}^{\times}_{\bm{j}}\bm{V}^{\ast}\mbox{Poly}_{f,\Psi_{\bm{j}},\mathscr{L}}(\underline{\bm{A}}_{j_1})\bm{V},\ldots,\sum\limits_{\bm{j}=\bm{1}}^{\bm{k}}\bm{w}^{\times}_{\bm{j}}\bm{V}^{\ast}\mbox{Poly}_{f,\Psi_{\bm{j}},\mathscr{L}}(\underline{\bm{A}}_{j_n})\bm{V}\Bigg) > \bm{0}$, then, we have the following upper bound for $\sum\limits_{\bm{j}=\bm{1}}^{\bm{k}}\bm{w}^{\times}_{\bm{j}}\Phi_{\bm{j}}(f(\underline{\bm{A}}_{\bm{j}}))$:
\begin{eqnarray}\label{eq1 UB: thm: 2.9 pos g}
\lefteqn{\sum\limits_{\bm{j}=\bm{1}}^{\bm{k}}\bm{w}^{\times}_{\bm{j}}\Phi_{\bm{j}}(f(\underline{\bm{A}}_{\bm{j}}))}\nonumber \\
&\leq&\left[\max\limits_{x \in \left(\bigcup\limits_{\bm{j}=\bm{1}}^{\bm{k}}\bm{w}^{\times}_{\bm{j}}\widetilde{\mbox{Poly}}_{f,\Psi_{\bm{j}},\mathscr{U}}\left(\bm{I}^{\bigtimes}\right)\right), \bm{x} \in \bigtimes\limits_{i=1}^n\left(\bigcup\limits_{\bm{j}=\bm{1}}^{\bm{k}}\bm{w}^{\times}_{\bm{j}}\widetilde{\mbox{Poly}}_{f,\Psi_{\bm{j}},\mathscr{U}}\left(\bm{I}^{\bigtimes}_i\right)\right)}xg^{-1}(\bm{x})\right] \nonumber \\
&&\times g\Bigg(\sum\limits_{\bm{j}=\bm{1}}^{\bm{k}}\bm{w}^{\times}_{\bm{j}}\bm{V}^{\ast}\mbox{Poly}_{f,\Psi_{\bm{j}},\mathscr{U}}(\underline{\bm{A}}_{j_1})\bm{V},\ldots,\sum\limits_{\bm{j}=\bm{1}}^{\bm{k}}\bm{w}^{\times}_{\bm{j}}\bm{V}^{\ast}\mbox{Poly}_{f,\Psi_{\bm{j}},\mathscr{U}}(\underline{\bm{A}}_{j_n})\bm{V}\Bigg);
\end{eqnarray}
and, the following lower bound for $\sum\limits_{\bm{j}=\bm{1}}^{\bm{k}}\bm{w}^{\times}_{\bm{j}}\Phi_{\bm{j}}(f(\underline{\bm{A}}_{\bm{j}}))$:
\begin{eqnarray}\label{eq1 LB: thm: 2.9 pos g}
\lefteqn{\sum\limits_{\bm{j}=\bm{1}}^{\bm{k}}\bm{w}^{\times}_{\bm{j}}\Phi_{\bm{j}}(f(\underline{\bm{A}}_{\bm{j}}))}\nonumber \\
&\geq&\left[\min\limits_{x \in \left(\bigcup\limits_{\bm{j}=\bm{1}}^{\bm{k}}\bm{w}^{\times}_{\bm{j}}\widetilde{\mbox{Poly}}_{f,\Psi_{\bm{j}},\mathscr{L}}\left(\bm{I}^{\bigtimes}\right)\right), \bm{x} \in \bigtimes\limits_{i=1}^n\left(\bigcup\limits_{\bm{j}=\bm{1}}^{\bm{k}}\bm{w}^{\times}_{\bm{j}}\widetilde{\mbox{Poly}}_{f,\Psi_{\bm{j}},\mathscr{L}}\left(\bm{I}^{\bigtimes}_i\right)\right)}xg^{-1}(\bm{x})\right] \nonumber \\
&&\times g\Bigg(\sum\limits_{\bm{j}=\bm{1}}^{\bm{k}}\bm{w}^{\times}_{\bm{j}}\bm{V}^{\ast}\mbox{Poly}_{f,\Psi_{\bm{j}},\mathscr{L}}(\underline{\bm{A}}_{j_1})\bm{V},\ldots,\sum\limits_{\bm{j}=\bm{1}}^{\bm{k}}\bm{w}^{\times}_{\bm{j}}\bm{V}^{\ast}\mbox{Poly}_{f,\Psi_{\bm{j}},\mathscr{L}}(\underline{\bm{A}}_{j_n})\bm{V}\Bigg).
\end{eqnarray}

(II) If we also assume that $g\Bigg(\sum\limits_{\bm{j}=\bm{1}}^{\bm{k}}\bm{w}^{\times}_{\bm{j}}\bm{V}^{\ast}\mbox{Poly}_{f,\Psi_{\bm{j}},\mathscr{U}}(\underline{\bm{A}}_{j_1})\bm{V},\ldots,\sum\limits_{\bm{j}=\bm{1}}^{\bm{k}}\bm{w}^{\times}_{\bm{j}}\bm{V}^{\ast}\mbox{Poly}_{f,\Psi_{\bm{j}},\mathscr{U}}(\underline{\bm{A}}_{j_n})\bm{V}\Bigg) < \bm{0}$ and $g\Bigg(\sum\limits_{\bm{j}=\bm{1}}^{\bm{k}}\bm{w}^{\times}_{\bm{j}}\bm{V}^{\ast}\mbox{Poly}_{f,\Psi_{\bm{j}},\mathscr{L}}(\underline{\bm{A}}_{j_1})\bm{V},\ldots,\sum\limits_{\bm{j}=\bm{1}}^{\bm{k}}\bm{w}^{\times}_{\bm{j}}\bm{V}^{\ast}\mbox{Poly}_{f,\Psi_{\bm{j}},\mathscr{L}}(\underline{\bm{A}}_{j_n})\bm{V}\Bigg) < \bm{0}$, then, we have the following upper bound for $\sum\limits_{\bm{j}=\bm{1}}^{\bm{k}}\bm{w}^{\times}_{\bm{j}}\Phi_{\bm{j}}(f(\underline{\bm{A}}_{\bm{j}}))$:
\begin{eqnarray}\label{eq1 UB: thm: 2.9 neg g}
\lefteqn{\sum\limits_{\bm{j}=\bm{1}}^{\bm{k}}\bm{w}^{\times}_{\bm{j}}\Phi_{\bm{j}}(f(\underline{\bm{A}}_{\bm{j}}))}\nonumber \\
&\leq&\left[\min\limits_{x \in \left(\bigcup\limits_{\bm{j}=\bm{1}}^{\bm{k}}\bm{w}^{\times}_{\bm{j}}\widetilde{\mbox{Poly}}_{f,\Psi_{\bm{j}},\mathscr{L}}\left(\bm{I}^{\bigtimes}\right)\right), \bm{x} \in \bigtimes\limits_{i=1}^n\left(\bigcup\limits_{\bm{j}=\bm{1}}^{\bm{k}}\bm{w}^{\times}_{\bm{j}}\widetilde{\mbox{Poly}}_{f,\Psi_{\bm{j}},\mathscr{L}}\left(\bm{I}^{\bigtimes}_i\right)\right)}xg^{-1}(\bm{x})\right] \nonumber \\
&&\times g\Bigg(\sum\limits_{\bm{j}=\bm{1}}^{\bm{k}}\bm{w}^{\times}_{\bm{j}}\bm{V}^{\ast}\mbox{Poly}_{f,\Psi_{\bm{j}},\mathscr{L}}(\underline{\bm{A}}_{j_1})\bm{V},\ldots,\sum\limits_{\bm{j}=\bm{1}}^{\bm{k}}\bm{w}^{\times}_{\bm{j}}\bm{V}^{\ast}\mbox{Poly}_{f,\Psi_{\bm{j}},\mathscr{L}}(\underline{\bm{A}}_{j_n})\bm{V}\Bigg);
\end{eqnarray}
and, the following lower bound for $\sum\limits_{\bm{j}=\bm{1}}^{\bm{k}}\bm{w}^{\times}_{\bm{j}}\Phi_{\bm{j}}(f(\underline{\bm{A}}_{\bm{j}}))$:
\begin{eqnarray}\label{eq1 LB: thm: 2.9 neg g}
\lefteqn{\sum\limits_{\bm{j}=\bm{1}}^{\bm{k}}\bm{w}^{\times}_{\bm{j}}\Phi_{\bm{j}}(f(\underline{\bm{A}}_{\bm{j}}))}\nonumber \\
&\geq&\left[\max\limits_{x \in \left(\bigcup\limits_{\bm{j}=\bm{1}}^{\bm{k}}\bm{w}^{\times}_{\bm{j}}\widetilde{\mbox{Poly}}_{f,\Psi_{\bm{j}},\mathscr{U}}\left(\bm{I}^{\bigtimes}\right)\right), \bm{x} \in \bigtimes\limits_{i=1}^n\left(\bigcup\limits_{\bm{j}=\bm{1}}^{\bm{k}}\bm{w}^{\times}_{\bm{j}}\widetilde{\mbox{Poly}}_{f,\Psi_{\bm{j}},\mathscr{U}}\left(\bm{I}^{\bigtimes}_i\right)\right)}xg^{-1}(\bm{x})\right] \nonumber \\
&&\times g\Bigg(\sum\limits_{\bm{j}=\bm{1}}^{\bm{k}}\bm{w}^{\times}_{\bm{j}}\bm{V}^{\ast}\mbox{Poly}_{f,\Psi_{\bm{j}},\mathscr{U}}(\underline{\bm{A}}_{j_1})\bm{V},\ldots,\sum\limits_{\bm{j}=\bm{1}}^{\bm{k}}\bm{w}^{\times}_{\bm{j}}\bm{V}^{\ast}\mbox{Poly}_{f,\Psi_{\bm{j}},\mathscr{U}}(\underline{\bm{A}}_{j_n})\bm{V}\Bigg).
\end{eqnarray}
\end{theorem}
\textbf{Proof:}
For part (I), we will apply $F(u,v)$ as 
\begin{eqnarray}\label{eq: F u v 2}
F(u,v)&=&v^{-1/2} u v^{-1/2},
\end{eqnarray}
to Eq.~\eqref{eq UB: thm:main 2.3} in Theorem~\ref{thm: main 2.3}, then, we will obtain 
\begin{eqnarray}\label{eq2 UB: thm: 2.9 pos g}
\lefteqn{\left(g\Bigg(\sum\limits_{\bm{j}=\bm{1}}^{\bm{k}}\bm{w}^{\times}_{\bm{j}}\bm{V}^{\ast}\mbox{Poly}_{f,\Psi_{\bm{j}},\mathscr{U}}(\underline{\bm{A}}_{j_1})\bm{V},\ldots,\sum\limits_{\bm{j}=\bm{1}}^{\bm{k}}\bm{w}^{\times}_{\bm{j}}\bm{V}^{\ast}\mbox{Poly}_{f,\Psi_{\bm{j}},\mathscr{U}}(\underline{\bm{A}}_{j_n})\bm{V}\Bigg)\right)^{-1/2}}\nonumber \\
&&\times \left(\sum\limits_{\bm{j}=\bm{1}}^{\bm{k}}\bm{w}^{\times}_{\bm{j}}\Phi_{\bm{j}}(f(\underline{\bm{A}}_{\bm{j}}))\right) \nonumber \\
&&\times \left(g\Bigg(\sum\limits_{\bm{j}=\bm{1}}^{\bm{k}}\bm{w}^{\times}_{\bm{j}}\bm{V}^{\ast}\mbox{Poly}_{f,\Psi_{\bm{j}},\mathscr{U}}(\underline{\bm{A}}_{j_1})\bm{V},\ldots,\sum\limits_{\bm{j}=\bm{1}}^{\bm{k}}\bm{w}^{\times}_{\bm{j}}\bm{V}^{\ast}\mbox{Poly}_{f,\Psi_{\bm{j}},\mathscr{U}}(\underline{\bm{A}}_{j_n})\bm{V}\Bigg)\right)^{-1/2}
\nonumber \\
&\leq&\left[\max\limits_{x \in \left(\bigcup\limits_{\bm{j}=\bm{1}}^{\bm{k}}\bm{w}^{\times}_{\bm{j}}\widetilde{\mbox{Poly}}_{f,\Psi_{\bm{j}},\mathscr{U}}\left(\bm{I}^{\bigtimes}\right)\right), \bm{x} \in \bigtimes\limits_{i=1}^n\left(\bigcup\limits_{\bm{j}=\bm{1}}^{\bm{k}}\bm{w}^{\times}_{\bm{j}}\widetilde{\mbox{Poly}}_{f,\Psi_{\bm{j}},\mathscr{U}}\left(\bm{I}^{\bigtimes}_i\right)\right)}xg^{-1}(\bm{x})\right]\bm{I}_{\mathfrak{K}}.
\end{eqnarray}
By multiplying $\left(g\Bigg(\sum\limits_{\bm{j}=\bm{1}}^{\bm{k}}\bm{w}^{\times}_{\bm{j}}\bm{V}^{\ast}\mbox{Poly}_{f,\Psi_{\bm{j}},\mathscr{U}}(\underline{\bm{A}}_{j_1})\bm{V},\ldots,\sum\limits_{\bm{j}=\bm{1}}^{\bm{k}}\bm{w}^{\times}_{\bm{j}}\bm{V}^{\ast}\mbox{Poly}_{f,\Psi_{\bm{j}},\mathscr{U}}(\underline{\bm{A}}_{j_n})\bm{V}\Bigg)\right)^{1/2}$ at both sides of Eq.~\eqref{eq2 UB: thm: 2.9 pos g}, we obtain the desired inequality provided by Eq.~\eqref{eq1 UB: thm: 2.9 pos g}. By applying $F(u,v)$ with Eq.~\eqref{eq: F u v 2} again to Eq.~\eqref{eq LB: thm:main 2.3} in Theorem~\ref{thm: main 2.3}, then, we will obtain 
\begin{eqnarray}\label{eq2 LB: thm: 2.9 pos g}
\lefteqn{\left(g\Bigg(\sum\limits_{\bm{j}=\bm{1}}^{\bm{k}}\bm{w}^{\times}_{\bm{j}}\bm{V}^{\ast}\mbox{Poly}_{f,\Psi_{\bm{j}},\mathscr{L}}(\underline{\bm{A}}_{j_1})\bm{V},\ldots,\sum\limits_{\bm{j}=\bm{1}}^{\bm{k}}\bm{w}^{\times}_{\bm{j}}\bm{V}^{\ast}\mbox{Poly}_{f,\Psi_{\bm{j}},\mathscr{L}}(\underline{\bm{A}}_{j_n})\bm{V}\Bigg)\right)^{-1/2}}\nonumber \\
&&\times \left(\sum\limits_{\bm{j}=\bm{1}}^{\bm{k}}\bm{w}^{\times}_{\bm{j}}\Phi_{\bm{j}}(f(\underline{\bm{A}}_{\bm{j}}))\right) \nonumber \\
&&\times \left(g\Bigg(\sum\limits_{\bm{j}=\bm{1}}^{\bm{k}}\bm{w}^{\times}_{\bm{j}}\bm{V}^{\ast}\mbox{Poly}_{f,\Psi_{\bm{j}},\mathscr{L}}(\underline{\bm{A}}_{j_1})\bm{V},\ldots,\sum\limits_{\bm{j}=\bm{1}}^{\bm{k}}\bm{w}^{\times}_{\bm{j}}\bm{V}^{\ast}\mbox{Poly}_{f,\Psi_{\bm{j}},\mathscr{L}}(\underline{\bm{A}}_{j_n})\bm{V}\Bigg)\right)^{-1/2}
\nonumber \\
&\geq&\left[\min\limits_{x \in \left(\bigcup\limits_{\bm{j}=\bm{1}}^{\bm{k}}\bm{w}^{\times}_{\bm{j}}\widetilde{\mbox{Poly}}_{f,\Psi_{\bm{j}},\mathscr{L}}\left(\bm{I}^{\bigtimes}\right)\right), \bm{x} \in \bigtimes\limits_{i=1}^n\left(\bigcup\limits_{\bm{j}=\bm{1}}^{\bm{k}}\bm{w}^{\times}_{\bm{j}}\widetilde{\mbox{Poly}}_{f,\Psi_{\bm{j}},\mathscr{L}}\left(\bm{I}^{\bigtimes}_i\right)\right)}xg^{-1}(\bm{x})\right]\bm{I}_{\mathfrak{K}}.
\end{eqnarray}
By multiplying $\left(g\Bigg(\sum\limits_{\bm{j}=\bm{1}}^{\bm{k}}\bm{w}^{\times}_{\bm{j}}\bm{V}^{\ast}\mbox{Poly}_{f,\Psi_{\bm{j}},\mathscr{L}}(\underline{\bm{A}}_{j_1})\bm{V},\ldots,\sum\limits_{\bm{j}=\bm{1}}^{\bm{k}}\bm{w}^{\times}_{\bm{j}}\bm{V}^{\ast}\mbox{Poly}_{f,\Psi_{\bm{j}},\mathscr{L}}(\underline{\bm{A}}_{j_n})\bm{V}\Bigg)\right)^{1/2}$ at both sides of Eq.~\eqref{eq2 LB: thm: 2.9 pos g}, we obtain the desired inequality provided by Eq.~\eqref{eq1 LB: thm: 2.9 pos g}. 

The proof of Part (II) is immediate obtained by setting $g(x)$ as $-g(x)$ in Eq.~\eqref{eq2 UB: thm: 2.9 pos g} and Eq.~\eqref{eq2 LB: thm: 2.9 pos g}. 
$\hfill \Box$

Next Corollary~\ref{cor: 2.11-14} is obtained by applying Theorem~\ref{thm: 2.9} to special types of function $g$. 

\begin{corollary}\label{cor: 2.11-14}
Let $\bm{A}_{j_i}$ be self-adjoint operators with $\Lambda(\bm{A}_{j_i}) \in [m_i, M_i]$ for real scalars $m_i <  M_i$. The mappings $\Phi_{j_1,\ldots,j_n}: \mathscr{B}(\mathfrak{H}) \rightarrow \mathscr{B}(\mathfrak{K})$ are defined by Eq.~\eqref{eq: new phi def}, where $j_i=1,2,\ldots,k_i$ for $i=1,2,\ldots,n$. We have $n$ probability vectors $\bm{w}_i =[w_{i,1},w_{i,2},\cdots, w_{i,k_i}]$ with the dimension $k_i$ for $i=1,2,\ldots,n$, i.e., $\sum\limits_{\ell=1}^{k_i}w_{i,\ell} = 1$. Let $f(\bm{x})$ be any real valued continuous functions with $n$ variables defined on the range $\bigtimes\limits_{i=1}^n [m_i, M_i] \in \mathbb{R}^n$, where $\bigtimes$ is the Cartesian product. Besides, given any $\epsilon>0$, we assume that the function $f(\bm{x})$ satisfies the following:
\begin{eqnarray}\label{eq: lower and upper Psi cor: 2.11-14}
0 &\leq& \Psi_{\mathscr{U}}(\bm{x}) - f(\bm{x})~~\leq~~\epsilon, \nonumber \\
0 &\leq& f(\bm{x})-\Psi_{\mathscr{L}}(\bm{x})~~\leq~~\epsilon,
\end{eqnarray}
for $\bm{x}\in \bm{I}^{\bigtimes}$. The function $g(\bm{x})$ is also a real valued continuous function with $n$ variables defined on the range $\bigtimes\limits_{i=1}^n [m_i, M_i]$.

(I) If we have $g(\bm{x})=\left(\sum\limits_{i=1}^n \beta_i x_i \right)^q$ with $\beta_i \geq 0$ and $q \in \mathbb{R}$, and we also assume that \\
$\left(\sum\limits_{\bm{j}=\bm{1}}^{\bm{k}}\bm{w}^{\times}_{\bm{j}}\bm{V}^{\ast}\mbox{Poly}_{f,\Psi_{\bm{j}},\mathscr{U}}(\underline{\bm{A}}_{j_i})\bm{V}\right) > \bm{0}$ and $\left(\sum\limits_{\bm{j}=\bm{1}}^{\bm{k}}\bm{w}^{\times}_{\bm{j}}\bm{V}^{\ast}\mbox{Poly}_{f,\Psi_{\bm{j}},\mathscr{L}}(\underline{\bm{A}}_{j_i})\bm{V}\right) > \bm{0}$, then, we have the following upper bound for $\sum\limits_{\bm{j}=\bm{1}}^{\bm{k}}\bm{w}^{\times}_{\bm{j}}\Phi_{\bm{j}}(f(\underline{\bm{A}}_{\bm{j}}))$:
\begin{eqnarray}\label{eq x q UB: cor: 2.11-14}
\sum\limits_{\bm{j}=\bm{1}}^{\bm{k}}\bm{w}^{\times}_{\bm{j}}\Phi_{\bm{j}}(f(\underline{\bm{A}}_{\bm{j}}))&\leq&\left[\max\limits_{x \in \left(\bigcup\limits_{\bm{j}=\bm{1}}^{\bm{k}}\bm{w}^{\times}_{\bm{j}}\widetilde{\mbox{Poly}}_{f,\Psi_{\bm{j}},\mathscr{U}}\left(\bm{I}^{\bigtimes}\right)\right),x_i \in \left(\bigcup\limits_{\bm{j}=\bm{1}}^{\bm{k}}\bm{w}^{\times}_{\bm{j}}\widetilde{\mbox{Poly}}_{f,\Psi_{\bm{j}},\mathscr{U}}\left(\bm{I}_i^{\bigtimes}\right)\right)}\frac{x}{\left(\sum\limits_{i=1}^n \beta_i x_i\right)^q}\right]\nonumber\\
&&\times \left(\sum\limits_{i=1}^n \beta_i \left(\sum\limits_{\bm{j}=\bm{1}}^{\bm{k}}\bm{w}^{\times}_{\bm{j}}\bm{V}^{\ast}\mbox{Poly}_{f,\Psi_{\bm{j}},\mathscr{U}}(\underline{\bm{A}}_{j_i})\bm{V}\right)\right)^q;
\end{eqnarray}
and, the following lower bound for $\sum\limits_{\bm{j}=\bm{1}}^{\bm{k}}\bm{w}^{\times}_{\bm{j}}\Phi_{\bm{j}}(f(\underline{\bm{A}}_{\bm{j}}))$:
\begin{eqnarray}\label{eq x q LB: cor: 2.11-14}
\sum\limits_{\bm{j}=\bm{1}}^{\bm{k}}\bm{w}^{\times}_{\bm{j}}\Phi_{\bm{j}}(f(\underline{\bm{A}}_{\bm{j}}))&\geq&\left[\min\limits_{x \in \left(\bigcup\limits_{\bm{j}=\bm{1}}^{\bm{k}}\bm{w}^{\times}_{\bm{j}}\widetilde{\mbox{Poly}}_{f,\Psi_{\bm{j}},\mathscr{L}}\left(\bm{I}^{\bigtimes}\right)\right),x_i \in \left(\bigcup\limits_{\bm{j}=\bm{1}}^{\bm{k}}\bm{w}^{\times}_{\bm{j}}\widetilde{\mbox{Poly}}_{f,\Psi_{\bm{j}},\mathscr{L}}\left(\bm{I}_i^{\bigtimes}\right)\right)}\frac{x}{\left(\sum\limits_{i=1}^n \beta_i x_i\right)^q}\right]\nonumber\\
&&\times \left(\sum\limits_{i=1}^n \beta_i \left(\sum\limits_{\bm{j}=\bm{1}}^{\bm{k}}\bm{w}^{\times}_{\bm{j}}\bm{V}^{\ast}\mbox{Poly}_{f,\Psi_{\bm{j}},\mathscr{L}}(\underline{\bm{A}}_{j_i})\bm{V}\right)\right)^q.
\end{eqnarray}

(II) If we have $g(\bm{x})=\log\left(\sum\limits_{i=1}^n \beta_i x_i \right)$ with $\sum\limits_{i=1}^n \beta_i x_i > 0$, and we also assume that \\
$\log\left(\sum\limits_{i=1}^n \beta_i \left(\sum\limits_{\bm{j}=\bm{1}}^{\bm{k}}\bm{w}^{\times}_{\bm{j}}\bm{V}^{\ast}\mbox{Poly}_{f,\Psi_{\bm{j}},\mathscr{U}}(\underline{\bm{A}}_{j_i})\bm{V}\right)\right)> \bm{0}$ and $\log\left(\sum\limits_{i=1}^n \beta_i \left(\sum\limits_{\bm{j}=\bm{1}}^{\bm{k}}\bm{w}^{\times}_{\bm{j}}\bm{V}^{\ast}\mbox{Poly}_{f,\Psi_{\bm{j}},\mathscr{L}}(\underline{\bm{A}}_{j_i})\bm{V}\right)\right)> \bm{0}$, then, we have the following upper bound for $\sum\limits_{\bm{j}=\bm{1}}^{\bm{k}}\bm{w}^{\times}_{\bm{j}}\Phi_{\bm{j}}(f(\underline{\bm{A}}_{\bm{j}}))$:
\begin{eqnarray}\label{eq log x l0 UB: cor: 2.11-14}
\sum\limits_{\bm{j}=\bm{1}}^{\bm{k}}\bm{w}^{\times}_{\bm{j}}\Phi_{\bm{j}}(f(\underline{\bm{A}}_{\bm{j}}))&\leq&\left[\max\limits_{x \in \left(\bigcup\limits_{\bm{j}=\bm{1}}^{\bm{k}}\bm{w}^{\times}_{\bm{j}}\widetilde{\mbox{Poly}}_{f,\Psi_{\bm{j}},\mathscr{U}}\left(\bm{I}^{\bigtimes}\right)\right),x_i \in \left(\bigcup\limits_{\bm{j}=\bm{1}}^{\bm{k}}\bm{w}^{\times}_{\bm{j}}\widetilde{\mbox{Poly}}_{f,\Psi_{\bm{j}},\mathscr{U}}\left(\bm{I}_i^{\bigtimes}\right)\right)}\frac{x}{\log\left(\sum\limits_{i=1}^n \beta_i x_i\right)}\right]\nonumber\\
&&\times \log\left(\sum\limits_{i=1}^n \beta_i \left(\sum\limits_{\bm{j}=\bm{1}}^{\bm{k}}\bm{w}^{\times}_{\bm{j}}\bm{V}^{\ast}\mbox{Poly}_{f,\Psi_{\bm{j}},\mathscr{U}}(\underline{\bm{A}}_{j_i})\bm{V}\right)\right);
\end{eqnarray}
and, the following lower bound for $\sum\limits_{\bm{j}=\bm{1}}^{\bm{k}}\bm{w}^{\times}_{\bm{j}}\Phi_{\bm{j}}(f(\underline{\bm{A}}_{\bm{j}}))$:
\begin{eqnarray}\label{eq log x l0 LB: cor: 2.11-14}
\sum\limits_{\bm{j}=\bm{1}}^{\bm{k}}\bm{w}^{\times}_{\bm{j}}\Phi_{\bm{j}}(f(\underline{\bm{A}}_{\bm{j}}))&\geq&\left[\min\limits_{x \in \left(\bigcup\limits_{\bm{j}=\bm{1}}^{\bm{k}}\bm{w}^{\times}_{\bm{j}}\widetilde{\mbox{Poly}}_{f,\Psi_{\bm{j}},\mathscr{L}}\left(\bm{I}^{\bigtimes}\right)\right),x_i \in \left(\bigcup\limits_{\bm{j}=\bm{1}}^{\bm{k}}\bm{w}^{\times}_{\bm{j}}\widetilde{\mbox{Poly}}_{f,\Psi_{\bm{j}},\mathscr{L}}\left(\bm{I}_i^{\bigtimes}\right)\right)}\frac{x}{\log\left(\sum\limits_{i=1}^n \beta_i x_i\right)}\right]\nonumber\\
&&\times \log\left(\sum\limits_{i=1}^n \beta_i \left(\sum\limits_{\bm{j}=\bm{1}}^{\bm{k}}\bm{w}^{\times}_{\bm{j}}\bm{V}^{\ast}\mbox{Poly}_{f,\Psi_{\bm{j}},\mathscr{L}}(\underline{\bm{A}}_{j_i})\bm{V}\right)\right).
\end{eqnarray}

(II') If we have $g(\bm{x})=\log\left(\sum\limits_{i=1}^n \beta_i x_i \right)$ with $\sum\limits_{i=1}^n \beta_i x_i > 0$, and we also assume that \\
$\log\left(\sum\limits_{i=1}^n \beta_i \left(\sum\limits_{\bm{j}=\bm{1}}^{\bm{k}}\bm{w}^{\times}_{\bm{j}}\bm{V}^{\ast}\mbox{Poly}_{f,\Psi_{\bm{j}},\mathscr{U}}(\underline{\bm{A}}_{j_i})\bm{V}\right)\right) < \bm{0}$ and $\log\left(\sum\limits_{i=1}^n \beta_i \left(\sum\limits_{\bm{j}=\bm{1}}^{\bm{k}}\bm{w}^{\times}_{\bm{j}}\bm{V}^{\ast}\mbox{Poly}_{f,\Psi_{\bm{j}},\mathscr{L}}(\underline{\bm{A}}_{j_i})\bm{V}\right)\right) < \bm{0}$, then, we have the following upper bound for $\sum\limits_{\bm{j}=\bm{1}}^{\bm{k}}\bm{w}^{\times}_{\bm{j}}\Phi_{\bm{j}}(f(\underline{\bm{A}}_{\bm{j}}))$:
\begin{eqnarray}\label{eq log x s0 UB: cor: 2.11-14}
\sum\limits_{\bm{j}=\bm{1}}^{\bm{k}}\bm{w}^{\times}_{\bm{j}}\Phi_{\bm{j}}(f(\underline{\bm{A}}_{\bm{j}}))&\leq&\left[\min\limits_{x \in \left(\bigcup\limits_{\bm{j}=\bm{1}}^{\bm{k}}\bm{w}^{\times}_{\bm{j}}\widetilde{\mbox{Poly}}_{f,\Psi_{\bm{j}},\mathscr{L}}\left(\bm{I}^{\bigtimes}\right)\right),x_i \in \left(\bigcup\limits_{\bm{j}=\bm{1}}^{\bm{k}}\bm{w}^{\times}_{\bm{j}}\widetilde{\mbox{Poly}}_{f,\Psi_{\bm{j}},\mathscr{L}}\left(\bm{I}_i^{\bigtimes}\right)\right)}\frac{x}{\log\left(\sum\limits_{i=1}^n \beta_i x_i\right)}\right]\nonumber\\
&&\times \log\left(\sum\limits_{i=1}^n \beta_i \left(\sum\limits_{\bm{j}=\bm{1}}^{\bm{k}}\bm{w}^{\times}_{\bm{j}}\bm{V}^{\ast}\mbox{Poly}_{f,\Psi_{\bm{j}},\mathscr{L}}(\underline{\bm{A}}_{j_i})\bm{V}\right)\right);
\end{eqnarray}
and, the following lower bound for $\sum\limits_{\bm{j}=\bm{1}}^{\bm{k}}\bm{w}^{\times}_{\bm{j}}\Phi_{\bm{j}}(f(\underline{\bm{A}}_{\bm{j}}))$:
\begin{eqnarray}\label{eq log x s0 LB: cor: 2.11-14}
\sum\limits_{\bm{j}=\bm{1}}^{\bm{k}}\bm{w}^{\times}_{\bm{j}}\Phi_{\bm{j}}(f(\underline{\bm{A}}_{\bm{j}}))&\geq&\left[\max\limits_{x \in \left(\bigcup\limits_{\bm{j}=\bm{1}}^{\bm{k}}\bm{w}^{\times}_{\bm{j}}\widetilde{\mbox{Poly}}_{f,\Psi_{\bm{j}},\mathscr{U}}\left(\bm{I}^{\bigtimes}\right)\right),x_i \in \left(\bigcup\limits_{\bm{j}=\bm{1}}^{\bm{k}}\bm{w}^{\times}_{\bm{j}}\widetilde{\mbox{Poly}}_{f,\Psi_{\bm{j}},\mathscr{U}}\left(\bm{I}_i^{\bigtimes}\right)\right)}\frac{x}{\log\left(\sum\limits_{i=1}^n \beta_i x_i\right)}\right]\nonumber\\
&&\times \log\left(\sum\limits_{i=1}^n \beta_i \left(\sum\limits_{\bm{j}=\bm{1}}^{\bm{k}}\bm{w}^{\times}_{\bm{j}}\bm{V}^{\ast}\mbox{Poly}_{f,\Psi_{\bm{j}},\mathscr{U}}(\underline{\bm{A}}_{j_i})\bm{V}\right)\right);
\end{eqnarray}

(III) If we have $g(\bm{x})=\exp\left(\sum\limits_{i=1}^n \beta_i x_i \right)$, then, we have the following upper bound for $\sum\limits_{\bm{j}=\bm{1}}^{\bm{k}}\bm{w}^{\times}_{\bm{j}}\Phi_{\bm{j}}(f(\underline{\bm{A}}_{\bm{j}}))$: 
\begin{eqnarray}\label{eq exp x UB: cor: 2.11-14}
\sum\limits_{\bm{j}=\bm{1}}^{\bm{k}}\bm{w}^{\times}_{\bm{j}}\Phi_{\bm{j}}(f(\underline{\bm{A}}_{\bm{j}}))&\leq&\left[\max\limits_{x \in \left(\bigcup\limits_{\bm{j}=\bm{1}}^{\bm{k}}\bm{w}^{\times}_{\bm{j}}\widetilde{\mbox{Poly}}_{f,\Psi_{\bm{j}},\mathscr{U}}\left(\bm{I}^{\bigtimes}\right)\right),x_i \in \left(\bigcup\limits_{\bm{j}=\bm{1}}^{\bm{k}}\bm{w}^{\times}_{\bm{j}}\widetilde{\mbox{Poly}}_{f,\Psi_{\bm{j}},\mathscr{U}}\left(\bm{I}_i^{\bigtimes}\right)\right)}\frac{x}{\exp\left(\sum\limits_{i=1}^n \beta_i x_i\right)}\right]\nonumber\\
&&\times \exp\left(\sum\limits_{i=1}^n \beta_i \left(\sum\limits_{\bm{j}=\bm{1}}^{\bm{k}}\bm{w}^{\times}_{\bm{j}}\bm{V}^{\ast}\mbox{Poly}_{f,\Psi_{\bm{j}},\mathscr{U}}(\underline{\bm{A}}_{j_i})\bm{V}\right)\right);
\end{eqnarray}
and, the following lower bound for $\sum\limits_{\bm{j}=\bm{1}}^{\bm{k}}\bm{w}^{\times}_{\bm{j}}\Phi_{\bm{j}}(f(\underline{\bm{A}}_{\bm{j}}))$: 
\begin{eqnarray}\label{eq exp x LB: cor: 2.11-14}
\sum\limits_{\bm{j}=\bm{1}}^{\bm{k}}\bm{w}^{\times}_{\bm{j}}\Phi_{\bm{j}}(f(\underline{\bm{A}}_{\bm{j}}))&\geq&\left[\min\limits_{x \in \left(\bigcup\limits_{\bm{j}=\bm{1}}^{\bm{k}}\bm{w}^{\times}_{\bm{j}}\widetilde{\mbox{Poly}}_{f,\Psi_{\bm{j}},\mathscr{L}}\left(\bm{I}^{\bigtimes}\right)\right),x_i \in \left(\bigcup\limits_{\bm{j}=\bm{1}}^{\bm{k}}\bm{w}^{\times}_{\bm{j}}\widetilde{\mbox{Poly}}_{f,\Psi_{\bm{j}},\mathscr{L}}\left(\bm{I}_i^{\bigtimes}\right)\right)}\frac{x}{\exp\left(\sum\limits_{i=1}^n \beta_i x_i\right)}\right]\nonumber\\
&&\times \exp\left(\sum\limits_{i=1}^n \beta_i \left(\sum\limits_{\bm{j}=\bm{1}}^{\bm{k}}\bm{w}^{\times}_{\bm{j}}\bm{V}^{\ast}\mbox{Poly}_{f,\Psi_{\bm{j}},\mathscr{L}}(\underline{\bm{A}}_{j_i})\bm{V}\right)\right).
\end{eqnarray}
\end{corollary}
\textbf{Proof:}
Part (I) of this Corollary is proved by applying Theorem~\ref{thm: 2.9} Part (I) with the function $g$ as $g(\bm{x})=\left(\sum\limits_{i=1}^n \beta_i x_i\right)^q$. Part (II) of this Corollary is proved by applying Theorem~\ref{thm: 2.9} Part (I) with the function $g$ as $g(\bm{x})=\log\left(\sum\limits_{i=1}^n \beta_i x_i\right)$, where $\log\left(\sum\limits_{i=1}^n \beta_i x_i\right) > 0$. Part (II') of this Corollary is proved by applying Theorem~\ref{thm: 2.9} Part (II) with the function $g$ as $g(\bm{x})=\log\left(\sum\limits_{i=1}^n \beta_i x_i\right)$, where $\log\left(\sum\limits_{i=1}^n \beta_i x_i\right) < 0$. Finally, Part (III) of this Corollary is proved by applying Theorem~\ref{thm: 2.9} Part (I) with the function $g$ as $g(\bm{x})=\exp\left(\sum\limits_{i=1}^n \beta_i x_i\right)$.      
$\hfill \Box$

\section{Multivariate Hypercomplex Functions Inequalities: Differences Kind}\label{sec: Multivariate Hypercomplex Functions Inequalities: Difference Kind}

In this section, we will derive the lower and upper bounds for $\sum\limits_{\bm{j}=\bm{1}}^{\bm{k}}\bm{w}^{\times}_{\bm{j}}\Phi_{\bm{j}}(f(\underline{\bm{A}}_{\bm{j}}))$ in terms of difference criteria related to the function $g$.
\begin{theorem}\label{thm: cor 2.15}
Let $\bm{A}_{j_i}$ be self-adjoint operators with $\Lambda(\bm{A}_{j_i}) \in [m_i, M_i]$ for real scalars $m_i <  M_i$. The mappings $\Phi_{j_1,\ldots,j_n}: \mathscr{B}(\mathfrak{H}) \rightarrow \mathscr{B}(\mathfrak{K})$ are defined by Eq.~\eqref{eq: new phi def}, where $j_i=1,2,\ldots,k_i$ for $i=1,2,\ldots,n$. We have $n$ probability vectors $\bm{w}_i =[w_{i,1},w_{i,2},\cdots, w_{i,k_i}]$ with the dimension $k_i$ for $i=1,2,\ldots,n$, i.e., $\sum\limits_{\ell=1}^{k_i}w_{i,\ell} = 1$. Let $f(\bm{x})$ be any real valued continuous functions with $n$ variables defined on the range $\bigtimes\limits_{i=1}^n [m_i, M_i] \in \mathbb{R}^n$, where $\bigtimes$ is the Cartesian product. Besides, given any $\epsilon>0$, we assume that the function $f(\bm{x})$ satisfies the following:
\begin{eqnarray}\label{eq: lower and upper Psi thm: cor 2.15}
0 &\leq& \Psi_{\mathscr{U}}(\bm{x}) - f(\bm{x})~~\leq~~\epsilon, \nonumber \\
0 &\leq& f(\bm{x})-\Psi_{\mathscr{L}}(\bm{x})~~\leq~~\epsilon,
\end{eqnarray}
for $\bm{x}\in \bm{I}^{\bigtimes}$. The function $g(\bm{x})$ is also a real valued continuous function with $n$ variables defined on the range $\bigtimes\limits_{i=1}^n [m_i, M_i]$. We also have a real valued function $F(u,v)$ defined by Eq.~\eqref{eq: F u v } with support domain on $U \times V$ such that $f(\bm{I}^{\bigtimes}) \subset U$, and \\
$g\left(\bigtimes\limits_{i=1}^n\left(\bigcup\limits_{\bm{j}=\bm{1}}^{\bm{k}}\bm{w}^{\times}_{\bm{j}}\widetilde{\mbox{Poly}}_{f,\Psi_{\bm{j}},\mathscr{U}}\left(\bm{I}^{\bigtimes}_i\right)\right)\bigcup\bigtimes\limits_{i=1}^n\left(\bigcup\limits_{\bm{j}=\bm{1}}^{\bm{k}}\bm{w}^{\times}_{\bm{j}}\widetilde{\mbox{Poly}}_{f,\Psi_{\bm{j}},\mathscr{L}}\left(\bm{I}^{\bigtimes}_i\right)\right)\right) \subset V$. 

Then, we have the following upper bound:
\begin{eqnarray}\label{eq UB: cor 2.15}
\lefteqn{\sum\limits_{\bm{j}=\bm{1}}^{\bm{k}}\bm{w}^{\times}_{\bm{j}}\Phi_{\bm{j}}(f(\underline{\bm{A}}_{\bm{j}})) - g\Bigg(\sum\limits_{\bm{j}=\bm{1}}^{\bm{k}}\bm{w}^{\times}_{\bm{j}}\bm{V}^{\ast}\mbox{Poly}_{f,\Psi_{\bm{j}},\mathscr{U}}(\underline{\bm{A}}_{j_1})\bm{V},\ldots,\sum\limits_{\bm{j}=\bm{1}}^{\bm{k}}\bm{w}^{\times}_{\bm{j}}\bm{V}^{\ast}\mbox{Poly}_{f,\Psi_{\bm{j}},\mathscr{U}}(\underline{\bm{A}}_{j_n})\bm{V}\Bigg)}\nonumber \\
&\leq& 
\max\limits_{x \in \left(\bigcup\limits_{\bm{j}=\bm{1}}^{\bm{k}}\bm{w}^{\times}_{\bm{j}}\widetilde{\mbox{Poly}}_{f,\Psi_{\bm{j}},\mathscr{U}}\left(\bm{I}^{\bigtimes}\right)\right), \bm{x} \in \bigtimes\limits_{i=1}^n\left(\bigcup\limits_{\bm{j}=\bm{1}}^{\bm{k}}\bm{w}^{\times}_{\bm{j}}\widetilde{\mbox{Poly}}_{f,\Psi_{\bm{j}},\mathscr{U}}\left(\bm{I}^{\bigtimes}_i\right)\right)}(x - g(\bm{x}))\bm{I}_{\mathfrak{K}}.
\end{eqnarray}
Similarly, we also have the following lower bound:
\begin{eqnarray}\label{eq LB: cor 2.15}
\lefteqn{\sum\limits_{\bm{j}=\bm{1}}^{\bm{k}}\bm{w}^{\times}_{\bm{j}}\Phi_{\bm{j}}(f(\underline{\bm{A}}_{\bm{j}})) - g\Bigg(\sum\limits_{\bm{j}=\bm{1}}^{\bm{k}}\bm{w}^{\times}_{\bm{j}}\bm{V}^{\ast}\mbox{Poly}_{f,\Psi_{\bm{j}},\mathscr{L}}(\underline{\bm{A}}_{j_1})\bm{V},\ldots,\sum\limits_{\bm{j}=\bm{1}}^{\bm{k}}\bm{w}^{\times}_{\bm{j}}\bm{V}^{\ast}\mbox{Poly}_{f,\Psi_{\bm{j}},\mathscr{L}}(\underline{\bm{A}}_{j_n})\bm{V}\Bigg)}\nonumber \\
&\geq&
\min\limits_{x \in \left(\bigcup\limits_{\bm{j}=\bm{1}}^{\bm{k}}\bm{w}^{\times}_{\bm{j}}\widetilde{\mbox{Poly}}_{f,\Psi_{\bm{j}},\mathscr{L}}\left(\bm{I}^{\bigtimes}\right)\right), \bm{x} \in \bigtimes\limits_{i=1}^n\left(\bigcup\limits_{\bm{j}=\bm{1}}^{\bm{k}}\bm{w}^{\times}_{\bm{j}}\widetilde{\mbox{Poly}}_{f,\Psi_{\bm{j}},\mathscr{L}}\left(\bm{I}^{\bigtimes}_i\right)\right)}(x - g(\bm{x}))\bm{I}_{\mathfrak{K}}.
\end{eqnarray}
\end{theorem}
\textbf{Proof:}
The upper bound of this theorem is proved by setting $\alpha=1$ in Eq.~\eqref{eq UB: thm:main 2.4} from Theorem~\ref{thm:main 2.4} and rearrangement of the term $g\Bigg(\sum\limits_{\bm{j}=\bm{1}}^{\bm{k}}\bm{w}^{\times}_{\bm{j}}\bm{V}^{\ast}\mbox{Poly}_{f,\Psi_{\bm{j}},\mathscr{U}}(\underline{\bm{A}}_{j_1})\bm{V},\ldots,\sum\limits_{\bm{j}=\bm{1}}^{\bm{k}}\bm{w}^{\times}_{\bm{j}}\bm{V}^{\ast}\mbox{Poly}_{f,\Psi_{\bm{j}},\mathscr{U}}(\underline{\bm{A}}_{j_n})\bm{V}\Bigg)$ to obtain Eq.~\eqref{eq UB: cor 2.15}. 

Similarly, the lower bound of this theorem is proved by setting $\alpha=1$ in Eq.~\eqref{eq LB: thm:main 2.4} from Theorem~\ref{thm:main 2.4} and rearrangement of the term $g\Bigg(\sum\limits_{\bm{j}=\bm{1}}^{\bm{k}}\bm{w}^{\times}_{\bm{j}}\bm{V}^{\ast}\mbox{Poly}_{f,\Psi_{\bm{j}},\mathscr{L}}(\underline{\bm{A}}_{j_1})\bm{V},\ldots,\sum\limits_{\bm{j}=\bm{1}}^{\bm{k}}\bm{w}^{\times}_{\bm{j}}\bm{V}^{\ast}\mbox{Poly}_{f,\Psi_{\bm{j}},\mathscr{L}}(\underline{\bm{A}}_{j_n})\bm{V}\Bigg)$ to obtain Eq.~\eqref{eq LB: cor 2.15}. 
$\hfill \Box$

Next Corollary~\ref{cor: 2.17 ext} is obtained by applying Theorem~\ref{thm: cor 2.15} to special types of function $g$. 
\begin{corollary}\label{cor: 2.17 ext}
Let $\bm{A}_{j_i}$ be self-adjoint operators with $\Lambda(\bm{A}_{j_i}) \in [m_i, M_i]$ for real scalars $m_i <  M_i$. The mappings $\Phi_{j_1,\ldots,j_n}: \mathscr{B}(\mathfrak{H}) \rightarrow \mathscr{B}(\mathfrak{K})$ are defined by Eq.~\eqref{eq: new phi def}, where $j_i=1,2,\ldots,k_i$ for $i=1,2,\ldots,n$. We have $n$ probability vectors $\bm{w}_i =[w_{i,1},w_{i,2},\cdots, w_{i,k_i}]$ with the dimension $k_i$ for $i=1,2,\ldots,n$, i.e., $\sum\limits_{\ell=1}^{k_i}w_{i,\ell} = 1$. Let $f(\bm{x})$ be any real valued continuous functions with $n$ variables defined on the range $\bigtimes\limits_{i=1}^n [m_i, M_i] \in \mathbb{R}^n$, where $\bigtimes$ is the Cartesian product. Besides, given any $\epsilon>0$, we assume that the function $f(\bm{x})$ satisfies the following:
\begin{eqnarray}\label{eq: lower and upper Psi thm: cor: 2.17 ext}
0 &\leq& \Psi_{\mathscr{U}}(\bm{x}) - f(\bm{x})~~\leq~~\epsilon, \nonumber \\
0 &\leq& f(\bm{x})-\Psi_{\mathscr{L}}(\bm{x})~~\leq~~\epsilon,
\end{eqnarray}
for $\bm{x}\in \bm{I}^{\bigtimes}$.

(I) If we also assume that $g(\bm{x})=\left(\sum\limits_{i=1}^n \beta_i x_i\right)^q$, where $q \in \mathbb{R}$, then, we have the following upper bound for $\sum\limits_{\bm{j}=\bm{1}}^{\bm{k}}\bm{w}^{\times}_{\bm{j}}\Phi_{\bm{j}}(f(\underline{\bm{A}}_{\bm{j}}))$:
\begin{eqnarray}\label{eq x q UB: cor: 2.17 diff}
\lefteqn{\sum\limits_{\bm{j}=\bm{1}}^{\bm{k}}\bm{w}^{\times}_{\bm{j}}\Phi_{\bm{j}}(f(\underline{\bm{A}}_{\bm{j}})) - \left(\sum\limits_{i=1}^n \beta_i \left(\sum\limits_{\bm{j}=\bm{1}}^{\bm{k}}\bm{w}^{\times}_{\bm{j}}\bm{V}^{\ast}\mbox{Poly}_{f,\Psi_{\bm{j}},\mathscr{U}}(\underline{\bm{A}}_{j_i})\bm{V}\right)\right)^q}\nonumber \\
&\leq& 
\max\limits_{x \in \left(\bigcup\limits_{\bm{j}=\bm{1}}^{\bm{k}}\bm{w}^{\times}_{\bm{j}}\widetilde{\mbox{Poly}}_{f,\Psi_{\bm{j}},\mathscr{U}}\left(\bm{I}^{\bigtimes}\right)\right),x_i \in \left(\bigcup\limits_{\bm{j}=\bm{1}}^{\bm{k}}\bm{w}^{\times}_{\bm{j}}\widetilde{\mbox{Poly}}_{f,\Psi_{\bm{j}},\mathscr{U}}\left(\bm{I}_i^{\bigtimes}\right)\right)}\left(x - \left(\sum\limits_{i=1}^n \beta_i x_i\right)^q\right)\bm{I}_{\mathfrak{K}};
\end{eqnarray}
and, the following lower bound for $\sum\limits_{\bm{j}=\bm{1}}^{\bm{k}}\bm{w}^{\times}_{\bm{j}}\Phi_{\bm{j}}(f(\underline{\bm{A}}_{\bm{j}}))$:
\begin{eqnarray}\label{eq x q LB: cor: 2.17 diff}
\lefteqn{\sum\limits_{\bm{j}=\bm{1}}^{\bm{k}}\bm{w}^{\times}_{\bm{j}}\Phi_{\bm{j}}(f(\underline{\bm{A}}_{\bm{j}})) - \left(\sum\limits_{i=1}^n \beta_i \left(\sum\limits_{\bm{j}=\bm{1}}^{\bm{k}}\bm{w}^{\times}_{\bm{j}}\bm{V}^{\ast}\mbox{Poly}_{f,\Psi_{\bm{j}},\mathscr{L}}(\underline{\bm{A}}_{j_i})\bm{V}\right)\right)^q}\nonumber \\
&\geq& 
\min\limits_{x \in \left(\bigcup\limits_{\bm{j}=\bm{1}}^{\bm{k}}\bm{w}^{\times}_{\bm{j}}\widetilde{\mbox{Poly}}_{f,\Psi_{\bm{j}},\mathscr{L}}\left(\bm{I}^{\bigtimes}\right)\right),x_i \in \left(\bigcup\limits_{\bm{j}=\bm{1}}^{\bm{k}}\bm{w}^{\times}_{\bm{j}}\widetilde{\mbox{Poly}}_{f,\Psi_{\bm{j}},\mathscr{L}}\left(\bm{I}_i^{\bigtimes}\right)\right)}\left(x - \left(\sum\limits_{i=1}^n \beta_i x_i\right)^q\right)\bm{I}_{\mathfrak{K}}.
\end{eqnarray}

(II) If we assume that $g(\bm{x})=\log\left(\sum\limits_{i=1}^n \beta_i x_i\right)$, then, we have the following upper bound for $\sum\limits_{\bm{j}=\bm{1}}^{\bm{k}}\bm{w}^{\times}_{\bm{j}}\Phi_{\bm{j}}(f(\underline{\bm{A}}_{\bm{j}}))$:
\begin{eqnarray}\label{eq log x UB: cor: 2.17 diff}
\lefteqn{\sum\limits_{\bm{j}=\bm{1}}^{\bm{k}}\bm{w}^{\times}_{\bm{j}}\Phi_{\bm{j}}(f(\underline{\bm{A}}_{\bm{j}})) - \log\left(\sum\limits_{i=1}^n \beta_i \left(\sum\limits_{\bm{j}=\bm{1}}^{\bm{k}}\bm{w}^{\times}_{\bm{j}}\bm{V}^{\ast}\mbox{Poly}_{f,\Psi_{\bm{j}},\mathscr{U}}(\underline{\bm{A}}_{j_i})\bm{V}\right)\right)}\nonumber \\
&\leq& 
\max\limits_{x \in \left(\bigcup\limits_{\bm{j}=\bm{1}}^{\bm{k}}\bm{w}^{\times}_{\bm{j}}\widetilde{\mbox{Poly}}_{f,\Psi_{\bm{j}},\mathscr{U}}\left(\bm{I}^{\bigtimes}\right)\right),x_i \in \left(\bigcup\limits_{\bm{j}=\bm{1}}^{\bm{k}}\bm{w}^{\times}_{\bm{j}}\widetilde{\mbox{Poly}}_{f,\Psi_{\bm{j}},\mathscr{U}}\left(\bm{I}_i^{\bigtimes}\right)\right)}\left(x - \log\left(\sum\limits_{i=1}^n \beta_i x_i\right)\right)\bm{I}_{\mathfrak{K}};
\end{eqnarray}
and, the following lower bound for $\sum\limits_{\bm{j}=\bm{1}}^{\bm{k}}\bm{w}^{\times}_{\bm{j}}\Phi_{\bm{j}}(f(\underline{\bm{A}}_{\bm{j}}))$:
\begin{eqnarray}\label{eq log x LB: cor: 2.17 diff}
\lefteqn{\sum\limits_{\bm{j}=\bm{1}}^{\bm{k}}\bm{w}^{\times}_{\bm{j}}\Phi_{\bm{j}}(f(\underline{\bm{A}}_{\bm{j}})) - \log\left(\sum\limits_{i=1}^n \beta_i \left(\sum\limits_{\bm{j}=\bm{1}}^{\bm{k}}\bm{w}^{\times}_{\bm{j}}\bm{V}^{\ast}\mbox{Poly}_{f,\Psi_{\bm{j}},\mathscr{L}}(\underline{\bm{A}}_{j_i})\bm{V}\right)\right)}\nonumber \\
&\geq& 
\min\limits_{x \in \left(\bigcup\limits_{\bm{j}=\bm{1}}^{\bm{k}}\bm{w}^{\times}_{\bm{j}}\widetilde{\mbox{Poly}}_{f,\Psi_{\bm{j}},\mathscr{L}}\left(\bm{I}^{\bigtimes}\right)\right),x_i \in \left(\bigcup\limits_{\bm{j}=\bm{1}}^{\bm{k}}\bm{w}^{\times}_{\bm{j}}\widetilde{\mbox{Poly}}_{f,\Psi_{\bm{j}},\mathscr{L}}\left(\bm{I}_i^{\bigtimes}\right)\right)}\left(x - \log\left(\sum\limits_{i=1}^n \beta_i x_i\right)\right)\bm{I}_{\mathfrak{K}}.
\end{eqnarray}

(III) If we assume that $g(\bm{x})=\exp\left(\sum\limits_{i=1}^n \beta_i x_i\right)$, then, we have the following upper bound for $\sum\limits_{\bm{j}=\bm{1}}^{\bm{k}}\bm{w}^{\times}_{\bm{j}}\Phi_{\bm{j}}(f(\underline{\bm{A}}_{\bm{j}}))$:
\begin{eqnarray}\label{eq exp x UB: cor: 2.17 diff}
\lefteqn{\sum\limits_{\bm{j}=\bm{1}}^{\bm{k}}\bm{w}^{\times}_{\bm{j}}\Phi_{\bm{j}}(f(\underline{\bm{A}}_{\bm{j}})) - \exp\left(\sum\limits_{i=1}^n \beta_i \left(\sum\limits_{\bm{j}=\bm{1}}^{\bm{k}}\bm{w}^{\times}_{\bm{j}}\bm{V}^{\ast}\mbox{Poly}_{f,\Psi_{\bm{j}},\mathscr{U}}(\underline{\bm{A}}_{j_i})\bm{V}\right)\right)}\nonumber \\
&\leq& 
\max\limits_{x \in \left(\bigcup\limits_{\bm{j}=\bm{1}}^{\bm{k}}\bm{w}^{\times}_{\bm{j}}\widetilde{\mbox{Poly}}_{f,\Psi_{\bm{j}},\mathscr{U}}\left(\bm{I}^{\bigtimes}\right)\right),x_i \in \left(\bigcup\limits_{\bm{j}=\bm{1}}^{\bm{k}}\bm{w}^{\times}_{\bm{j}}\widetilde{\mbox{Poly}}_{f,\Psi_{\bm{j}},\mathscr{U}}\left(\bm{I}_i^{\bigtimes}\right)\right)}\left(x - \exp\left(\sum\limits_{i=1}^n \beta_i x_i\right)\right)\bm{I}_{\mathfrak{K}};
\end{eqnarray}
and, the following lower bound for $\sum\limits_{\bm{j}=\bm{1}}^{\bm{k}}\bm{w}^{\times}_{\bm{j}}\Phi_{\bm{j}}(f(\underline{\bm{A}}_{\bm{j}}))$:
\begin{eqnarray}\label{eq exp x LB: cor: 2.17 diff}
\lefteqn{\sum\limits_{\bm{j}=\bm{1}}^{\bm{k}}\bm{w}^{\times}_{\bm{j}}\Phi_{\bm{j}}(f(\underline{\bm{A}}_{\bm{j}})) - \exp\left(\sum\limits_{i=1}^n \beta_i \left(\sum\limits_{\bm{j}=\bm{1}}^{\bm{k}}\bm{w}^{\times}_{\bm{j}}\bm{V}^{\ast}\mbox{Poly}_{f,\Psi_{\bm{j}},\mathscr{L}}(\underline{\bm{A}}_{j_i})\bm{V}\right)\right)}\nonumber \\
&\geq& 
\min\limits_{x \in \left(\bigcup\limits_{\bm{j}=\bm{1}}^{\bm{k}}\bm{w}^{\times}_{\bm{j}}\widetilde{\mbox{Poly}}_{f,\Psi_{\bm{j}},\mathscr{L}}\left(\bm{I}^{\bigtimes}\right)\right),x_i \in \left(\bigcup\limits_{\bm{j}=\bm{1}}^{\bm{k}}\bm{w}^{\times}_{\bm{j}}\widetilde{\mbox{Poly}}_{f,\Psi_{\bm{j}},\mathscr{L}}\left(\bm{I}_i^{\bigtimes}\right)\right)}\left(x - \exp\left(\sum\limits_{i=1}^n \beta_i x_i\right)\right)\bm{I}_{\mathfrak{K}}.
\end{eqnarray}
\end{corollary}
\textbf{Proof:}
For Part (I), we will use $g(\bm{x})=\left(\sum\limits_{i=1}^n \beta_i x_i\right)^q$ in Eq.~\eqref{eq UB: cor 2.15} in Theorem~\ref{thm: cor 2.15} to obtain Eq.~\eqref{eq x q UB: cor: 2.17 diff}. We will also use $g(\bm{x})=\left(\sum\limits_{i=1}^n \beta_i x_i\right)^q$ in Eq.~\eqref{eq LB: cor 2.15} in Theorem~\ref{thm: cor 2.15} to obtain Eq.~\eqref{eq x q LB: cor: 2.17 diff}. 

For Part (II), we will use $g(\bm{x})=\log\left(\sum\limits_{i=1}^n \beta_i x_i\right)$ in Eq.~\eqref{eq UB: cor 2.15} in Theorem~\ref{thm: cor 2.15} to obtain Eq.~\eqref{eq log x UB: cor: 2.17 diff}. We will also use $g(\bm{x})=\log\left(\sum\limits_{i=1}^n \beta_i x_i\right)$ in Eq.~\eqref{eq LB: cor 2.15} in Theorem~\ref{thm: cor 2.15} to obtain Eq.~\eqref{eq log x LB: cor: 2.17 diff}. 

For Part (III), we will use $g(\bm{x})=\exp\left(\sum\limits_{i=1}^n \beta_i x_i\right)$ in Eq.~\eqref{eq UB: cor 2.15} in Theorem~\ref{thm: cor 2.15} to obtain Eq.~\eqref{eq exp x UB: cor: 2.17 diff}. We will also use $g(\bm{x})=\exp\left(\sum\limits_{i=1}^n \beta_i x_i\right)$ in Eq.~\eqref{eq LB: cor 2.15} in Theorem~\ref{thm: cor 2.15} to obtain Eq.~\eqref{eq exp x LB: cor: 2.17 diff}. 
$\hfill \Box$

\section{$W$-boundedness with Loewner Order}\label{sec: W-boundedness with Loewner Order}

In \cite{hytonen2016analysis_I,hytonen2017analysis_II}, the authors initiated their attempt by proving an operator-valued version of the Mihlin multiplier theorem. They replaced the uniform boundedness condition on certain operator families with the more stringent requirement of $R$-boundedness. This stricter condition is highly versatile and applicable to numerous operator-theoretic results in Hilbert spaces. The current volume provides a comprehensive set of analytical methods to verify the $R$-boundedness of many classical operators, which are essential for applications in harmonic analysis and stochastic analysis.

Given a family of linear mappings, represented by $\mathfrak{T}$, from $\mathbb{B}(\mathfrak{H})$ to $\mathbb{B}(\mathfrak{K})$, we say $\mathfrak{T}$ is a $R$-boundness if there exists a constant $C \geq 0$ such that for all finite set of linear operators $T_j \in \mathfrak{T}$ for $j=1,2,\ldots,k$ and $\bm{A}_j \in \mathbb{B}(\mathfrak{H})$ for $j=1,2,\ldots,k$, we have
\begin{eqnarray}\label{eq: R bound def}
\left\Vert\sum\limits_{j=1}^k r_j T_j \bm{A}_j \right\Vert_{L^2(\mathfrak{K})} \leq C \left\Vert\sum\limits_{j=1}^k r_j \bm{A}_j \right\Vert_{L^2(\mathfrak{H})},
\end{eqnarray}
where $r_j$ are Rademacher random variables. 

From Theorem~\ref{thm: 2.9}, we have the following Corollary~\ref{cor: thm 2.9}.
\begin{corollary}\label{cor: thm 2.9}
Let $\bm{A}_j$ for $j=1,2,\ldots,k$ be self-adjoint operators with $\Lambda(\bm{A}_j) \in [m,M]$, where $m<M$, $\Phi_j$ for $j=1,2,\ldots,k$ be nomralized positive linear maps, $\bm{w}=[w_1,w_2,\ldots,w_k]$ be a probability vector with $\sum\limits_{j=1}^k w_j = 1$. Given $k$ convex functions $f_j$ for $j=1,2,\ldots,k$ defined on $[m,M]$ and a continuous function $g$ defined on $[m,M]$ with $g(x) > 0$ on $[m,M]$, then we have
\begin{eqnarray}\label{eq1: cor: thm 2.9}
\sum\limits_{j=1}^k w_j \Phi_j (f_j(\bm{A}_j))&\leq&\left(\max\limits_{m \leq x \leq M}\frac{\max\limits_{j=1,2,\ldots,k}(a_j x+b_j)}{g(x)} \right)g\left(\sum\limits_{j=1}^k w_j \Phi_j (\bm{A}_j)\right),
\end{eqnarray}
where $a_j = \frac{f_j(M) - f_j(m)}{M-m}$ and $b_j = \frac{Mf_j(m) - m f_j(M)}{M-m}$.
\end{corollary}
From Corollary~\ref{cor: thm 2.9}, we can define $W$-boundedness with Loewner order with respect to a set of hypercomplex complex functions. Given a family of convex functions defined on $[m,M]$, denoted by $\mathfrak{F}(m.M)$, we say $\mathfrak{F}$ is a $W$-boundness if there exists a constant $C \geq 0$ such that for all finite set of convex hypercomplex functions $f_j \in \mathfrak{F}$ and $\bm{A}_j \in \mathbb{B}(\mathfrak{H})$ for $j=1,2,\ldots,k$, we have
\begin{eqnarray}\label{eq: W bound def}
\sum\limits_{j=1}^k w_j \Phi (f_j(\bm{A}_j))&\leq&C g\left(\sum\limits_{j=1}^k w_j \Phi (\bm{A}_j)\right),
\end{eqnarray}
where $w_j$ are entries of a probability vector. We use $\mathscr{W}_{\Phi, g}(\mathfrak{F}(m.M))$ to represent the least admissible constant $C$ in Eq.~\eqref{eq: W bound def} with respect to the normalized positive linear map $\Phi$ and the continuous function $g$ defined on $[m,M]$. $\mathscr{W}_{\Phi, g}(\mathfrak{F}(m.M))$ will be named as the \emph{$W$-bound} with respect to the map $\Phi$, the continuous function $g$, and a family of convex functions $\mathfrak{F}(m.M)$. We have the following two propositions about $\mathscr{W}_{\Phi, g}(\mathfrak{F}(m.M))$.

\begin{proposition}\label{prop: W bound}
If two functions $g_1$ and $g_2$ as $g$ in Eq~\eqref{eq: W bound def} have the relationship $g_1(x)=\kappa g_2(x)$ for any positive constant $\kappa$ and $x \in [m,M]$, we have
\begin{eqnarray}\label{eq1: prop: W bound}
\mathscr{W}_{\Phi, g_1}(\mathfrak{F}(m.M))=\frac{1}{\kappa}\mathscr{W}_{\Phi, g_2}(\mathfrak{F}(m.M)). 
\end{eqnarray}
\end{proposition}
\textbf{Proof:}
This is immediatedly from Eq.~\eqref{eq1: cor: thm 2.9} and the defnition about $\mathscr{W}_{\Phi, g}(\mathfrak{F}(m.M))$.
$\hfill \Box$

\begin{proposition}(Domination by a linear operator)\label{prop: Domination by a linear operator}
Let $\bm{A}_j$ for $j=1,2,\ldots,k$ be self-adjoint operators with $\Lambda(\bm{A}_j) \in [m,M]$, where $m<M$, $\Phi_j$ for $j=1,2,\ldots,k$ be nomralized positive linear maps, $\bm{w}=[w_1,w_2,\ldots,w_k]$ be a probability vector with $\sum\limits_{j=1}^k w_j = 1$. Given $k$ convex functions $f_j$ for $j=1,2,\ldots,k$ defined on $[m,M]$ and a continuous function $g$ defined on $[m,M]$ with $g(x) > 0$ on $[m,M]$. Besides, we introduce the operator $\bm{U}$, and assume that $\Phi_j(\bm{U}\bm{A}_j)=\bm{U}\Phi_j(\bm{A}_j)$ and $g(\bm{U}\bm{Y}) \leq \upsilon_g\left\Vert\bm{U}\right\Vert g(\bm{Y})$ for any self-adjoint operator $\bm{Y}$, where $\upsilon_g$ is a positive constant depends on $g$ only. Then we have
\begin{eqnarray}\label{eq0: prop: Domination by a linear operator}
\mathscr{W}_{\Phi, g}(\mathfrak{F}_{\bm{U}}(m.M)) \leq \upsilon_g\left\Vert\bm{U}\right\Vert,
\end{eqnarray}
where $\mathfrak{F}_{\bm{U}}(m.M)$ represents the family of the convex functions. 
\end{proposition}
\textbf{Proof:}
From Corollary~\ref{cor: thm 2.9}, we have
\begin{eqnarray}\label{eq1: prop: Domination by a linear operator}
\sum\limits_{j=1}^k w_j \Phi_j (f_j(\bm{A}_j))&\leq&\left(\max\limits_{m \leq x \leq M}\frac{\max\limits_{j=1,2,\ldots,k}(a_j x+b_j)}{g(x)} \right)g\left(\sum\limits_{j=1}^k w_j \Phi_j (\bm{U}\bm{A}_j)\right) \nonumber \\
&=&\left(\max\limits_{m \leq x \leq M}\frac{\max\limits_{j=1,2,\ldots,k}(a_j x+b_j)}{g(x)} \right)g\left(\bm{U}\sum\limits_{j=1}^k w_j \Phi_j (\bm{A}_j)\right) \nonumber \\
&\leq&\left(\max\limits_{m \leq x \leq M}\frac{\max\limits_{j=1,2,\ldots,k}(a_j x+b_j)}{g(x)} \right)\upsilon_g\left\Vert\bm{U}\right\Vert g\left(\sum\limits_{j=1}^k w_j \Phi_j (\bm{A}_j)\right),
\end{eqnarray}
this proposition is proved by $\mathscr{W}_{\Phi, g}(\mathfrak{F}(m.M))$ definition. 
$\hfill \Box$

\begin{remark}
In this section, we assume that the function $f$ is convex with a single input variable. The concept of $W$-boundedness, discussed here, is analogous to the conventional $R$-boundedness. However, the inequalities established in Section~\ref{sec: Multivariate Hypercomplex Functions Inequalities: Ratio Kind} can be utilized to extend the current $W$-boundedness framework from one input variable to multiple variables.
\end{remark}

\section{Multivariate Hypercomplex Function Approximation}\label{sec: Multivariate Hypercomplex Function Approximations}

In this section, we will consider multivariate hypercomplex function approximation problem in terms of the ratio error discussed in Section~\ref{sec: Ratio Type Approximation} and the difference error in Section~\ref{sec: Difference Type Approximation}.

\subsection{Ratio Type Approximation}\label{sec: Ratio Type Approximation}

Let $\bm{A}_1, \bm{A}_2, \ldots, \bm{A}_n$ be $n$ self-adjoint operators with $\Lambda(\bm{A}_1) \times \Lambda(\bm{A}_2) \times \ldots \times \Lambda(\bm{A}_n) \in \bigtimes\limits_{i=1}^n [m_i, M_i] \in \mathbb{R}^n$ for real scalars $m_i <  M_i$. The mapping $\Phi: \mathbb{B}(\mathfrak{H}) \rightarrow \mathbb{B}(\mathfrak{K})$ is defined by Eq.~\eqref{eq: new phi def}. Let $f(\bm{x})$ be any real valued continuous functions with $n$ input variables defined on the range $\bigtimes\limits_{i=1}^n [m_i, M_i]$, represented by $f(\bm{x}) \in \mathcal{C}\left(\bigtimes\limits_{i=1}^n [m_i, M_i]\right)$. The \emph{ratio type approximation} problem is to find the invertible function $g(\bm{x})$ and $n$ polynomial type functions, represented by $\mbox{Poly}_{f,\Psi,\mathscr{U}}(\underline{\bm{A}}_{i})$ for $i=1,2,\ldots,n$, to satisfy the following:
\begin{eqnarray}\label{eq1: UP Ratio Type Approximation}
\frac{\Phi(f(\underline{\bm{A}}))}{g\left(\bm{V}^{\ast}\mbox{Poly}_{f,\Psi,\mathscr{U}}(\underline{\bm{A}}_{1})\bm{V},\ldots,\bm{V}^{\ast}\mbox{Poly}_{f,\Psi,\mathscr{U}}(\underline{\bm{A}}_{n})\bm{V}\right))}\leq \alpha_1 \bm{I}_{\mathfrak{K}},
\end{eqnarray}
where $\alpha_1$ is some specified positive real number. Similarly, we also can find the invertible function $g(\bm{x})$ and $n$ polynomial type functions, represented by $\mbox{Poly}_{f,\Psi,\mathscr{L}}(\underline{\bm{A}}_{i})$ for $i=1,2,\ldots,n$, to satisfy the following:
\begin{eqnarray}\label{eq1: LO Ratio Type Approximation}
\frac{\Phi(f(\underline{\bm{A}}))}{g\left(\bm{V}^{\ast}\mbox{Poly}_{f,\Psi,\mathscr{L}}(\underline{\bm{A}}_{1})\bm{V},\ldots,\bm{V}^{\ast}\mbox{Poly}_{f,\Psi,\mathscr{L}}(\underline{\bm{A}}_{n})\bm{V}\right))}\geq \alpha_2 \bm{I}_{\mathfrak{K}},
\end{eqnarray}
where $\alpha_2$ is some specified positive real number.

From Theorem~\ref{thm: 2.9}, if $g\left(\bm{V}^{\ast}\mbox{Poly}_{f,\Psi,\mathscr{U}}(\underline{\bm{A}}_{1})\bm{V},\ldots,\bm{V}^{\ast}\mbox{Poly}_{f,\Psi,\mathscr{U}}(\underline{\bm{A}}_{n})\bm{V}\right)) > \bm{0}$, Eq.~\eqref{eq1: UP Ratio Type Approximation} can be established by setting 
\begin{eqnarray}\label{eq2: UP Ratio Type Approximation}
\alpha_1 &\geq& \left[\max\limits_{x \in \left(\widetilde{\mbox{Poly}}_{f,\Psi,\mathscr{U}}\left(\bm{I}^{\bigtimes}\right)\right), \bm{x} \in \bigtimes\limits_{i=1}^n\left(\widetilde{\mbox{Poly}}_{f,\Psi,\mathscr{U}}\left(\bm{I}^{\bigtimes}_i\right)\right)}xg^{-1}(\bm{x})\right] 
\end{eqnarray}
Similarly,  if $g\left(\bm{V}^{\ast}\mbox{Poly}_{f,\Psi,\mathscr{L}}(\underline{\bm{A}}_{1})\bm{V},\ldots,\bm{V}^{\ast}\mbox{Poly}_{f,\Psi,\mathscr{L}}(\underline{\bm{A}}_{n})\bm{V}\right))> \bm{0}$, Eq.~\eqref{eq1: LO Ratio Type Approximation} can be established by setting 
\begin{eqnarray}\label{eq2: LO Ratio Type Approximation}
\alpha_2 &\leq& \left[\max\limits_{x \in \left(\widetilde{\mbox{Poly}}_{f,\Psi,\mathscr{L}}\left(\bm{I}^{\bigtimes}\right)\right), \bm{x} \in \bigtimes\limits_{i=1}^n\left(\widetilde{\mbox{Poly}}_{f,\Psi,\mathscr{L}}\left(\bm{I}^{\bigtimes}_i\right)\right)}xg^{-1}(\bm{x})\right].
\end{eqnarray}

On the other hand, from Theorem~\ref{thm: 2.9}, if $g\left(\bm{V}^{\ast}\mbox{Poly}_{f,\Psi,\mathscr{L}}(\underline{\bm{A}}_{1})\bm{V},\ldots,\bm{V}^{\ast}\mbox{Poly}_{f,\Psi,\mathscr{L}}(\underline{\bm{A}}_{n})\bm{V}\right)) < \bm{0}$, Eq.~\eqref{eq1: UP Ratio Type Approximation} can be established by setting 
\begin{eqnarray}\label{eq2-1: UP Ratio Type Approximation}
\alpha_1 &\geq&\left[\min\limits_{x \in \left(\widetilde{\mbox{Poly}}_{f,\Psi,\mathscr{L}}\left(\bm{I}^{\bigtimes}\right)\right), \bm{x} \in \bigtimes\limits_{i=1}^n\left(\widetilde{\mbox{Poly}}_{f,\Psi,\mathscr{L}}\left(\bm{I}^{\bigtimes}_i\right)\right)}xg^{-1}(\bm{x})\right] . 
\end{eqnarray}
Similarly,  if  $g\left(\bm{V}^{\ast}\mbox{Poly}_{f,\Psi,\mathscr{U}}(\underline{\bm{A}}_{1})\bm{V},\ldots,\bm{V}^{\ast}\mbox{Poly}_{f,\Psi,\mathscr{U}}(\underline{\bm{A}}_{n})\bm{V}\right)) <  \bm{0}$, Eq.~\eqref{eq1: LO Ratio Type Approximation} can be established by setting 
\begin{eqnarray}\label{eq2-1: LO Ratio Type Approximation}
\alpha_2 &\leq& \left[\max\limits_{x \in \left(\widetilde{\mbox{Poly}}_{f,\Psi,\mathscr{U}}\left(\bm{I}^{\bigtimes}\right)\right), \bm{x} \in \bigtimes\limits_{i=1}^n\left(\widetilde{\mbox{Poly}}_{f,\Psi,\mathscr{U}}\left(\bm{I}^{\bigtimes}_i\right)\right)}xg^{-1}(\bm{x})\right].
\end{eqnarray}

Let $\bm{A}_{j_i}$ be self-adjoint operators with $\Lambda(\bm{A}_{j_i}) \in [m_i, M_i]$ for real scalars $m_i <  M_i$. The mappings $\Phi_{j_1,\ldots,j_n}: \mathscr{B}(\mathfrak{H}) \rightarrow \mathscr{B}(\mathfrak{K})$ are defined by Eq.~\eqref{eq: new phi def}, where $j_i=1,2,\ldots,k_i$ for $i=1,2,\ldots,n$. We have $n$ probability vectors $\bm{w}_i =[w_{i,1},w_{i,2},\cdots, w_{i,k_i}]$ with the dimension $k_i$ for $i=1,2,\ldots,n$, i.e., $\sum\limits_{\ell=1}^{k_i}w_{i,\ell} = 1$. Let $f(\bm{x})$ be any real valued continuous functions with $n$ variables defined on the range $\bigtimes\limits_{i=1}^n [m_i, M_i] \in \mathbb{R}^n$, where $\bigtimes$ is the Cartesian product. Besides, given any $\epsilon>0$, we assume that the function $f(\bm{x})$ satisfies the following:
\begin{eqnarray}\label{eq: lower and upper Psi Ratio Type Approximation ratio}
0 &\leq& \Psi_{\mathscr{U}}(\bm{x}) - f(\bm{x})~~\leq~~\epsilon, \nonumber \\
0 &\leq& f(\bm{x})-\Psi_{\mathscr{L}}(\bm{x})~~\leq~~\epsilon,
\end{eqnarray}
for $\bm{x}\in \bm{I}^{\bigtimes}$. The function $g(\bm{x})$ is also a real valued continuous function with $n$ variables defined on the range $\bigtimes\limits_{i=1}^n [m_i, M_i]$. The \emph{ratio type approximation} problem is to find the invertible function $g(\bm{x})$ and $n$ arguments of the function $g(\bm{x})$ with the format $\sum\limits_{\bm{j}=\bm{1}}^{\bm{k}}\bm{w}^{\times}_{\bm{j}}\bm{V}^{\ast}\mbox{Poly}_{f,\Psi_{\bm{j}},\mathscr{U}}(\underline{\bm{A}}_{j_i})\bm{V}$ for $i=1,2,\ldots,n$, to satisfy the following:
\begin{eqnarray}\label{eq3: UP Ratio Type Approximation}
\frac{\sum\limits_{\bm{j}=\bm{1}}^{\bm{k}}\bm{w}^{\times}_{\bm{j}}\Phi_{\bm{j}}(f(\underline{\bm{A}}_{\bm{j}}))}{g\Bigg(\sum\limits_{\bm{j}=\bm{1}}^{\bm{k}}\bm{w}^{\times}_{\bm{j}}\bm{V}^{\ast}\mbox{Poly}_{f,\Psi_{\bm{j}},\mathscr{U}}(\underline{\bm{A}}_{j_1})\bm{V},\ldots,\sum\limits_{\bm{j}=\bm{1}}^{\bm{k}}\bm{w}^{\times}_{\bm{j}}\bm{V}^{\ast}\mbox{Poly}_{f,\Psi_{\bm{j}},\mathscr{U}}(\underline{\bm{A}}_{j_n})\bm{V}\Bigg)}\leq\alpha_1 \bm{I}_{\mathfrak{K}},
\end{eqnarray}
where $\alpha_1$ is some specfied positive real number. Similarly, we also can find the invertible function $g(\bm{x})$ and $n$ arguments of the function $g(\bm{x})$ with the format $\sum\limits_{\bm{j}=\bm{1}}^{\bm{k}}\bm{w}^{\times}_{\bm{j}}\bm{V}^{\ast}\mbox{Poly}_{f,\Psi_{\bm{j}},\mathscr{L}}(\underline{\bm{A}}_{j_i})\bm{V}$ for $i=1,2,\ldots,n$, to satisfy the following:
\begin{eqnarray}\label{eq3: LO Ratio Type Approximation}
\frac{\sum\limits_{\bm{j}=\bm{1}}^{\bm{k}}\bm{w}^{\times}_{\bm{j}}\Phi_{\bm{j}}(f(\underline{\bm{A}}_{\bm{j}}))}{g\Bigg(\sum\limits_{\bm{j}=\bm{1}}^{\bm{k}}\bm{w}^{\times}_{\bm{j}}\bm{V}^{\ast}\mbox{Poly}_{f,\Psi_{\bm{j}},\mathscr{L}}(\underline{\bm{A}}_{j_1})\bm{V},\ldots,\sum\limits_{\bm{j}=\bm{1}}^{\bm{k}}\bm{w}^{\times}_{\bm{j}}\bm{V}^{\ast}\mbox{Poly}_{f,\Psi_{\bm{j}},\mathscr{L}}(\underline{\bm{A}}_{j_n})\bm{V}\Bigg)}\geq\alpha_2 \bm{I}_{\mathfrak{K}},
\end{eqnarray}
where $\alpha_2$ is some specfied positive real number.

From Theorem~\ref{thm: 2.9}, if $g\Bigg(\sum\limits_{\bm{j}=\bm{1}}^{\bm{k}}\bm{w}^{\times}_{\bm{j}}\bm{V}^{\ast}\mbox{Poly}_{f,\Psi_{\bm{j}},\mathscr{U}}(\underline{\bm{A}}_{j_1})\bm{V},\ldots,\sum\limits_{\bm{j}=\bm{1}}^{\bm{k}}\bm{w}^{\times}_{\bm{j}}\bm{V}^{\ast}\mbox{Poly}_{f,\Psi_{\bm{j}},\mathscr{U}}(\underline{\bm{A}}_{j_n})\bm{V}\Bigg) > \bm{0}$, Eq.~\eqref{eq3: UP Ratio Type Approximation} can be established by setting 
\begin{eqnarray}\label{eq4: UP Ratio Type Approximation}
\alpha_1 &\geq& \left[\max\limits_{x \in \left(\bigcup\limits_{\bm{j}=\bm{1}}^{\bm{k}}\bm{w}^{\times}_{\bm{j}}\widetilde{\mbox{Poly}}_{f,\Psi_{\bm{j}},\mathscr{U}}\left(\bm{I}^{\bigtimes}\right)\right), \bm{x} \in \bigtimes\limits_{i=1}^n\left(\bigcup\limits_{\bm{j}=\bm{1}}^{\bm{k}}\bm{w}^{\times}_{\bm{j}}\widetilde{\mbox{Poly}}_{f,\Psi_{\bm{j}},\mathscr{U}}\left(\bm{I}^{\bigtimes}_i\right)\right)}xg^{-1}(\bm{x})\right]. 
\end{eqnarray}
Similarly,  if $g\Bigg(\sum\limits_{\bm{j}=\bm{1}}^{\bm{k}}\bm{w}^{\times}_{\bm{j}}\bm{V}^{\ast}\mbox{Poly}_{f,\Psi_{\bm{j}},\mathscr{L}}(\underline{\bm{A}}_{j_1})\bm{V},\ldots,\sum\limits_{\bm{j}=\bm{1}}^{\bm{k}}\bm{w}^{\times}_{\bm{j}}\bm{V}^{\ast}\mbox{Poly}_{f,\Psi_{\bm{j}},\mathscr{L}}(\underline{\bm{A}}_{j_n})\bm{V}\Bigg) > \bm{0}$, Eq.~\eqref{eq3: LO Ratio Type Approximation} can be established by setting 
\begin{eqnarray}\label{eq4: LO Ratio Type Approximation}
\alpha_2 &\leq& \left[\min\limits_{x \in \left(\bigcup\limits_{\bm{j}=\bm{1}}^{\bm{k}}\bm{w}^{\times}_{\bm{j}}\widetilde{\mbox{Poly}}_{f,\Psi_{\bm{j}},\mathscr{L}}\left(\bm{I}^{\bigtimes}\right)\right), \bm{x} \in \bigtimes\limits_{i=1}^n\left(\bigcup\limits_{\bm{j}=\bm{1}}^{\bm{k}}\bm{w}^{\times}_{\bm{j}}\widetilde{\mbox{Poly}}_{f,\Psi_{\bm{j}},\mathscr{L}}\left(\bm{I}^{\bigtimes}_i\right)\right)}xg^{-1}(\bm{x})\right].
\end{eqnarray}

On the other hand, from Theorem~\ref{thm: 2.9}, if \\
$g\Bigg(\sum\limits_{\bm{j}=\bm{1}}^{\bm{k}}\bm{w}^{\times}_{\bm{j}}\bm{V}^{\ast}\mbox{Poly}_{f,\Psi_{\bm{j}},\mathscr{U}}(\underline{\bm{A}}_{j_1})\bm{V},\ldots,\sum\limits_{\bm{j}=\bm{1}}^{\bm{k}}\bm{w}^{\times}_{\bm{j}}\bm{V}^{\ast}\mbox{Poly}_{f,\Psi_{\bm{j}},\mathscr{U}}(\underline{\bm{A}}_{j_n})\bm{V}\Bigg) < \bm{0}$, Eq.~\eqref{eq3: UP Ratio Type Approximation} can be established by setting 
\begin{eqnarray}\label{eq4-1: UP Ratio Type Approximation}
\alpha_1 &\geq& \left[\min\limits_{x \in \left(\bigcup\limits_{\bm{j}=\bm{1}}^{\bm{k}}\bm{w}^{\times}_{\bm{j}}\widetilde{\mbox{Poly}}_{f,\Psi_{\bm{j}},\mathscr{L}}\left(\bm{I}^{\bigtimes}\right)\right), \bm{x} \in \bigtimes\limits_{i=1}^n\left(\bigcup\limits_{\bm{j}=\bm{1}}^{\bm{k}}\bm{w}^{\times}_{\bm{j}}\widetilde{\mbox{Poly}}_{f,\Psi_{\bm{j}},\mathscr{L}}\left(\bm{I}^{\bigtimes}_i\right)\right)}xg^{-1}(\bm{x})\right].  
\end{eqnarray}
Similarly,  if $g\Bigg(\sum\limits_{\bm{j}=\bm{1}}^{\bm{k}}\bm{w}^{\times}_{\bm{j}}\bm{V}^{\ast}\mbox{Poly}_{f,\Psi_{\bm{j}},\mathscr{L}}(\underline{\bm{A}}_{j_1})\bm{V},\ldots,\sum\limits_{\bm{j}=\bm{1}}^{\bm{k}}\bm{w}^{\times}_{\bm{j}}\bm{V}^{\ast}\mbox{Poly}_{f,\Psi_{\bm{j}},\mathscr{L}}(\underline{\bm{A}}_{j_n})\bm{V}\Bigg)  < \bm{0}$, Eq.~\eqref{eq3: LO Ratio Type Approximation} can be established by setting 
\begin{eqnarray}\label{eq4-1: LO Ratio Type Approximation}
\alpha_2 &\leq&\left[\max\limits_{x \in \left(\bigcup\limits_{\bm{j}=\bm{1}}^{\bm{k}}\bm{w}^{\times}_{\bm{j}}\widetilde{\mbox{Poly}}_{f,\Psi_{\bm{j}},\mathscr{U}}\left(\bm{I}^{\bigtimes}\right)\right), \bm{x} \in \bigtimes\limits_{i=1}^n\left(\bigcup\limits_{\bm{j}=\bm{1}}^{\bm{k}}\bm{w}^{\times}_{\bm{j}}\widetilde{\mbox{Poly}}_{f,\Psi_{\bm{j}},\mathscr{U}}\left(\bm{I}^{\bigtimes}_i\right)\right)}xg^{-1}(\bm{x})\right].
\end{eqnarray}

The following Example~\ref{exp: thm 2.9} is provided to evaluate the upper ratio bound given by Eq.~\eqref{eq3: UP Ratio Type Approximation} and evaluate the lower ratio bound given by Eq.~\eqref{eq3: LO Ratio Type Approximation}. 
\begin{example}\label{exp: thm 2.9}
Let $\bm{A}_{j_i}$ be self-adjoint operators with $\Lambda(\bm{A}_{j_i}) \in [m_i, M_i]$ for real scalars $m_i <  M_i$. The mappings $\Phi_{j_1,\ldots,j_n}: \mathscr{B}(\mathfrak{H}) \rightarrow \mathscr{B}(\mathfrak{K})$ are normalized positive linear maps, where $j_i=1,2,\ldots,k_i$ for $i=1,2,\ldots,n$. We have $n$ probability vectors $\bm{w}_i =[w_{i,1},w_{i,2},\cdots, w_{i,k_i}]$ with the dimension $k_i$ for $i=1,2,\ldots,n$, i.e., $\sum\limits_{\ell=1}^{k_i}w_{i,\ell} = 1$. Let $f(x_1,x_2,\ldots,x_n)$ be any real valued continuous functions with $n$ variables defined on the range $\bigtimes\limits_{i=1}^n [m_i, M_i] \in \mathbb{R}^n$, where $\times$ is the Cartesian product. Besides, we assume that the function $f$ satisfies the following:
\begin{eqnarray}\label{eq f bounds: cor: 2.11-14}
\sum\limits_{i=1}^n a_i x_i + b \leq f(x_1,x_2,\ldots,x_n) \leq \sum\limits_{i=1}^n c_i x_i + d,
\end{eqnarray}
for $[x_1,x_2,\ldots,x_n] \in \bigtimes\limits_{i=1}^n [m_i,M_i]$. We will discuss the followingfour cases.      

(I) If we also assume that $g(x_1,\ldots,x_n) = \Big(\sum\limits_{i=1}^n \beta_i x_i\Big)^q$, where $\beta_i \geq 0$, $q \in \mathbb{R}$ and $m_i \geq 0$, then, we have the following upper bound for $\sum\limits_{j_{1}=1,\ldots,j_{n}=1}^{k_1,\ldots,k_n}  w_{j_1}\ldots w_{j_n}\Phi_{j_1,\ldots,j_n}(f(\bm{A}_{j_1},\bm{A}_{j_2},\ldots,\bm{A}_{j_n}))$:
\begin{eqnarray}\label{eq x q UB: cor: 2.11-14}
\lefteqn{\sum\limits_{j_{1}=1,\ldots,j_{n}=1}^{k_1,\ldots,k_n}  w_{j_1}\ldots w_{j_n}\Phi_{j_1,\ldots,j_n}(f(\bm{A}_{j_1},\bm{A}_{j_2},\ldots,\bm{A}_{j_n}))}\nonumber \\
&\leq&\left[\max\limits_{\substack{m_i \leq x_i \leq M_i \\ i=1,2,\ldots,n}}\frac{\sum\limits_{i=1}^n c_i x_i +d}{\Big(\sum\limits_{i=1}^n \beta_i x_i\Big)^q}\right]\left(\sum\limits_{i=1}^n \beta_i \sum\limits_{j_{1}=1,\ldots,j_{n}=1}^{k_1,\ldots,k_n}w_{j_1}\ldots w_{j_n}\Phi_{j_1,\ldots,j_n}(\bm{A}_{j_i})\right)^q;
\end{eqnarray}
and, the following lower bound for $\sum\limits_{j_{1}=1,\ldots,j_{n}=1}^{k_1,\ldots,k_n}  w_{j_1}\ldots w_{j_n}\Phi_{j_1,\ldots,j_n}(f(\bm{A}_{j_1},\bm{A}_{j_2},\ldots,\bm{A}_{j_n}))$:
\begin{eqnarray}\label{eq x q LB: cor: 2.11-14}
\lefteqn{\sum\limits_{j_{1}=1,\ldots,j_{n}=1}^{k_1,\ldots,k_n}  w_{j_1}\ldots w_{j_n}\Phi_{j_1,\ldots,j_n}(f(\bm{A}_{j_1},\bm{A}_{j_2},\ldots,\bm{A}_{j_n}))}\nonumber \\
&\geq&\left[\min\limits_{\substack{m_i \leq x_i \leq M_i \\ i=1,2,\ldots,n}}\frac{\sum\limits_{i=1}^n a_i x_i +b}{\Big(\sum\limits_{i=1}^n \beta_i x_i\Big)^q}\right]\left(\sum\limits_{i=1}^n \beta_i \sum\limits_{j_{1}=1,\ldots,j_{n}=1}^{k_1,\ldots,k_n}w_{j_1}\ldots w_{j_n}\Phi_{j_1,\ldots,j_n}(\bm{A}_{j_i})\right)^q.
\end{eqnarray}

Part (I) of this example is proved by applying Theorem~\ref{thm: 2.9} Part (I) with the function $g$ as $g(x_1,\ldots,x_n) = \Big(\sum\limits_{i=1}^n \beta_i x_i\Big)^q$. Then, we have
\begin{eqnarray}\label{eq1: exp: thm 2.9}
\alpha_1 &\geq& \left[\max\limits_{\substack{m_i \leq x_i \leq M_i \\ i=1,2,\ldots,n}}\frac{\sum\limits_{i=1}^n c_i x_i +d}{\Big(\sum\limits_{i=1}^n \beta_i x_i\Big)^q}\right], \nonumber \\
\alpha_2 &\leq& \left[\min\limits_{\substack{m_i \leq x_i \leq M_i \\ i=1,2,\ldots,n}}\frac{\sum\limits_{i=1}^n a_i x_i +b}{\Big(\sum\limits_{i=1}^n \beta_i x_i\Big)^q}\right]. 
\end{eqnarray}

(II) If $g(x_1,\ldots,x_n) = \log\Big(\sum\limits_{i=1}^n \beta_i x_i\Big)$ with $\log\Big(\sum\limits_{i=1}^n \beta_i x_i\Big)>0$ for $[x_1,\ldots,x_n]\in\bigtimes\limits_{i=1}^n [m_i, M_i]$, we have the upper bound for $\sum\limits_{j_{1}=1,\ldots,j_{n}=1}^{k_1,\ldots,k_n}  w_{j_1}\ldots w_{j_n}\Phi_{j_1,\ldots,j_n}(f(\bm{A}_{j_1},\bm{A}_{j_2},\ldots,\bm{A}_{j_n}))$:
\begin{eqnarray}\label{eq log x l0 UB: cor: 2.11-14}
\lefteqn{\sum\limits_{j_{1}=1,\ldots,j_{n}=1}^{k_1,\ldots,k_n}  w_{j_1}\ldots w_{j_n}\Phi_{j_1,\ldots,j_n}(f(\bm{A}_{j_1},\bm{A}_{j_2},\ldots,\bm{A}_{j_n}))}\nonumber \\
&\leq&\left[\max\limits_{\substack{m_i \leq x_i \leq M_i \\ i=1,2,\ldots,n}}\frac{\sum\limits_{i=1}^n c_i x_i +d}{\log\Big(\sum\limits_{i=1}^n \beta_i x_i\Big)}\right]\log\left(\sum\limits_{i=1}^n \beta_i \sum\limits_{j_{1}=1,\ldots,j_{n}=1}^{k_1,\ldots,k_n}w_{j_1}\ldots w_{j_n}\Phi_{j_1,\ldots,j_n}(\bm{A}_{j_i})\right);
\end{eqnarray}
and, the following lower bound for $\sum\limits_{j=1}^k w_j\Phi(f(\bm{A}_j))$:
\begin{eqnarray}\label{eq log x l0 LB: cor: 2.11-14}
\lefteqn{\sum\limits_{j_{1}=1,\ldots,j_{n}=1}^{k_1,\ldots,k_n}  w_{j_1}\ldots w_{j_n}\Phi_{j_1,\ldots,j_n}(f(\bm{A}_{j_1},\bm{A}_{j_2},\ldots,\bm{A}_{j_n}))}\nonumber \\
&\geq&\left[\min\limits_{\substack{m_i \leq x_i \leq M_i \\ i=1,2,\ldots,n}}\frac{\sum\limits_{i=1}^n a_i x_i +b}{\log\Big(\sum\limits_{i=1}^n \beta_i x_i\Big)}\right]\log\left(\sum\limits_{i=1}^n \beta_i \sum\limits_{j_{1}=1,\ldots,j_{n}=1}^{k_1,\ldots,k_n}w_{j_1}\ldots w_{j_n}\Phi_{j_1,\ldots,j_n}(\bm{A}_{j_i})\right).
\end{eqnarray}

Part (II) of this example is proved by applying Theorem~\ref{thm: 2.9} Part (I) with the function $g$ as $g(x_1,\ldots,x_n) = \log\Big(\sum\limits_{i=1}^n \beta_i x_i\Big)$, where $\log\Big(\sum\limits_{i=1}^n \beta_i x_i\Big)> 0$. Then, we have
\begin{eqnarray}\label{eq2: exp: thm 2.9}
\alpha_1 &\geq& \left[\max\limits_{\substack{m_i \leq x_i \leq M_i \\ i=1,2,\ldots,n}}\frac{\sum\limits_{i=1}^n c_i x_i +d}{\log\Big(\sum\limits_{i=1}^n \beta_i x_i\Big)}\right], \nonumber \\
\alpha_2 &\leq& \left[\min\limits_{\substack{m_i \leq x_i \leq M_i \\ i=1,2,\ldots,n}}\frac{\sum\limits_{i=1}^n a_i x_i +b}{\log\Big(\sum\limits_{i=1}^n \beta_i x_i\Big)}\right]. 
\end{eqnarray}

(III) If $g(x_1,\ldots,x_n) = \log\Big(\sum\limits_{i=1}^n \beta_i x_i\Big)$ with $\log\Big(\sum\limits_{i=1}^n \beta_i x_i\Big)<0$ for $[x_1,\ldots,x_n]\in\bigtimes\limits_{i=1}^n [m_i, M_i]$, we have the upper bound for $\sum\limits_{j_{1}=1,\ldots,j_{n}=1}^{k_1,\ldots,k_n}  w_{j_1}\ldots w_{j_n}\Phi_{j_1,\ldots,j_n}(f(\bm{A}_{j_1},\bm{A}_{j_2},\ldots,\bm{A}_{j_n}))$:
\begin{eqnarray}\label{eq log x s0 UB: cor: 2.11-14}
\lefteqn{\sum\limits_{j_{1}=1,\ldots,j_{n}=1}^{k_1,\ldots,k_n}  w_{j_1}\ldots w_{j_n}\Phi_{j_1,\ldots,j_n}(f(\bm{A}_{j_1},\bm{A}_{j_2},\ldots,\bm{A}_{j_n}))}\nonumber \\
&\leq&\left[\min\limits_{\substack{m_i \leq x_i \leq M_i \\ i=1,2,\ldots,n}}\frac{\sum\limits_{i=1}^n a_i x_i +b}{\log\Big(\sum\limits_{i=1}^n \beta_i x_i\Big)}\right]\log\left(\sum\limits_{i=1}^n \beta_i \sum\limits_{j_{1}=1,\ldots,j_{n}=1}^{k_1,\ldots,k_n}w_{j_1}\ldots w_{j_n}\Phi_{j_1,\ldots,j_n}(\bm{A}_{j_i})\right);
\end{eqnarray}
and, the following lower bound for $\sum\limits_{j=1}^k w_j\Phi(f(\bm{A}_j))$:
\begin{eqnarray}\label{eq log x s0 LB: cor: 2.11-14}
\lefteqn{\sum\limits_{j_{1}=1,\ldots,j_{n}=1}^{k_1,\ldots,k_n}  w_{j_1}\ldots w_{j_n}\Phi_{j_1,\ldots,j_n}(f(\bm{A}_{j_1},\bm{A}_{j_2},\ldots,\bm{A}_{j_n}))}\nonumber \\
&\geq&\left[\max\limits_{\substack{m_i \leq x_i \leq M_i \\ i=1,2,\ldots,n}}\frac{\sum\limits_{i=1}^n c_i x_i +d}{\log\Big(\sum\limits_{i=1}^n \beta_i x_i\Big)}\right]\log\left(\sum\limits_{i=1}^n \beta_i \sum\limits_{j_{1}=1,\ldots,j_{n}=1}^{k_1,\ldots,k_n}w_{j_1}\ldots w_{j_n}\Phi_{j_1,\ldots,j_n}(\bm{A}_{j_i})\right).
\end{eqnarray}

Part (III) of this example is proved by applying Theorem~\ref{thm: 2.9} Part (II) with the function $g$ as $g(x_1,\ldots,x_n) = \log\Big(\sum\limits_{i=1}^n \beta_i x_i\Big)$, where $\log\Big(\sum\limits_{i=1}^n \beta_i x_i\Big)< 0$. Then, we have
\begin{eqnarray}\label{eq3: exp: thm 2.9}
\alpha_1 &\geq& \left[\min\limits_{\substack{m_i \leq x_i \leq M_i \\ i=1,2,\ldots,n}}\frac{\sum\limits_{i=1}^n a_i x_i +b}{\log\Big(\sum\limits_{i=1}^n \beta_i x_i\Big)}\right], \nonumber \\
\alpha_2 &\leq& \left[\max\limits_{\substack{m_i \leq x_i \leq M_i \\ i=1,2,\ldots,n}}\frac{\sum\limits_{i=1}^n c_i x_i +d}{\log\Big(\sum\limits_{i=1}^n \beta_i x_i\Big)}\right]. 
\end{eqnarray}

(IV) If $g(x_1,\ldots,x_n) = \exp\Big(\sum\limits_{i=1}^n \beta_i x_i\Big)$, we have the upper bound for \\
$\sum\limits_{j_{1}=1,\ldots,j_{n}=1}^{k_1,\ldots,k_n}  w_{j_1}\ldots w_{j_n}\Phi_{j_1,\ldots,j_n}(f(\bm{A}_{j_1},\bm{A}_{j_2},\ldots,\bm{A}_{j_n}))$:
\begin{eqnarray}\label{eq exp x UB: cor: 2.11-14}
\lefteqn{\sum\limits_{j_{1}=1,\ldots,j_{n}=1}^{k_1,\ldots,k_n}  w_{j_1}\ldots w_{j_n}\Phi_{j_1,\ldots,j_n}(f(\bm{A}_{j_1},\bm{A}_{j_2},\ldots,\bm{A}_{j_n}))}\nonumber \\
&\leq&\left[\max\limits_{\substack{m_i \leq x_i \leq M_i \\ i=1,2,\ldots,n}}\frac{\sum\limits_{i=1}^n c_i x_i +d}{\exp\Big(\sum\limits_{i=1}^n \beta_i x_i\Big)}\right]\exp\left(\sum\limits_{i=1}^n \beta_i \sum\limits_{j_{1}=1,\ldots,j_{n}=1}^{k_1,\ldots,k_n}w_{j_1}\ldots w_{j_n}\Phi_{j_1,\ldots,j_n}(\bm{A}_{j_i})\right);
\end{eqnarray}
and, the following lower bound for  \\
$\sum\limits_{j_{1}=1,\ldots,j_{n}=1}^{k_1,\ldots,k_n}  w_{j_1}\ldots w_{j_n}\Phi_{j_1,\ldots,j_n}(f(\bm{A}_{j_1},\bm{A}_{j_2},\ldots,\bm{A}_{j_n}))$:
\begin{eqnarray}\label{eq exp x LB: cor: 2.11-14}
\lefteqn{\sum\limits_{j_{1}=1,\ldots,j_{n}=1}^{k_1,\ldots,k_n}  w_{j_1}\ldots w_{j_n}\Phi_{j_1,\ldots,j_n}(f(\bm{A}_{j_1},\bm{A}_{j_2},\ldots,\bm{A}_{j_n}))}\nonumber \\
&\geq&\left[\min\limits_{\substack{m_i \leq x_i \leq M_i \\ i=1,2,\ldots,n}}\frac{\sum\limits_{i=1}^n a_i x_i +b}{\exp\Big(\sum\limits_{i=1}^n \beta_i x_i\Big)}\right]\exp\left(\sum\limits_{i=1}^n \beta_i \sum\limits_{j_{1}=1,\ldots,j_{n}=1}^{k_1,\ldots,k_n}w_{j_1}\ldots w_{j_n}\Phi_{j_1,\ldots,j_n}(\bm{A}_{j_i})\right).
\end{eqnarray}

Finally, Part (IV) of this example is proved by applying Theorem~\ref{thm: 2.9} Part (I) with the function $g$ as $g(x_1,\ldots,x_n) = \exp\Big(\sum\limits_{i=1}^n \beta_i x_i\Big)$.   Then, we have
\begin{eqnarray}\label{eq4: exp: thm 2.9}
\alpha_1 &\geq& \left[\max\limits_{\substack{m_i \leq x_i \leq M_i \\ i=1,2,\ldots,n}}\frac{\sum\limits_{i=1}^n c_i x_i +d}{\exp\Big(\sum\limits_{i=1}^n \beta_i x_i\Big)}\right], \nonumber \\
\alpha_2 &\leq& \left[\min\limits_{\substack{m_i \leq x_i \leq M_i \\ i=1,2,\ldots,n}}\frac{\sum\limits_{i=1}^n a_i x_i +b}{\exp\Big(\sum\limits_{i=1}^n \beta_i x_i\Big)}\right]. 
\end{eqnarray}
\end{example}

\begin{remark}\label{rmk: ratio approx conjecture}
Given requirements of Eq.~\eqref{eq1: UP Ratio Type Approximation}, Eq.~\eqref{eq1: LO Ratio Type Approximation}, Eq.~\eqref{eq3: UP Ratio Type Approximation}, and Eq.~\eqref{eq3: LO Ratio Type Approximation}, we conjecture the existence of function $g(\bm{x})$ and polynomials in terms of $\Psi$ as arguments: $\mbox{Poly}_{f,\Psi,\mathscr{U}}(\bm{x})$, $\mbox{Poly}_{f,\Psi_{\bm{j}},\mathscr{U}}(\bm{x})$ and $\mbox{Poly}_{f,\Psi,\mathscr{L}}(\bm{x})$, $\mbox{Poly}_{f,\Psi_{\bm{j}}, \mathscr{L}}(\bm{x})$ for any given $\alpha_1$ or $\alpha_2$. If existence, how to find these functions?
\end{remark}

\subsection{Difference Type Approximation}\label{sec: Difference Type Approximation}

Let $\bm{A}_1, \bm{A}_2, \ldots, \bm{A}_n$ be $n$ self-adjoint operators with $\Lambda(\bm{A}_1) \times \Lambda(\bm{A}_2) \times \ldots \times \Lambda(\bm{A}_n) \in \bigtimes\limits_{i=1}^n [m_i, M_i] \in \mathbb{R}^n$ for real scalars $m_i <  M_i$. The mapping $\Phi: \mathbb{B}(\mathfrak{H}) \rightarrow \mathbb{B}(\mathfrak{K})$ is defined by Eq.~\eqref{eq: new phi def}. Let $f(\bm{x})$ be any real valued continuous functions with $n$ input variables defined on the range $\bigtimes\limits_{i=1}^n [m_i, M_i]$, represented by $f(\bm{x}) \in \mathcal{C}\left(\bigtimes\limits_{i=1}^n [m_i, M_i]\right)$. The \emph{difference type approximation} problem is to find the invertible function $g(\bm{x})$ and $n$ polynomial function in terms of $\Psi$ arguments, represented by $\mbox{Poly}_{f,\Psi,\mathscr{U}}(\underline{\bm{A}}_{i})$ for $i=1,2,\ldots,n$, to satisfy the following:
\begin{eqnarray}\label{eq1: UP Difference Type Approximation}
\Phi(f(\underline{\bm{A}})) - g\Bigg(\bm{V}^{\ast}\mbox{Poly}_{f,\Psi,\mathscr{U}}(\underline{\bm{A}}_{1})\bm{V},\ldots,\bm{V}^{\ast}\mbox{Poly}_{f,\Psi,\mathscr{U}}(\underline{\bm{A}}_{n})\bm{V}\Bigg)
\leq 
\beta_1 (x - g(\bm{x}))\bm{I}_{\mathfrak{K}}.
\end{eqnarray}
where $\beta_1$ is some specified real number. Similarly, we also can find $n$ polynomial function in terms of $\Psi$ arguments, represented by $\mbox{Poly}_{f,\Psi,\mathscr{L}}(\underline{\bm{A}}_{i})$ for $i=1,2,\ldots,n$, to satisfy the following:
\begin{eqnarray}\label{eq1: LO Difference Type Approximation}
\Phi(f(\underline{\bm{A}})) - g\Bigg(\bm{V}^{\ast}\mbox{Poly}_{f,\Psi,\mathscr{L}}(\underline{\bm{A}}_{1})\bm{V},\ldots,\bm{V}^{\ast}\mbox{Poly}_{f,\Psi,\mathscr{L}}(\underline{\bm{A}}_{n})\bm{V}\Bigg)\geq 
\beta_2 (x - g(\bm{x}))\bm{I}_{\mathfrak{K}},
\end{eqnarray}
where $\beta_2$ is some specified real number.

From Theorem~\ref{thm: cor 2.15}, Eq.~\eqref{eq1: UP Difference Type Approximation} can be established by setting 
\begin{eqnarray}\label{eq2: UP Difference Type Approximation}
\beta_1 &\geq& \max\limits_{x \in \left(\bigcup\limits_{\bm{j}=\bm{1}}^{\bm{k}}\bm{w}^{\times}_{\bm{j}}\widetilde{\mbox{Poly}}_{f,\Psi_{\bm{j}},\mathscr{U}}\left(\bm{I}^{\bigtimes}\right)\right), \bm{x} \in \bigtimes\limits_{i=1}^n\left(\bigcup\limits_{\bm{j}=\bm{1}}^{\bm{k}}\bm{w}^{\times}_{\bm{j}}\widetilde{\mbox{Poly}}_{f,\Psi_{\bm{j}},\mathscr{U}}\left(\bm{I}^{\bigtimes}_i\right)\right)}(x - g(\bm{x}))\bm{I}_{\mathfrak{K}}.
\end{eqnarray}
Similarly, Eq.~\eqref{eq1: LO Difference Type Approximation} can be established by setting 
\begin{eqnarray}\label{eq2: LO Difference Type Approximation}
\beta_2 &\leq& \min\limits_{x \in \left(\bigcup\limits_{\bm{j}=\bm{1}}^{\bm{k}}\bm{w}^{\times}_{\bm{j}}\widetilde{\mbox{Poly}}_{f,\Psi_{\bm{j}},\mathscr{L}}\left(\bm{I}^{\bigtimes}\right)\right), \bm{x} \in \bigtimes\limits_{i=1}^n\left(\bigcup\limits_{\bm{j}=\bm{1}}^{\bm{k}}\bm{w}^{\times}_{\bm{j}}\widetilde{\mbox{Poly}}_{f,\Psi_{\bm{j}},\mathscr{L}}\left(\bm{I}^{\bigtimes}_i\right)\right)}(x - g(\bm{x}))\bm{I}_{\mathfrak{K}}. 
\end{eqnarray}

Let $\bm{A}_{j_i}$ be self-adjoint operators with $\Lambda(\bm{A}_{j_i}) \in [m_i, M_i]$ for real scalars $m_i <  M_i$. The mappings $\Phi_{j_1,\ldots,j_n}: \mathscr{B}(\mathfrak{H}) \rightarrow \mathscr{B}(\mathfrak{K})$ are defined by Eq.~\eqref{eq: new phi def}, where $j_i=1,2,\ldots,k_i$ for $i=1,2,\ldots,n$. We have $n$ probability vectors $\bm{w}_i =[w_{i,1},w_{i,2},\cdots, w_{i,k_i}]$ with the dimension $k_i$ for $i=1,2,\ldots,n$, i.e., $\sum\limits_{\ell=1}^{k_i}w_{i,\ell} = 1$. Let $f(\bm{x})$ be any real valued continuous functions with $n$ variables defined on the range $\bigtimes\limits_{i=1}^n [m_i, M_i] \in \mathbb{R}^n$, where $\bigtimes$ is the Cartesian product. Besides, given any $\epsilon>0$, we assume that the function $f(\bm{x})$ satisfies the following:
\begin{eqnarray}\label{eq: lower and upper Psi Ratio Type Approximation difference}
0 &\leq& \Psi_{\mathscr{U}}(\bm{x}) - f(\bm{x})~~\leq~~\epsilon, \nonumber \\
0 &\leq& f(\bm{x})-\Psi_{\mathscr{L}}(\bm{x})~~\leq~~\epsilon,
\end{eqnarray}
for $\bm{x}\in \bm{I}^{\bigtimes}$. The function $g(\bm{x})$ is also a real valued continuous function with $n$ variables defined on the range $\bigtimes\limits_{i=1}^n [m_i, M_i]$. The \emph{difference type approximation} problem is to find the invertible function $g(\bm{x})$ and $n$ polynomial function in terms of $\Psi$ arguments, represented by $\mbox{Poly}_{f,\Psi_{\bm{j}},\mathscr{U}}(\underline{\bm{A}}_{j_i})$ for $i=1,2,\ldots,n$, to satisfy the following:
\begin{eqnarray}\label{eq3: UP Difference Type Approximation}
\lefteqn{\sum\limits_{\bm{j}=\bm{1}}^{\bm{k}}\bm{w}^{\times}_{\bm{j}}\Phi_{\bm{j}}(f(\underline{\bm{A}}_{\bm{j}}))}\nonumber \\
&&- g\Bigg(\sum\limits_{\bm{j}=\bm{1}}^{\bm{k}}\bm{w}^{\times}_{\bm{j}}\bm{V}^{\ast}\mbox{Poly}_{f,\Psi_{\bm{j}},\mathscr{U}}(\underline{\bm{A}}_{j_1})\bm{V},\ldots,\sum\limits_{\bm{j}=\bm{1}}^{\bm{k}}\bm{w}^{\times}_{\bm{j}}\bm{V}^{\ast}\mbox{Poly}_{f,\Psi_{\bm{j}},\mathscr{U}}(\underline{\bm{A}}_{j_n})\bm{V}\Bigg)\nonumber \\
&\leq&
\beta_1 (x - g(\bm{x}))\bm{I}_{\mathfrak{K}},
\end{eqnarray}
where $\beta_1$ is some specfied real number. Similarly, we also can find the polynomial function $p_{2,1}(x),\cdots,p_{2,k}(x)$  to satisfy the following:
\begin{eqnarray}\label{eq3: LO Difference Type Approximation}
\lefteqn{\sum\limits_{\bm{j}=\bm{1}}^{\bm{k}}\bm{w}^{\times}_{\bm{j}}\Phi_{\bm{j}}(f(\underline{\bm{A}}_{\bm{j}}))}\nonumber \\
&&- g\Bigg(\sum\limits_{\bm{j}=\bm{1}}^{\bm{k}}\bm{w}^{\times}_{\bm{j}}\bm{V}^{\ast}\mbox{Poly}_{f,\Psi_{\bm{j}},\mathscr{L}}(\underline{\bm{A}}_{j_1})\bm{V},\ldots,\sum\limits_{\bm{j}=\bm{1}}^{\bm{k}}\bm{w}^{\times}_{\bm{j}}\bm{V}^{\ast}\mbox{Poly}_{f,\Psi_{\bm{j}},\mathscr{L}}(\underline{\bm{A}}_{j_n})\bm{V}\Bigg)\nonumber \\
&\geq& 
\beta_2 (x - g(\bm{x}))\bm{I}_{\mathfrak{K}},
\end{eqnarray}
where $\beta_2$ is some specfied real number. 

The following Example~\ref{exp: cor 2.15} is provided to evaluate the upper difference bound given by Eq.~\eqref{eq3: UP Difference Type Approximation} and evaluate the lower difference bound given by Eq.~\eqref{eq3: LO Difference Type Approximation}. 
\begin{example}\label{exp: cor 2.15}
Let $\bm{A}_{j_i}$ be self-adjoint operators with $\Lambda(\bm{A}_{j_i}) \in [m_i, M_i]$ for real scalars $m_i <  M_i$. The mappings $\Phi_{j_1,\ldots,j_n}: \mathscr{B}(\mathfrak{H}) \rightarrow \mathscr{B}(\mathfrak{K})$ are normalized positive linear maps, where $j_i=1,2,\ldots,k_i$ for $i=1,2,\ldots,n$. We have $n$ probability vectors $\bm{w}_i =[w_{i,1},w_{i,2},\cdots, w_{i,k_i}]$ with the dimension $k_i$ for $i=1,2,\ldots,n$, i.e., $\sum\limits_{\ell=1}^{k_i}w     _{i,\ell} = 1$. Let $f(x_1,x_2,\ldots,x_n)$ be any real valued continuous functions with $n$ variables defined on the range $\bigtimes\limits_{i=1}^n [m_i, M_i] \in \mathbb{R}^n$, where $\times$ is the Cartesian product. Besides, we assume that the function $f$ satisfies the following:
\begin{eqnarray}\label{eq f bounds: cor: 2.17 ext}
\sum\limits_{i=1}^n a_i x_i + b \leq f(x_1,x_2,\ldots,x_n) \leq \sum\limits_{i=1}^n c_i x_i + d,
\end{eqnarray}
for $[x_1,x_2,\ldots,x_n] \in \bigtimes\limits_{i=1}^n [m_i,M_i]$. We will discuss three cases. 

(I) If we also assume that $g(x_1,\ldots,x_n) = \Big(\sum\limits_{i=1}^n \beta_i x_i\Big)^q$, where $\beta_i \geq 0$, $q \in \mathbb{R}$ and $m_i \geq 0$, then, we have the following upper bound for $\sum\limits_{j_{1}=1,\ldots,j_{n}=1}^{k_1,\ldots,k_n}  w_{j_1}\ldots w_{j_n}\Phi_{j_1,\ldots,j_n}(f(\bm{A}_{j_1},\bm{A}_{j_2},\ldots,\bm{A}_{j_n}))$:
\begin{eqnarray}\label{eq x q UB: cor: 2.17 ext}
\lefteqn{\sum\limits_{j_{1}=1,\ldots,j_{n}=1}^{k_1,\ldots,k_n}  w_{j_1}\ldots w_{j_n}\Phi_{j_1,\ldots,j_n}(f(\bm{A}_{j_1},\bm{A}_{j_2},\ldots,\bm{A}_{j_n}))}\nonumber \\
&&- \left(\sum\limits_{i=1}^n \beta_i \sum\limits_{j_{1}=1,\ldots,j_{n}=1}^{k_1,\ldots,k_n}w_{j_1}\ldots w_{j_n}\Phi_{j_1,\ldots,j_n}(\bm{A}_{j_i})\right)^q
\nonumber \\
&\leq&
\max\limits_{\substack{m_i \leq x_i \leq M_i \\ i=1,2,\ldots,n}}\Big(\sum\limits_{i=1}^n c_i x_i + d - \Big(\sum\limits_{i=1}^n \beta_i x_i\Big)^q\Big)\bm{I}_{\mathfrak{K}};
\end{eqnarray}
and, the following lower bound for $\sum\limits_{j_{1}=1,\ldots,j_{n}=1}^{k_1,\ldots,k_n}  w_{j_1}\ldots w_{j_n}\Phi_{j_1,\ldots,j_n}(f(\bm{A}_{j_1},\bm{A}_{j_2},\ldots,\bm{A}_{j_n}))$:
\begin{eqnarray}\label{eq x q LB: cor: 2.17 ext}
\lefteqn{\sum\limits_{j_{1}=1,\ldots,j_{n}=1}^{k_1,\ldots,k_n}  w_{j_1}\ldots w_{j_n}\Phi_{j_1,\ldots,j_n}(f(\bm{A}_{j_1},\bm{A}_{j_2},\ldots,\bm{A}_{j_n}))}\nonumber \\
&&- \left(\sum\limits_{i=1}^n \beta_i  \sum\limits_{j_{1}=1,\ldots,j_{n}=1}^{k_1,\ldots,k_n}w_{j_1}\ldots w_{j_n}\Phi_{j_1,\ldots,j_n}(\bm{A}_{j_i})\right)^q
\nonumber \\
&\geq&
\min\limits_{\substack{m_i \leq x_i \leq M_i \\ i=1,2,\ldots,n}}\Big(\sum\limits_{i=1}^n a_i x_i + b - \Big(\sum\limits_{i=1}^n \beta_i x_i\Big)^q\Big)\bm{I}_{\mathfrak{K}}.
\end{eqnarray}

For Part (I), we will use $g(x_1,\ldots,x_n) = \Big(\sum\limits_{i=1}^n \beta_i x_i\Big)^q$ in Eq.~\eqref{eq UB: cor 2.15} in Theorem~\ref{thm: cor 2.15} to obtain Eq.~\eqref{eq x q UB: cor: 2.17 ext}. We will also use $g(\bm{x})=\Big(\sum\limits_{i=1}^n \beta_i x_i\Big)^q$ in Eq.~\eqref{eq LB: cor 2.15} in Theorem~\ref{thm: cor 2.15} to obtain Eq.~\eqref{eq x q LB: cor: 2.17 ext}.  Then, we have
\begin{eqnarray}\label{eq1: exp: cor 2.15}
\beta_1 &\geq& \max\limits_{\substack{m_i \leq x_i \leq M_i \\ i=1,2,\ldots,n}}\Big(\sum\limits_{i=1}^n c_i x_i + d - \Big(\sum\limits_{i=1}^n \beta_i x_i\Big)^q\Big)\bm{I}_{\mathfrak{K}}, \nonumber \\
\beta_2 &\leq& \min\limits_{\substack{m_i \leq x_i \leq M_i \\ i=1,2,\ldots,n}}\Big(\sum\limits_{i=1}^n a_i x_i + b - \Big(\sum\limits_{i=1}^n \beta_i x_i\Big)^q\Big)\bm{I}_{\mathfrak{K}}.
\end{eqnarray}

(II) If $g(x_1,\ldots,x_n) = \log\Big(\sum\limits_{i=1}^n \beta_i x_i\Big)$ with $\beta_i \geq 0$ and $m_i >0$, we have the upper bound for $\sum\limits_{j_{1}=1,\ldots,j_{n}=1}^{k_1,\ldots,k_n}  w_{j_1}\ldots w_{j_n}\Phi_{j_1,\ldots,j_n}(f(\bm{A}_{j_1},\bm{A}_{j_2},\ldots,\bm{A}_{j_n}))$:
\begin{eqnarray}\label{eq log x UB: cor: 2.17 ext}
\lefteqn{\sum\limits_{j_{1}=1,\ldots,j_{n}=1}^{k_1,\ldots,k_n}  w_{j_1}\ldots w_{j_n}\Phi_{j_1,\ldots,j_n}(f(\bm{A}_{j_1},\bm{A}_{j_2},\ldots,\bm{A}_{j_n}))}\nonumber \\
&&- \log\left(\sum\limits_{i=1}^n \beta_i \sum\limits_{j_{1}=1,\ldots,j_{n}=1}^{k_1,\ldots,k_n}w_{j_1}\ldots w_{j_n}\Phi_{j_1,\ldots,j_n}(\bm{A}_{j_i})\right)
\nonumber \\
&\leq&
\max\limits_{\substack{m_i \leq x_i \leq M_i \\ i=1,2,\ldots,n}}\Big(\sum\limits_{i=1}^n c_i x_i + d - \log\Big(\sum\limits_{i=1}^n \beta_i x_i\Big)\Big)\bm{I}_{\mathfrak{K}};
\end{eqnarray}
and, the following lower bound for $\sum\limits_{j_{1}=1,\ldots,j_{n}=1}^{k_1,\ldots,k_n}  w_{j_1}\ldots w_{j_n}\Phi_{j_1,\ldots,j_n}(f(\bm{A}_{j_1},\bm{A}_{j_2},\ldots,\bm{A}_{j_n}))$:
\begin{eqnarray}\label{eq log x LB: cor: 2.17 ext}
\lefteqn{\sum\limits_{j_{1}=1,\ldots,j_{n}=1}^{k_1,\ldots,k_n}  w_{j_1}\ldots w_{j_n}\Phi_{j_1,\ldots,j_n}(f(\bm{A}_{j_1},\bm{A}_{j_2},\ldots,\bm{A}_{j_n}))}\nonumber \\
&&- \log\left(\sum\limits_{i=1}^n \beta_i \sum\limits_{j_{1}=1,\ldots,j_{n}=1}^{k_1,\ldots,k_n}w_{j_1}\ldots w_{j_n}\Phi_{j_1,\ldots,j_n}(\bm{A}_{j_i})\right)
\nonumber \\
&\geq&
\min\limits_{\substack{m_i \leq x_i \leq M_i \\ i=1,2,\ldots,n}}\Big(\sum\limits_{i=1}^n a_i x_i + b - \log\Big(\sum\limits_{i=1}^n \beta_i x_i\Big)\Big)\bm{I}_{\mathfrak{K}}.
\end{eqnarray}

For Part (II), we will use  $g(x_1,\ldots,x_n) = \log\Big(\sum\limits_{i=1}^n \beta_i x_i\Big)$ in Eq.~\eqref{eq UB: cor 2.15} in Theorem~\ref{thm: cor 2.15} to obtain Eq.~\eqref{eq log x UB: cor: 2.17 ext}. We will also use $g(\bm{x})=\log\Big(\sum\limits_{i=1}^n \beta_i x_i\Big)$ in Eq.~\eqref{eq LB: cor 2.15} in Theorem~\ref{thm: cor 2.15} to obtain Eq.~\eqref{eq log x LB: cor: 2.17 ext}. Then, we have
\begin{eqnarray}\label{eq2: exp: cor 2.15}
\beta_1 &\geq& \max\limits_{\substack{m_i \leq x_i \leq M_i \\ i=1,2,\ldots,n}}\Big(\sum\limits_{i=1}^n c_i x_i + d - \log\Big(\sum\limits_{i=1}^n \beta_i x_i\Big)\Big)\bm{I}_{\mathfrak{K}}, \nonumber \\
\beta_2 &\leq& \min\limits_{\substack{m_i \leq x_i \leq M_i \\ i=1,2,\ldots,n}}\Big(\sum\limits_{i=1}^n a_i x_i + b - \log\Big(\sum\limits_{i=1}^n \beta_i x_i\Big)\Big)\bm{I}_{\mathfrak{K}}.
\end{eqnarray}

(III) If $g(x_1,\ldots,x_n) = \exp\Big(\sum\limits_{i=1}^n \beta_i x_i\Big)$, we have the upper bound for \\
$\sum\limits_{j_{1}=1,\ldots,j_{n}=1}^{k_1,\ldots,k_n}  w_{j_1}\ldots w_{j_n}\Phi_{j_1,\ldots,j_n}(f(\bm{A}_{j_1},\bm{A}_{j_2},\ldots,\bm{A}_{j_n}))$:
\begin{eqnarray}\label{eq exp x UB: cor: 2.17 ext}
\lefteqn{\sum\limits_{j_{1}=1,\ldots,j_{n}=1}^{k_1,\ldots,k_n}  w_{j_1}\ldots w_{j_n}\Phi_{j_1,\ldots,j_n}(f(\bm{A}_{j_1},\bm{A}_{j_2},\ldots,\bm{A}_{j_n}))}\nonumber \\
&&- \exp\left(\sum\limits_{i=1}^n \beta_i \sum\limits_{j_{1}=1,\ldots,j_{n}=1}^{k_1,\ldots,k_n}w_{j_1}\ldots w_{j_n}\Phi_{j_1,\ldots,j_n}(\bm{A}_{j_i})\right)
\nonumber \\
&\leq&
\max\limits_{\substack{m_i \leq x_i \leq M_i \\ i=1,2,\ldots,n}}\Big(\sum\limits_{i=1}^n c_i x_i + d - \exp\Big(\sum\limits_{i=1}^n \beta_i x_i\Big)\Big)\bm{I}_{\mathfrak{K}};
\end{eqnarray}
and, the following lower bound for $\sum\limits_{j_{1}=1,\ldots,j_{n}=1}^{k_1,\ldots,k_n}  w_{j_1}\ldots w_{j_n}\Phi_{j_1,\ldots,j_n}(f(\bm{A}_{j_1},\bm{A}_{j_2},\ldots,\bm{A}_{j_n}))$:
\begin{eqnarray}\label{eq exp x LB: cor: 2.17 ext}
\lefteqn{\sum\limits_{j_{1}=1,\ldots,j_{n}=1}^{k_1,\ldots,k_n}  w_{j_1}\ldots w_{j_n}\Phi_{j_1,\ldots,j_n}(f(\bm{A}_{j_1},\bm{A}_{j_2},\ldots,\bm{A}_{j_n}))}\nonumber \\
&&- \exp\left(\sum\limits_{i=1}^n \beta_i \sum\limits_{j_{1}=1,\ldots,j_{n}=1}^{k_1,\ldots,k_n}w_{j_1}\ldots w_{j_n}\Phi_{j_1,\ldots,j_n}(\bm{A}_{j_i})\right)
\nonumber \\
&\geq&
\min\limits_{\substack{m_i \leq x_i \leq M_i \\ i=1,2,\ldots,n}}\Big(\sum\limits_{i=1}^n a_i x_i + b - \exp\Big(\sum\limits_{i=1}^n \beta_i x_i\Big)\Big)\bm{I}_{\mathfrak{K}}.
\end{eqnarray}

For Part (III), we will use $g(x_1,\ldots,x_n) = \exp\Big(\sum\limits_{i=1}^n \beta_i x_i\Big)$ in Eq.~\eqref{eq UB: cor 2.15} in Theorem~\ref{thm: cor 2.15} to obtain Eq.~\eqref{eq exp x UB: cor: 2.17 ext}. We will also use $g(\bm{x})=\exp\Big(\sum\limits_{i=1}^n \beta_i x_i\Big)$ in Eq.~\eqref{eq LB: cor 2.15} in Theorem~\ref{thm: cor 2.15} to obtain Eq.~\eqref{eq exp x LB: cor: 2.17 ext}. Then, we have
\begin{eqnarray}\label{eq2: exp: cor 2.15}
\beta_1 &\geq& \max\limits_{\substack{m_i \leq x_i \leq M_i \\ i=1,2,\ldots,n}}\Big(\sum\limits_{i=1}^n c_i x_i + d - \exp\Big(\sum\limits_{i=1}^n \beta_i x_i\Big)\Big)\bm{I}_{\mathfrak{K}}, \nonumber \\
\beta_2 &\leq& \min\limits_{\substack{m_i \leq x_i \leq M_i \\ i=1,2,\ldots,n}}\Big(\sum\limits_{i=1}^n a_i x_i + b - \exp\Big(\sum\limits_{i=1}^n \beta_i x_i\Big)\Big)\bm{I}_{\mathfrak{K}}.
\end{eqnarray}
\end{example}

\begin{remark}\label{rmk: difference approx conjecture}
Given requirements of Eq.~\eqref{eq1: UP Difference Type Approximation}, Eq.~\eqref{eq1: LO Difference Type Approximation}, Eq.~\eqref{eq3: UP Difference Type Approximation}, and Eq.~\eqref{eq3: LO Difference Type Approximation}, we conjecture the existence of function $g(\bm{x})$ and polynomials in terms of $\Psi$as arguments $\mbox{Poly}_{f,\Psi,\mathscr{U}}(\bm{x})$, $\mbox{Poly}_{f,\Psi_{\bm{j}},\mathscr{U}}(\bm{x})$ and $\mbox{Poly}_{f,\Psi,\mathscr{L}}(\bm{x})$, $\mbox{Poly}_{f,\Psi_{\bm{j}}, \mathscr{L}}(\bm{x})$ for any given $\beta_1$ or $\beta_2$. If existence, how to find these functions?
\end{remark}

\section{Bounds Algebra}\label{sec: Bounds Algebra}

\subsection{Multivariate Hypercomplex Function Bounds Algebra}\label{sec: Multivariate Function of Operator Bounds Algebra}

Let $\bm{A}_{j_i}$ be self-adjoint operators with $\Lambda(\bm{A}_{j_i}) \in [m_i, M_i]$ for real scalars $m_i <  M_i$. The mappings $\Phi_{j_1,\ldots,j_n}: \mathscr{B}(\mathfrak{H}) \rightarrow \mathscr{B}(\mathfrak{K})$ are defined by Eq.~\eqref{eq: new phi def}, where $j_i=1,2,\ldots,k_i$ for $i=1,2,\ldots,n$. We have $n$ probability vectors $\bm{w}_i =[w_{i,1},w_{i,2},\cdots, w_{i,k_i}]$ with the dimension $k_i$ for $i=1,2,\ldots,n$, i.e., $\sum\limits_{\ell=1}^{k_i}w_{i,\ell} = 1$. Let $f(\bm{x})$ be any real valued continuous functions with $n$ variables defined on the range $\bigtimes\limits_{i=1}^n [m_i, M_i] \in \mathbb{R}^n$, where $\bigtimes$ is the Cartesian product. Besides, given any $\epsilon>0$, we assume that the function $f(\bm{x})=f(x_1,x_2,\ldots,x_n)$ satisfies the following:
\begin{eqnarray}\label{eq1: sec: Function of Operator Bounds Algebra}
\sum\limits_{i=1}^n a_{f,i} x_i + b_f \leq f(x_1,x_2,\ldots,x_n) \leq \sum\limits_{i=1}^n c_{f,i} x_i + d_f,
\end{eqnarray}
for $\bm{x}\in \bm{I}^{\bigtimes}$, and another function $h(\bm{x})=h(x_1,x_2,\ldots,x_n)$ also satisfies the following:
\begin{eqnarray}\label{eq2: sec: Function of Operator Bounds Algebra}
\sum\limits_{i=1}^n a_{h,i} x_i + b_h \leq h(x_1,x_2,\ldots,x_n) \leq \sum\limits_{i=1}^n c_{h,i} x_i + d_h,
\end{eqnarray}
for $\bm{x}\in \bm{I}^{\bigtimes}$. If we require $g(\bm{x})>0$ for all $\bm{x} \in \bigtimes\limits_{i=1}^n [m_i, M_i]$, from Theorem~\ref{thm: 2.9}, we have

\begin{eqnarray}\label{eq2f: UP sec: Function of Operator Bounds Algebra}
\lefteqn{\sum\limits_{\bm{j}=\bm{1}}^{\bm{k}}\bm{w}^{\times}_{\bm{j}}\Phi_{\bm{j}}(f(\underline{\bm{A}}_{\bm{j}}))}\nonumber \\
&\leq&\underbrace{\left[\max\limits_{\substack{m_i \leq x_i \leq M_i \\ i=1,2,\ldots,n}}\frac{\sum\limits_{i=1}^n c_{f,i} x_i +d_f}{g(x_1,\ldots,x_n)}\right]}_{:=\alpha_{f,\mathscr{U}}}g\Bigg(\sum\limits_{\bm{j}=\bm{1}}^{\bm{k}}\bm{w}^{\times}_{\bm{j}}\Phi_{\bm{j}}(\bm{A}_{j_1}),\ldots,\sum\limits_{\bm{j}=\bm{1}}^{\bm{k}}\bm{w}^{\times}_{\bm{j}}\Phi_{\bm{j}}(\bm{A}_{j_n})\Bigg),
\end{eqnarray}
and 
\begin{eqnarray}\label{eq2f: LO sec: Function of Operator Bounds Algebra}
\lefteqn{\sum\limits_{j=1}^k w_j \Phi(f(\bm{A}_j))}\nonumber \\
&\geq& \underbrace{\left[\min\limits_{\substack{m_i \leq x_i \leq M_i \\ i=1,2,\ldots,n}}\frac{\sum\limits_{i=1}^n a_{f,i} x_i +b_f}{g(x_1,\ldots,x_n)}\right]}_{:=\alpha_{f,\mathscr{L}}}g\Bigg(\sum\limits_{\bm{j}=\bm{1}}^{\bm{k}}\bm{w}^{\times}_{\bm{j}}\Phi_{\bm{j}}(\bm{A}_{j_1}),\ldots,\sum\limits_{\bm{j}=\bm{1}}^{\bm{k}}\bm{w}^{\times}_{\bm{j}}\Phi_{\bm{j}}(\bm{A}_{j_n})\Bigg).
\end{eqnarray}
Similarly,  from Theorem~\ref{thm: 2.9}, we also have
\begin{eqnarray}\label{eq2h: UP sec: Function of Operator Bounds Algebra}
\lefteqn{\sum\limits_{\bm{j}=\bm{1}}^{\bm{k}}\bm{w}^{\times}_{\bm{j}}\Phi_{\bm{j}}(h(\underline{\bm{A}}_{\bm{j}}))}\nonumber \\
&\leq&\underbrace{\left[\max\limits_{\substack{m_i \leq x_i \leq M_i \\ i=1,2,\ldots,n}}\frac{\sum\limits_{i=1}^n c_{h,i} x_i +d_h}{g(x_1,\ldots,x_n)}\right]}_{:=\alpha_{h,\mathscr{U}}}g\Bigg(\sum\limits_{\bm{j}=\bm{1}}^{\bm{k}}\bm{w}^{\times}_{\bm{j}}\Phi_{\bm{j}}(\bm{A}_{j_1}),\ldots,\sum\limits_{\bm{j}=\bm{1}}^{\bm{k}}\bm{w}^{\times}_{\bm{j}}\Phi_{\bm{j}}(\bm{A}_{j_n})\Bigg),
\end{eqnarray}
and 
\begin{eqnarray}\label{eq2h: LO sec: Function of Operator Bounds Algebra}
\lefteqn{\sum\limits_{\bm{j}=\bm{1}}^{\bm{k}}\bm{w}^{\times}_{\bm{j}}\Phi_{\bm{j}}(h(\underline{\bm{A}}_{\bm{j}}))}\nonumber \\
&\geq& \underbrace{\left[\min\limits_{\substack{m_i \leq x_i \leq M_i \\ i=1,2,\ldots,n}}\frac{\sum\limits_{i=1}^n a_{h,i} x_i +b_h}{g(x_1,\ldots,x_n)}\right]}_{:=\alpha_{h,\mathscr{L}}}g\Bigg(\sum\limits_{\bm{j}=\bm{1}}^{\bm{k}}\bm{w}^{\times}_{\bm{j}}\Phi_{\bm{j}}(\bm{A}_{j_1}),\ldots,\sum\limits_{\bm{j}=\bm{1}}^{\bm{k}}\bm{w}^{\times}_{\bm{j}}\Phi_{\bm{j}}(\bm{A}_{j_n})\Bigg).
\end{eqnarray}

Consider the addition between $\sum\limits_{\bm{j}=\bm{1}}^{\bm{k}}\bm{w}^{\times}_{\bm{j}}\Phi_{\bm{j}}(f(\underline{\bm{A}}_{\bm{j}}))$ and $\sum\limits_{\bm{j}=\bm{1}}^{\bm{k}}\bm{w}^{\times}_{\bm{j}}\Phi_{\bm{j}}(h(\underline{\bm{A}}_{\bm{j}}))$, from Eq.~\eqref{eq2f: UP sec: Function of Operator Bounds Algebra} and Eq.~\eqref{eq2h: UP sec: Function of Operator Bounds Algebra}, then, we have
\begin{eqnarray}\label{eq2f h: UP sec: Function of Operator Bounds Algebra} 
\lefteqn{\sum\limits_{\bm{j}=\bm{1}}^{\bm{k}}\bm{w}^{\times}_{\bm{j}}\Phi_{\bm{j}}(f(\underline{\bm{A}}_{\bm{j}}))+\sum\limits_{\bm{j}=\bm{1}}^{\bm{k}}\bm{w}^{\times}_{\bm{j}}\Phi_{\bm{j}}(h(\underline{\bm{A}}_{\bm{j}}))}
\nonumber \\
&\leq&(\alpha_{f,\mathscr{U}} + \alpha_{h,\mathscr{U}})g\Bigg(\sum\limits_{\bm{j}=\bm{1}}^{\bm{k}}\bm{w}^{\times}_{\bm{j}}\Phi_{\bm{j}}(\bm{A}_{j_1}),\ldots,\sum\limits_{\bm{j}=\bm{1}}^{\bm{k}}\bm{w}^{\times}_{\bm{j}}\Phi_{\bm{j}}(\bm{A}_{j_n})\Bigg),
\end{eqnarray}
and, from Eq.~\eqref{eq2f: LO sec: Function of Operator Bounds Algebra} and Eq.~\eqref{eq2h: LO sec: Function of Operator Bounds Algebra}, we also have
\begin{eqnarray}\label{eq2f h: LO sec: Function of Operator Bounds Algebra} 
\lefteqn{\sum\limits_{\bm{j}=\bm{1}}^{\bm{k}}\bm{w}^{\times}_{\bm{j}}\Phi_{\bm{j}}(f(\underline{\bm{A}}_{\bm{j}}))+\sum\limits_{\bm{j}=\bm{1}}^{\bm{k}}\bm{w}^{\times}_{\bm{j}}\Phi_{\bm{j}}(h(\underline{\bm{A}}_{\bm{j}}))}
\nonumber \\
&\geq&(\alpha_{f,\mathscr{L}} + \alpha_{h,\mathscr{L}})g\Bigg(\sum\limits_{\bm{j}=\bm{1}}^{\bm{k}}\bm{w}^{\times}_{\bm{j}}\Phi_{\bm{j}}(\bm{A}_{j_1}),\ldots,\sum\limits_{\bm{j}=\bm{1}}^{\bm{k}}\bm{w}^{\times}_{\bm{j}}\Phi_{\bm{j}}(\bm{A}_{j_n})\Bigg).
\end{eqnarray}

If any two functions $f,h$ satisfying Eq.~\eqref{eq1: sec: Function of Operator Bounds Algebra} and Eq.~\eqref{eq2: sec: Function of Operator Bounds Algebra}, from Eq.~\eqref{eq2f h: UP sec: Function of Operator Bounds Algebra} and Eq.~\eqref{eq2f h: LO sec: Function of Operator Bounds Algebra}, we can bound the addition of $\sum\limits_{\bm{j}=\bm{1}}^{\bm{k}}\bm{w}^{\times}_{\bm{j}}\Phi_{\bm{j}}(f(\underline{\bm{A}}_{\bm{j}}))$ and $\sum\limits_{\bm{j}=\bm{1}}^{\bm{k}}\bm{w}^{\times}_{\bm{j}}\Phi_{\bm{j}}(h(\underline{\bm{A}}_{\bm{j}}))$ by coefficients of $g\left(\sum\limits_{j=1}^k w_j \Phi(\bm{A}_j)\right)$S, which are independent of the functions $f$ and $h$. Therefore, the ranges of bounding coefficients of \\
 $g\Bigg(\sum\limits_{\bm{j}=\bm{1}}^{\bm{k}}\bm{w}^{\times}_{\bm{j}}\Phi_{\bm{j}}(\bm{A}_{j_1}),\ldots,\sum\limits_{\bm{j}=\bm{1}}^{\bm{k}}\bm{w}^{\times}_{\bm{j}}\Phi_{\bm{j}}(\bm{A}_{j_n})\Bigg)$ form an \emph{algebraic system of interval numbers}, which form an \emph{Abelian monoid} with respect to the opration \emph{addition} accoding to Theorem 2.14~\cite{dawood2011theories}. 

Consider the multiplication between $\sum\limits_{\bm{j}=\bm{1}}^{\bm{k}}\bm{w}^{\times}_{\bm{j}}\Phi_{\bm{j}}(f(\underline{\bm{A}}_{\bm{j}}))$ and $\sum\limits_{\bm{j}=\bm{1}}^{\bm{k}}\bm{w}^{\times}_{\bm{j}}\Phi_{\bm{j}}(h(\underline{\bm{A}}_{\bm{j}}))$~\footnote{We assume that the multiplication of the self-adjoint $\sum\limits_{j=1}^k w_j \Phi(f(\bm{A}_j))$ operator and the self-adjoint $\sum\limits_{j=1}^k w_j \Phi(h(\bm{A}_j))$ operator is still a self-adjoint operator.}, from Eq.~\eqref{eq2f: UP sec: Function of Operator Bounds Algebra} and Eq.~\eqref{eq2h: UP sec: Function of Operator Bounds Algebra} with assumptions of positive $\alpha_{f,\mathscr{U}}, \alpha_{h,\mathscr{U}}, \alpha_{f,\mathscr{L}}$ and $\alpha_{h,\mathscr{L}}$ and $g(\bm{x})>0$ for all $\bm{x} \in \bigtimes\limits_{i=1}^n [m_i, M_i]$, then, we have
\begin{eqnarray}\label{eq3f h: UP sec: Function of Operator Bounds Algebra} 
\lefteqn{\left(\sum\limits_{\bm{j}=\bm{1}}^{\bm{k}}\bm{w}^{\times}_{\bm{j}}\Phi_{\bm{j}}(f(\underline{\bm{A}}_{\bm{j}}))\right)\times\left(\sum\limits_{\bm{j}=\bm{1}}^{\bm{k}}\bm{w}^{\times}_{\bm{j}}\Phi_{\bm{j}}(h(\underline{\bm{A}}_{\bm{j}}))\right)}\nonumber \\
&\leq&(\alpha_{f,\mathscr{U}} \times \alpha_{h,\mathscr{U}})\left(g\Bigg(\sum\limits_{\bm{j}=\bm{1}}^{\bm{k}}\bm{w}^{\times}_{\bm{j}}\Phi_{\bm{j}}(\bm{A}_{j_1}),\ldots,\sum\limits_{\bm{j}=\bm{1}}^{\bm{k}}\bm{w}^{\times}_{\bm{j}}\Phi_{\bm{j}}(\bm{A}_{j_n})\Bigg)\right)^2,
\end{eqnarray}
and, from Eq.~\eqref{eq2f: LO sec: Function of Operator Bounds Algebra} and Eq.~\eqref{eq2h: LO sec: Function of Operator Bounds Algebra}, we also have
\begin{eqnarray}\label{eq3f h: LO sec: Function of Operator Bounds Algebra} 
\lefteqn{\left(\sum\limits_{\bm{j}=\bm{1}}^{\bm{k}}\bm{w}^{\times}_{\bm{j}}\Phi_{\bm{j}}(f(\underline{\bm{A}}_{\bm{j}}))\right)\times\left(\sum\limits_{\bm{j}=\bm{1}}^{\bm{k}}\bm{w}^{\times}_{\bm{j}}\Phi_{\bm{j}}(h(\underline{\bm{A}}_{\bm{j}}))\right)}\nonumber \\
&\geq&(\alpha_{f,\mathscr{L}} \times \alpha_{h,\mathscr{L}})\left(g\Bigg(\sum\limits_{\bm{j}=\bm{1}}^{\bm{k}}\bm{w}^{\times}_{\bm{j}}\Phi_{\bm{j}}(\bm{A}_{j_1}),\ldots,\sum\limits_{\bm{j}=\bm{1}}^{\bm{k}}\bm{w}^{\times}_{\bm{j}}\Phi_{\bm{j}}(\bm{A}_{j_n})\Bigg)\right)^2.
\end{eqnarray}

Analogly, if any two functions $f,h$ satisfying Eq.~\eqref{eq1: sec: Function of Operator Bounds Algebra} and Eq.~\eqref{eq2: sec: Function of Operator Bounds Algebra}, from Eq.~\eqref{eq3f h: UP sec: Function of Operator Bounds Algebra} and Eq.~\eqref{eq3f h: LO sec: Function of Operator Bounds Algebra}, we can bound the multiplication between $\sum\limits_{\bm{j}=\bm{1}}^{\bm{k}}\bm{w}^{\times}_{\bm{j}}\Phi_{\bm{j}}(f(\underline{\bm{A}}_{\bm{j}}))$ and $\sum\limits_{\bm{j}=\bm{1}}^{\bm{k}}\bm{w}^{\times}_{\bm{j}}\Phi_{\bm{j}}(h(\underline{\bm{A}}_{\bm{j}}))$ by coefficients of $g\Bigg(\sum\limits_{\bm{j}=\bm{1}}^{\bm{k}}\bm{w}^{\times}_{\bm{j}}\Phi_{\bm{j}}(\bm{A}_{j_1}),\ldots,\sum\limits_{\bm{j}=\bm{1}}^{\bm{k}}\bm{w}^{\times}_{\bm{j}}\Phi_{\bm{j}}(\bm{A}_{j_n})\Bigg)$, which are independent of the functions $f$ and $h$. Therefore, the ranges of bounding coefficients of $g\Bigg(\sum\limits_{\bm{j}=\bm{1}}^{\bm{k}}\bm{w}^{\times}_{\bm{j}}\Phi_{\bm{j}}(\bm{A}_{j_1}),\ldots,\sum\limits_{\bm{j}=\bm{1}}^{\bm{k}}\bm{w}^{\times}_{\bm{j}}\Phi_{\bm{j}}(\bm{A}_{j_n})\Bigg)$, which form an \emph{Abelian monoid} with respect to the opration \emph{multiplication} accoding to Theorem 2.14~\cite{dawood2011theories}. 

\begin{remark}\label{rmk: norm bounds}
The norms of the addition (or multiplication) between $\sum\limits_{\bm{j}=\bm{1}}^{\bm{k}}\bm{w}^{\times}_{\bm{j}}\Phi_{\bm{j}}(f(\underline{\bm{A}}_{\bm{j}}))$ and $\sum\limits_{\bm{j}=\bm{1}}^{\bm{k}}\bm{w}^{\times}_{\bm{j}}\Phi_{\bm{j}}(h(\underline{\bm{A}}_{\bm{j}}))$ can also be bounded by coefficients of $g\Bigg(\sum\limits_{\bm{j}=\bm{1}}^{\bm{k}}\bm{w}^{\times}_{\bm{j}}\Phi_{\bm{j}}(\bm{A}_{j_1}),\ldots,\sum\limits_{\bm{j}=\bm{1}}^{\bm{k}}\bm{w}^{\times}_{\bm{j}}\Phi_{\bm{j}}(\bm{A}_{j_n})\Bigg)$\\
(or $\left(g\Bigg(\sum\limits_{\bm{j}=\bm{1}}^{\bm{k}}\bm{w}^{\times}_{\bm{j}}\Phi_{\bm{j}}(\bm{A}_{j_1}),\ldots,\sum\limits_{\bm{j}=\bm{1}}^{\bm{k}}\bm{w}^{\times}_{\bm{j}}\Phi_{\bm{j}}(\bm{A}_{j_n})\Bigg)\right)^2$) with Abelian monoid algebraic structure and norms of $g\Bigg(\sum\limits_{\bm{j}=\bm{1}}^{\bm{k}}\bm{w}^{\times}_{\bm{j}}\Phi_{\bm{j}}(\bm{A}_{j_1}),\ldots,\sum\limits_{\bm{j}=\bm{1}}^{\bm{k}}\bm{w}^{\times}_{\bm{j}}\Phi_{\bm{j}}(\bm{A}_{j_n})\Bigg)$ (or $\left(g\Bigg(\sum\limits_{\bm{j}=\bm{1}}^{\bm{k}}\bm{w}^{\times}_{\bm{j}}\Phi_{\bm{j}}(\bm{A}_{j_1}),\ldots,\sum\limits_{\bm{j}=\bm{1}}^{\bm{k}}\bm{w}^{\times}_{\bm{j}}\Phi_{\bm{j}}(\bm{A}_{j_n})\Bigg)\right)^2$).
\end{remark}

\subsection{Multivariate Random Tensors Tail Bounds Algebra}\label{sec: Multivariate Random Tensors Tail Bounds Algebra}

Let $\bm{A}_{j_i}$ be random Hermitian tensors with $\Lambda(\bm{A}_{j_i}) \in [m_i, M_i]$ for real scalars $m_i <  M_i$. We have $n$ probability vectors $\bm{w}_i =[w_{i,1},w_{i,2},\cdots, w_{i,k_i}]$ with the dimension $k_i$ for $i=1,2,\ldots,n$, i.e., $\sum\limits_{\ell=1}^{k_i}w_{i,\ell} = 1$. Let $f(\bm{x})$ be any real valued continuous functions with $n$ variables defined on the range $\bigtimes\limits_{i=1}^n [m_i, M_i] \in \mathbb{R}^n$, where $\bigtimes$ is the Cartesian product. Besides, given any $\epsilon>0$, we assume that the function $f(\bm{x})=f(x_1,x_2,\ldots,x_n)$ satisfies the following:
\begin{eqnarray}\label{eq1: sec: Multivariate Random Tensors Tail Bounds Algebra}
f(x_1,x_2,\ldots,x_n) \leq \sum\limits_{i=1}^n c_{f,i} x_i + d_f,
\end{eqnarray}
for $\bm{x}\in \bm{I}^{\bigtimes}$, and another function $h(\bm{x})=h(x_1,x_2,\ldots,x_n)$ also satisfies the following:
\begin{eqnarray}\label{eq2: sec: Multivariate Random Tensors Tail Bounds Algebra}
h(x_1,x_2,\ldots,x_n) \leq \sum\limits_{i=1}^n c_{h,i} x_i + d_h,
\end{eqnarray}
for $\bm{x}\in \bm{I}^{\bigtimes}$. From Theorem~\ref{thm: 2.9}, we have
\begin{eqnarray}\label{eq2f: UP sec: Random Tensor Tail Bounds Algebra}
\lefteqn{\sum\limits_{j=1}^k w_j \Phi(f(\bm{A}_j))}\nonumber \\
&\leq& \underbrace{\left[\min\limits_{\substack{m_i \leq x_i \leq M_i \\ i=1,2,\ldots,n}}\frac{\sum\limits_{i=1}^n c_{f,i} x_i +d_f}{g(x_1,\ldots,x_n)}\right]}_{:=\alpha_{f,\mathscr{L}}}g\Bigg(\sum\limits_{\bm{j}=\bm{1}}^{\bm{k}}\bm{w}^{\times}_{\bm{j}}\Phi_{\bm{j}}(\bm{A}_{j_1}),\ldots,\sum\limits_{\bm{j}=\bm{1}}^{\bm{k}}\bm{w}^{\times}_{\bm{j}}\Phi_{\bm{j}}(\bm{A}_{j_n})\Bigg).
\end{eqnarray}
Similarly,  from Theorem~\ref{thm: 2.9}, we also have
\begin{eqnarray}\label{eq2h: UP sec: Random Tensor Tail Bounds Algebra}
\lefteqn{\sum\limits_{\bm{j}=\bm{1}}^{\bm{k}}\bm{w}^{\times}_{\bm{j}}\Phi_{\bm{j}}(h(\underline{\bm{A}}_{\bm{j}}))}\nonumber \\
&\leq&\underbrace{\left[\max\limits_{\substack{m_i \leq x_i \leq M_i \\ i=1,2,\ldots,n}}\frac{\sum\limits_{i=1}^n c_{h,i} x_i +d_h}{g(x_1,\ldots,x_n)}\right]}_{:=\alpha_{h,\mathscr{U}}}g\Bigg(\sum\limits_{\bm{j}=\bm{1}}^{\bm{k}}\bm{w}^{\times}_{\bm{j}}\Phi_{\bm{j}}(\bm{A}_{j_1}),\ldots,\sum\limits_{\bm{j}=\bm{1}}^{\bm{k}}\bm{w}^{\times}_{\bm{j}}\Phi_{\bm{j}}(\bm{A}_{j_n})\Bigg).
\end{eqnarray}

Given $g(\bm{x})>0$ for all $\bm{x}\in \bm{I}^{\bigtimes}$ and any positive number $\theta$, we consider tail bounds for the addition between $\sum\limits_{\bm{j}=\bm{1}}^{\bm{k}}\bm{w}^{\times}_{\bm{j}}\Phi_{\bm{j}}(f(\underline{\bm{A}}_{\bm{j}}))$ and $\sum\limits_{\bm{j}=\bm{1}}^{\bm{k}}\bm{w}^{\times}_{\bm{j}}\Phi_{\bm{j}}(h(\underline{\bm{A}}_{\bm{j}}))$, from Eq.~\eqref{eq2f: UP sec: Random Tensor Tail Bounds Algebra} and Eq.~\eqref{eq2h: UP sec: Random Tensor Tail Bounds Algebra}, then, we have
\begin{eqnarray}\label{eq2f h: UP sec: Function of Operator Bounds Algebra}
\lefteqn{\mathrm{Pr}\left(\left\Vert\left(\sum\limits_{\bm{j}=\bm{1}}^{\bm{k}}\bm{w}^{\times}_{\bm{j}}\Phi_{\bm{j}}(f(\underline{\bm{A}}_{\bm{j}}))\right)+\left(\sum\limits_{\bm{j}=\bm{1}}^{\bm{k}}\bm{w}^{\times}_{\bm{j}}\Phi_{\bm{j}}(h(\underline{\bm{A}}_{\bm{j}}))\right)\right\Vert_{\ell}\geq\theta\right)}\nonumber \\
&\leq&\mathrm{Pr}\left( \left\Vert(\alpha_{f,\mathscr{U}} + \alpha_{h,\mathscr{U}})g\Bigg(\sum\limits_{\bm{j}=\bm{1}}^{\bm{k}}\bm{w}^{\times}_{\bm{j}}\Phi_{\bm{j}}(\bm{A}_{j_1}),\ldots,\sum\limits_{\bm{j}=\bm{1}}^{\bm{k}}\bm{w}^{\times}_{\bm{j}}\Phi_{\bm{j}}(\bm{A}_{j_n})\Bigg)\right\Vert_{\ell}\geq\theta \right),
\end{eqnarray}
where $\left\Vert \cdot \right\Vert_{(\ell)}$ is Ky Fan $\ell$-norm. R.H.S. of Eq.~\eqref{eq2f h: UP sec: Function of Operator Bounds Algebra}, where the random tensors summation part is independent of functions $f$ and $h$ can be upper bounded by those theorems in Section IV in~\cite{chang2022generalMA} with proper assumptions of the function $g$.

Consider the multiplication between $\sum\limits_{\bm{j}=\bm{1}}^{\bm{k}}\bm{w}^{\times}_{\bm{j}}\Phi_{\bm{j}}(f(\underline{\bm{A}}_{\bm{j}}))$ and $\sum\limits_{\bm{j}=\bm{1}}^{\bm{k}}\bm{w}^{\times}_{\bm{j}}\Phi_{\bm{j}}(h(\underline{\bm{A}}_{\bm{j}}))$~\footnote{We assume that the multiplication of the Hermtian tensor $\sum\limits_{j=1}^k w_j \Phi(f(\bm{A}_j))$ and the Hermtian tensor $\sum\limits_{j=1}^k w_j \Phi(h(\bm{A}_j))$ is still a Hermitian tensor.}, from Eq.~\eqref{eq2f: UP sec: Random Tensor Tail Bounds Algebra} and Eq.~\eqref{eq2h: UP sec: Random Tensor Tail Bounds Algebra} with assumptions of positive $\alpha_{f,\mathscr{U}}$ and $\alpha_{h,\mathscr{U}}$ and $g\Bigg(\sum\limits_{\bm{j}=\bm{1}}^{\bm{k}}\bm{w}^{\times}_{\bm{j}}\Phi_{\bm{j}}(\bm{A}_{j_1}),\ldots,\sum\limits_{\bm{j}=\bm{1}}^{\bm{k}}\bm{w}^{\times}_{\bm{j}}\Phi_{\bm{j}}(\bm{A}_{j_n})\Bigg) > \bm{0}$, then, we have
\begin{eqnarray}\label{eq3f h: UP sec: Random Tensor Tail Bounds Algebra} 
\lefteqn{\mathrm{Pr}\left(\left\Vert\left(\sum\limits_{\bm{j}=\bm{1}}^{\bm{k}}\bm{w}^{\times}_{\bm{j}}\Phi_{\bm{j}}(f(\underline{\bm{A}}_{\bm{j}}))\right)\times\left(\sum\limits_{\bm{j}=\bm{1}}^{\bm{k}}\bm{w}^{\times}_{\bm{j}}\Phi_{\bm{j}}(h(\underline{\bm{A}}_{\bm{j}}))\right)\right\Vert_{\ell}\geq \theta\right)}\nonumber \\
&\leq&\mathrm{Pr}\left(\left\Vert (\alpha_{f,\mathscr{U}} \times \alpha_{h,\mathscr{U}})\left(g\Bigg(\sum\limits_{\bm{j}=\bm{1}}^{\bm{k}}\bm{w}^{\times}_{\bm{j}}\Phi_{\bm{j}}(\bm{A}_{j_1}),\ldots,\sum\limits_{\bm{j}=\bm{1}}^{\bm{k}}\bm{w}^{\times}_{\bm{j}}\Phi_{\bm{j}}(\bm{A}_{j_n})\Bigg)\right)^2\right\Vert_{\ell}\geq \theta\right)
\nonumber \\
&=&\mathrm{Pr}\left(\left\Vert g\Bigg(\sum\limits_{\bm{j}=\bm{1}}^{\bm{k}}\bm{w}^{\times}_{\bm{j}}\Phi_{\bm{j}}(\bm{A}_{j_1}),\ldots,\sum\limits_{\bm{j}=\bm{1}}^{\bm{k}}\bm{w}^{\times}_{\bm{j}}\Phi_{\bm{j}}(\bm{A}_{j_n})\Bigg)\right\Vert_{\ell}\geq \sqrt{\frac{\theta}{(\alpha_{f,\mathscr{U}} \times \alpha_{h,\mathscr{U}})}}\right),
\end{eqnarray}
where the last term of Eq.~\eqref{eq3f h: UP sec: Random Tensor Tail Bounds Algebra} can be upper bounded by those theorems in Section IV in~\cite{chang2022generalMA} with proper assumptions of the function $g$.

\bibliographystyle{IEEETran}
\bibliography{MultivariateMP_NonLinear_Bib}

\begin{thebibliography}{10}
\providecommand{\url}[1]{#1}
\csname url@samestyle\endcsname
\providecommand{\newblock}{\relax}
\providecommand{\bibinfo}[2]{#2}
\providecommand{\BIBentrySTDinterwordspacing}{\spaceskip=0pt\relax}
\providecommand{\BIBentryALTinterwordstretchfactor}{4}
\providecommand{\BIBentryALTinterwordspacing}{\spaceskip=\fontdimen2\font plus
\BIBentryALTinterwordstretchfactor\fontdimen3\font minus
  \fontdimen4\font\relax}
\providecommand{\BIBforeignlanguage}[2]{{%
\expandafter\ifx\csname l@#1\endcsname\relax
\typeout{** WARNING: IEEEtran.bst: No hyphenation pattern has been}%
\typeout{** loaded for the language `#1'. Using the pattern for}%
\typeout{** the default language instead.}%
\else
\language=\csname l@#1\endcsname
\fi
#2}}
\providecommand{\BIBdecl}{\relax}
\BIBdecl

\bibitem{chang2024generalizedCDJ}
S.~Y. Chang and Y.~Wei, ``Generalized choi-davis-jensen's operator inequalities
  and their applications,'' \emph{arXiv preprint arXiv:2403.04892}, 2024.

\bibitem{hutson2005applications}
V.~Hutson, J.~Pym, and M.~Cloud, \emph{Applications of functional analysis and
  operator theory}.\hskip 1em plus 0.5em minus 0.4em\relax Elsevier, 2005.

\bibitem{naylor1982linear}
A.~W. Naylor and G.~R. Sell, \emph{Linear operator theory in engineering and
  science}.\hskip 1em plus 0.5em minus 0.4em\relax Springer Science \& Business
  Media, 1982.

\bibitem{aron2022operator}
R.~M. Aron, M.~S. Moslehian, I.~M. Spitkovsky, and H.~J. Woerdeman,
  \emph{Operator and Norm Inequalities and Related Topics}.\hskip 1em plus
  0.5em minus 0.4em\relax Springer, 2022.

\bibitem{furuta2005mond}
T.~Furuta, J.~Mi{\'c}i{\'c}, J.~Pe{\v{c}}ari{\'c}, and Y.~Seo,
  ``Mond-pe{\v{c}}ari{\'c} method in operator inequalities,'' 2005.

\bibitem{fujii2012recentMP}
M.~Fujii, J.~Mi{\'c}i{\'c}~Hot, J.~Pe{\v{c}}ari{\'c}, and Y.~Seo, \emph{Recent
  developments of Mond-Pe{\v{c}}ari{\'c} method in operator inequalities},
  2012.

\bibitem{chang2024generalizedJensen}
S.-Y. Chang, ``Generalized converses of operator {J}ensens inequalities with
  applications to hypercomplex function approximations and bounds algebra,''
  \emph{arXiv preprint arXiv:2404.11880}, 2024.

\bibitem{chang2024multivariateMP}
------, ``Multivariate {M}ond-{P}ecaric method with applications to
  hypercomplex function {S}obolev embedding,'' \emph{arXiv preprint
  arXiv:2405.17203}, 2024.

\bibitem{cybenko1989approximation}
G.~Cybenko, ``Approximation by superpositions of a sigmoidal function,''
  \emph{Mathematics of control, signals and systems}, vol.~2, no.~4, pp.
  303--314, 1989.

\bibitem{gozlan2019high}
N.~Gozlan, R.~Lata{\l}a, K.~Lounici, and M.~Madiman, \emph{High Dimensional
  Probability VIII: The Oaxaca Volume}.\hskip 1em plus 0.5em minus 0.4em\relax
  Springer Nature, 2019, vol.~74.

\bibitem{chang2022convenient}
S.~Y. Chang and W.-W. Lin, ``Convenient tail bounds for sums of random
  tensors,'' \emph{Taiwanese Journal of Mathematics}, vol.~26, no.~3, pp.
  571--606, 2022.

\bibitem{chang2022randomPDT}
S.~Y. Chang, ``Random parametrization double tensors integrals and their
  applications,'' \emph{arXiv preprint arXiv:2205.03523}, 2022.

\bibitem{chang2022randomMOI}
S.-Y. Chang, ``Random multiple operator integrals,'' \emph{arXiv preprint
  arXiv:2210.09392}, 2022.

\bibitem{chang2022randomDTI}
S.~Y. Chang, ``Random double tensors integrals,'' \emph{arXiv preprint
  arXiv:2204.01927}, 2022.

\bibitem{chang2022generalMaj}
S.~Y. Chang and Y.~Wei, ``General tail bounds for random tensors summation:
  majorization approach,'' \emph{Journal of Computational and Applied
  Mathematics}, vol. 416, p. 114533, 2022.

\bibitem{chang2023tailTRP}
S.-Y. Chang, ``Tail bounds for tensor-valued random process,'' \emph{arXiv
  preprint arXiv:2302.00602}, 2023.

\bibitem{chang2023BiRanTenPartI}
------, ``Random tensor inequalities and tail bounds for bivariate random
  tensor means, part i,'' \emph{arXiv preprint arXiv:2305.03301}, 2023.

\bibitem{chang2023tailMulRanTenMeans}
------, ``Tail bounds for multivariate random tensor means,'' \emph{arXiv
  preprint arXiv:2308.06478}, 2023.

\bibitem{chang2023algebraicConn}
------, ``Algebraic connectivity characterization of ensemble random
  hypergraphs,'' \emph{arXiv preprint arXiv:2310.08700}, 2023.

\bibitem{floudas2008encyclopedia}
C.~A. Floudas and P.~M. Pardalos, \emph{Encyclopedia of optimization}.\hskip
  1em plus 0.5em minus 0.4em\relax Springer Science \& Business Media, 2008.

\bibitem{anastassiou2019frontiers}
G.~A. Anastassiou and J.~M. Rassias, \emph{Frontiers in functional equations
  and analytic inequalities}.\hskip 1em plus 0.5em minus 0.4em\relax Springer,
  2019.

\bibitem{hytonen2016analysis_I}
T.~Hyt{\"o}nen, J.~Van~Neerven, M.~Veraar, and L.~Weis, \emph{Analysis in
  {B}anach spaces Volume I: Martingales and {L}ittlewood-{P}aley Theory}.\hskip
  1em plus 0.5em minus 0.4em\relax Springer, 2016, vol.~63.

\bibitem{hytonen2017analysis_II}
T.~Hyt{\"o}nen, J.~van Neerven, M.~Veraar, and L.~Weis, ``Analysis in banach
  spaces volume ii: Probabilistic methods and operator theory,'' vol.~67, 2017.

\bibitem{dawood2011theories}
H.~Dawood, \emph{Theories of interval arithmetic: mathematical foundations and
  applications}.\hskip 1em plus 0.5em minus 0.4em\relax LAP Lambert Academic
  Publishing, 2011.

\bibitem{chang2022generalMA}
S.~Y. Chang and Y.~Wei, ``General tail bounds for random tensors summation:
  majorization approach,'' \emph{Journal of Computational and Applied
  Mathematics}, vol. 416, p. 114533, 2022.

\end{thebibliography}

\end{document}